\documentclass[11pt]{article}
\usepackage[T1]{fontenc} 
\usepackage[utf8]{inputenc}
\usepackage{lmodern}

\usepackage{amssymb} 
\usepackage{latexsym}
\usepackage{epsfig}
\usepackage{lscape,graphicx,color,psfrag}
\usepackage{calrsfs}
\usepackage[cmex10]{amsmath}
\usepackage{amsthm}

\theoremstyle{plain}
\newtheorem{lemma}{Lemma}
\newtheorem{assumption}{Assumption}
\newtheorem{proposition}{Proposition}
\newtheorem{theorem}{Theorem}
\newtheorem{corollary}{Corollary}
\theoremstyle{definition}
\newtheorem{definition}{Definition}
\theoremstyle{remark}
\newtheorem{remark}{Remark}

\usepackage{appendix}

\usepackage[top=1in, bottom=1in, left=0.9in, right=0.9in]{geometry}

\usepackage[colorlinks=true,urlcolor=blue,citecolor=blue,linkcolor=blue]{hyperref}
\usepackage[noabbrev,nameinlink]{cleveref}     

\def\EMAIL#1{\href{mailto:#1}{#1}}% When hyperref is used, otherwise outcomment 
\begin{document}

\title{Dynamic allocation indices for restless projects and queueing admission control: a
  polyhedral approach\thanks{Work partly supported by 
the Spanish  Ministry of Science and Technology (grant BEC2000-1027),
NATO (Collaborative Linkage Grant PST.CLG.976568),  and
the Joint Spanish-US (Fulbright) Commission for Scientific and Technical
Exchange (project 2000-20132).}}

% The thanks line in the title should be filled in if there is
% any support acknowledgement for the overall work to be included
% This \thanks is also used for the received by date info, but
% authors are not expected to provide this.

\author{Jos\'e Ni\~no-Mora 
\\ Department of Statistics \\
    Carlos III University of Madrid  \\  \EMAIL{jose.nino@uc3m.es}, \href{http://alum.mit.edu/www/jnimora}{http://alum.mit.edu/www/jnimora} \\
      ORCID: \href{http://orcid.org/0000-0002-2172-3983}{0000-0002-2172-3983}}
 
\date{Published in \textit{Mathematical Programming, Series A}, vol.\  93, pp.\ 361--413, 2002 \\ \vspace{.1in}
DOI: \href{https://doi.org/10.1007/s10107-002-0362-6}{10.1007/s10107-002-0362-6}}

\maketitle

% Here is the abstract:

\begin{abstract}%
This paper develops a polyhedral approach to the design, analysis, and 
computation of dynamic allocation indices for scheduling
binary-action (engage/rest) Markovian stochastic projects which can
change state when rested
(\emph{restless bandits (RBs)}), based on  \emph{partial
  conservation laws (PCLs)}.
This extends previous work by the author 
[J. Ni\~no-Mora (2001): Restless bandits, partial
conservation laws and indexability. \textit{Adv. Appl. Probab.} 
\textbf{33}, 76--98], where PCLs were shown to
 imply the optimality of
index policies \emph{with a postulated structure} in 
stochastic scheduling problems, 
under \emph{admissible linear objectives}, and they were deployed
to obtain simple sufficient
conditions for the existence of Whittle's (1988)  RB index (\emph{indexability}),
along with an adaptive-greedy index algorithm.
The new contributions include: 
(i) we develop the polyhedral foundation of the PCL framework, based
on the structural and algorithmic properties of a new
polytope associated with an \emph{accessible set
system} $(J, \mathcal{F})$ (\emph{$\mathcal{F}$-extended
  polymatroid});
(ii) we present new dynamic allocation indices for RBs,
motivated by an admission control model, which
 extend Whittle's and have a significantly
increased scope; 
(iii) we deploy PCLs to obtain both sufficient conditions
for the existence of the new indices (\emph{PCL-indexability}), and 
a new adaptive-greedy index algorithm; 
(iv) we interpret PCL-indexability as a form of the classic
economics law of \emph{diminishing marginal returns}, and 
characterize the index as an \emph{optimal marginal cost rate};
we further solve a related optimal \emph{constrained} control problem;
(v) we carry out a PCL-indexability analysis of the motivating 
admission  control model, under time-discounted and 
long-run average criteria; this gives, under mild conditions, a new index
characterization of optimal threshold policies; and 
(vi) we apply the latter to 
present new heuristic index policies for two hard queueing control 
problems: admission control and routing to parallel queues; and
scheduling a multiclass make-to-stock queue with lost sales, both under
state-dependent holding cost rates and  birth--death dynamics.
\end{abstract}%


\tableofcontents
\newpage

% The body of the paper starts here:
\section{Introduction}
\label{s:intro}
This paper develops a polyhedral approach to the design, analysis, and 
computation of dynamic allocation indices for scheduling
binary-action (engage/rest) Markovian stochastic projects
 which can
change state when rested, or \emph{restless bandits (RBs)}.
The work draws on and contributes to three research areas which have evolved
with substantial autonomy: (1) index policies in stochastic
scheduling; (2) monotone optimal policies in
\emph{Markov decision processes (MDPs)};
and (3) polyhedral methods in resource allocation problems.
We
next briefly discuss each area's relevant background. 

\subsection*{Index policies in stochastic scheduling} 
\emph{Stochastic scheduling} (cf. \cite{nmeoo})
is concerned with the dynamic resource allocation to competing, randomly evolving
 activities.
An important model class concerns the design of a 
\emph{scheduling policy} for optimal dynamic \emph{effort}  allocation 
 to a collection of Markovian
 \emph{stochastic projects}, which can be either
\emph{engaged} or \emph{rested}. 
Following Whittle \cite{whit}, we shall call such
projects \emph{restless bandits (RBs)}, and refer to a
 corresponding multi-project model as a \emph{restless bandit problem (RBP)}.
The term ``project''  is  used here in a lax sense, befitting
 the
application at hand. 
Thus, a project may represent, e.g., a
 \emph{queue} subject to admission control,
whose evolution depends on the \emph{policy} adopted to decide
whether each arriving customer should be admitted into or rejected from
 the system. 

\emph{Index policies} are particularly appealing for such
 problems:
 an index $\nu_k(j_k)$ is
attached to the \emph{states} $j_k$ of each project $k$; then, 
the required number of projects
with larger indices are engaged at each time.
The quest for models with optimal
index policies drew major
research efforts in the 1960s and 1970s, yielding a
classic body of work. This includes the celebrated \emph{$c
  \mu$-rule} \cite{coxsm} for scheduling a multiclass $M/G/1$ queue, 
\emph{Klimov's index rule} \cite{kl} for the corresponding model with
 feedback, 
and the \emph{Gittins index rule} \cite{gijo,gi79} for the
\emph{multiarmed bandit problem (MBP)}.

The MBP is a paradigm among 
such well-solved models, yielding unifying insights. 
In it, rested projects  do not
change state, one project is engaged at
each time, and a discounted criterion is employed.
The optimal \emph{Gittins index} $\nu_k(j_k)$ has an insightful interpretation. 
It was introduced in \cite{gijo} via  a \emph{single-project
subproblem}, where at each time one can continue or abandon
operation, earning in the latter
case a pension at constant rate $\nu$; $\nu_k(j_k)$ 
is then the 
\emph{fair passivity subsidy} in state $j_k$, i.e., 
the minimum value of $\nu$ one
should be willing to accept to rest project $k$. 
Gittins \cite{gi79} further
characterized his index as the maximal rate of expected discounted
reward per unit of expected discounted time, or
\emph{maximal reward rate}, starting at each state.

As for the general RBP, its increased modeling power comes at the expense of
tractability: it is  
\emph{P-SPACE HARD} \cite{patsi}.
The research focus must hence shift to the design of well-grounded,
tractable \emph{heuristic index policies}. 
Whittle \cite{whit} first proposed such a policy, which recovers
Gittins' in the MBP case, and enjoys a form of asymptotic optimality. 
See \cite{wewe90}.
The \emph{Whittle index} is also defined as a \emph{fair passivity
  subsidy} via a single-project
subproblem, precisely as outlined above for the Gittins index.
The \emph{Whittle index policy} 
prescribes to assign higher priority to projects with larger indices.

In contrast to the Gittins index,
passive transitions cause the
 Whittle index to have a \emph{limited scope}: it is  defined only for
an \emph{indexable} project, whose \emph{optimal active set}
 (states where it should be engaged) decreases as the
passive subsidy grows.
The lack of simple sufficient conditions for 
\emph{indexability} hindered the application of such index 
in the 1990s.
An alternative index, free from such scope limitation, was proposed in 
\cite{benirb}. 
Regarding indexability, we first presented in \cite{nmpcl01} a
tractable set of
sufficient conditions, based on the notion of \emph{partial
  conservation laws}, along with a one-pass
index algorithm. 

\subsection*{Monotone optimal policies in MDPs}
In MDP applications intuition often leads to \emph{postulate}
 qualitative properties on
optimal policies. 
The optimal action, e.g., may be \emph{monotone} on the state. 
Thus, in a model for control of admission to a 
queue, one might postulate that arriving customers should be accepted
iff the queue length exceeds a critical \emph{threshold}. 
Establishing the optimality of such policies
 can lead to 
 efficient special algorithms. 

The most developed approach for such purpose is
grounded on the theory of submodular functions on
lattices. See \cite{topk78}.
One must establish \emph{submodularity} properties
on the problem's value function, exploiting the 
\emph{dynamic programming (DP)} equations, by induction on
the \emph{finite} horizon. Infinite-horizon models inherit such properties.
See, e.g.,  \cite[Ch. 8]{heymsob} and \cite{stidh85,altsti}.

Related yet distinct qualitative properties are suggested by the 
\emph{indexability} analysis of RB models,
given in terms of the \emph{monotonicity of the index on 
the state}. 
Consider a queueing admission control RB model, where the active
action corresponds to \emph{shutting} the entry gate, and the passive 
action to \emph{opening} it. 
The Whittle index   would then represent a fair subsidy for
keeping the gate open
\emph{per unit time}; or, equivalently, a fair charge
 for keeping the gate shut per unit time, in each state. 
To be consistent with the optimality of \emph{threshold policies} (see
above), the index should increase monotonically on the queue length. 
Such requirement is critical when the index is used to define a
heuristic policy for related RBPs, such as those discussed in Section \ref{s:appl}.

Yet we have found that, when such model incorporates state-dependent
 arrival rates,
 the Whittle index can
fail to possess the required monotonicity.  See Appendix \ref{a:1}. 
Such considerations motivate us in this paper to develop
extensions of the Whittle index which are consistent with a postulated
structure on optimal policies.

\subsection*{Polyhedral methods in resource allocation
  problems}
The application of polyhedral methods to resource allocation 
originated in \emph{combinatorial optimization}, within
the area of \emph{polyhedral combinatorics}. See, e.g., 
\cite{nemwol}. 
Edmonds \cite{ed70,ed71} first explained the optimality of the
classic \emph{greedy algorithm} ---the simplest 
\emph{index rule} for resource allocation---
from properties of underlying
polyhedra, termed \emph{polymatroids},
arising in the problem's \emph{linear programming (LP)}
formulation. 

The application of LP to MDPs started with the LP formulation of a
general finite-state and -action MDP in
\cite{depenoux,manne}. The seminal application of LP
to stochastic scheduling is due to
Klimov \cite{kl}.
He formulated the problem of optimal
scheduling of a multiclass $M/G/1$ queue with feedback as an LP, whose
constraints represent \emph{flow conservation laws}. He 
solved such LP by an
\emph{adaptive-greedy algorithm}, giving an optimal index rule.

Coffman and Mitrani \cite{cofmi} formulated a simpler
model
 ---\emph{without} feedback--- as an LP,
whose constraints formulate \emph{work conservation laws}. 
These characterize the 
\emph{region of achievable (expected delay) performance} as a
\emph{polymatroid}, thus giving a polyhedral account
for the optimality of the classic 
\emph{$c \mu$  rule}. 
The relation between conservation laws and  polymatroids was
clarified in \cite{fedgro,shaya}. 

Tsoucas \cite{tsouc} applied work conservation laws 
to Klimov's model, obtaining a new LP formulation over an \emph{extended
  polymatroid} (cf. \cite{bhagetso}).
His analysis was extended into the \emph{generalized
  conservation
laws (GCLs)} framework in \cite{beni}, 
giving a polyhedral account of the optimality of Gittins'
index rule
for the MBP and extensions. 
\emph{Approximate GCLs} were  deployed in
\cite{glanm01} to establish the 
near-optimality of Klimov's rule in the
parallel-server case. 
See \cite{arjrss} for an overview
of
such \emph{achievable region approach}.

The theory of conservation laws was extended in \cite{nmpcl01}, through  the notion of 
\emph{partial conservation laws (PCLs)}, 
which were brought to bear on the
analysis of Whittle's  RB index. 
PCLs imply
the optimality of index policies \emph{with a postulated structure}
under \emph{admissible linear objectives}.
Their application  yielded the  class of \emph{PCL-indexable}
RBs,  where the Whittle index exists and is
calculated by an extension of Klimov's algorithm.

\subsection*{Goals, contributions,  and structure}
\label{s:scopep}
The prime goal of this paper is the development, analysis,
 and
application
of
well-grounded extensions of 
Whittle's RB index, which significantly
increase its scope.
For such purpose, we shall deepen the understanding of
 the PCL framework and its polyhedral foundation,
which is the paper's second goal. 

The contributions include: 
(i) we develop the polyhedral foundation of the PCLs, based
on properties of a new
polytope associated with a \emph{set
system} $(J, \mathcal{F})$ (\emph{$\mathcal{F}$-extended
  polymatroid});
(ii) we present new dynamic allocation indices for RBs,
motivated by an admission control model, which
 extend Whittle's and have a significantly
increased scope; 
(iii) we deploy PCLs to obtain both sufficient conditions
for the existence of the new indices (\emph{PCL-indexability}), and 
a new adaptive-greedy index algorithm; 
(iv) we interpret PCL-indexability as a form of the classic
economics law of \emph{diminishing marginal returns}, and 
characterize the index as an \emph{optimal marginal cost rate};
we further solve a related optimal \emph{constrained} control problem;
(v) we carry out a PCL-indexability analysis of the motivating 
admission  control model, under time-discounted and 
long-run average criteria; this gives, under mild conditions, a new index
characterization of optimal threshold policies; and 
(vi) we apply the latter to 
present new heuristic index policies for two hard queueing control 
problems: admission control and routing to parallel queues; and
scheduling a  multiclass make-to-stock queue with lost sales.

The rest of the paper is organized as follows. 
Section \ref{s:tp} describes the motivating admission control
model, and introduces
the new solution approach. 
Section \ref{s:srbosip} describes a general RB model, 
introduces new indices, and formulates the issues to be resolved.
Section \ref{s:egp} introduces $\mathcal{F}$-extended polymatroids,
and studies their properties.
Section \ref{s:papisp} reviews the PCL framework.
Section \ref{s:iprb} applies PCLs  to the analysis of
RBs, 
 yielding sufficient indexability conditions  and an index algorithm.
Section \ref{s:qaca} deploys such results in the 
  admission 
control model. 
Section \ref{s:appl} applies the new indices to present new
policies for two 
complex queueing control models.
Section \ref{s:c} ends the paper with  some concluding remarks.
Three appendices contain important yet ancillary
material. 

\section{Motivating problem: optimal control of admission to a 
birth--death queue}
\label{s:tp}
This section discusses a model for the optimal control
 of admission to a birth--death queue, a fundamental problem
which has drawn extensive research attention. 
 See  \cite{naor,stidh85,altsti,hordspi,chenyao}.
 We shall use the model to motivate our approach, by
introducing  a novel
 analysis grounded on an intuitive index characterization of
optimal \emph{threshold} policies.

\subsection{Model description}
\label{s:md}
\begin{figure}[ht]
\centering
\psfragscanon
\psfrag{B}[dc]{$\mu_{i}$}
\psfrag{D}[dl]{$\lambda_{i}$}
\psfrag{G}{Entry gate}
\includegraphics[height=3in,width=4in,keepaspectratio]{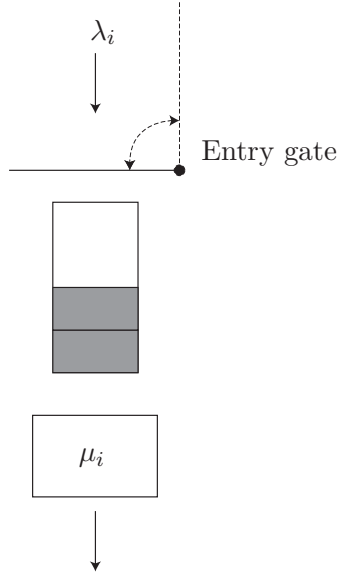}
\caption{Control of admission  to a single queue.}
\label{fig:acrpq2}
\end{figure}

Consider the system portrayed in Figure \ref{fig:acrpq2}, 
which represents a single-server facility catering to an incoming
customer stream, endowed with
a finite buffer capable of holding $n$ customers, 
 waiting or in service.
Customer flow is regulated by a \emph{gatekeeper}, who
dynamically opens or shuts an \emph{entry gate} which customers must
cross to enter the buffer; those finding a shut gate, or a full
buffer, on arrival are rejected and lost.

The \emph{state}
$L(t)$, recording the number in system
at times $t \geq 0$, evolves as a controlled
\emph{birth--death process} 
over \emph{state space} $N = \{0, \ldots, n\}$.
While in state $i$, customers arrive  at
rate
 $\lambda_i$ (being then admitted or rejected),
 and the server works at rate $\mu_i$.

We assume that \emph{holding costs} are continuously incurred in state
$i$ at rate $h_i$, and  a \emph{charge} $\nu$
is incurred \emph{per customer rejection}. Costs are discounted at rate $\alpha > 0$.

The system is  governed by an  \emph{admission control
 policy} $u$, 
prescribing
the action $a(t) \in \{0, 1\}$ to take at each time $t$. 
Policies are chosen from the class
$\mathcal{U}$ of \emph{stationary} policies, basing action choice
 on the  state.
Given policy $u$ and state $j$, 
we denote by $u(j) \in [0, 1]$ the \emph{probability} of taking action 
$a = 1$ (\emph{shut} the entry gate), so that $1-u(j)$ is
 the probability of action $a = 0$ (\emph{open} it). 
We shall refer to $a = 1$ as the \emph{active} action, and to 
$a = 0$ as the \emph{passive} action; 
one may imagine that the gate is naturally open, unless the
 gatekeeper
intervenes to shut it. 
We shall adopt the convention that $u(n) \equiv 1$, so that action
 choice is effectively limited to the set 
$N^{\{0, 1\}} = \{0, \ldots,
 n-1\}$ of 
\emph{controllable states}. 
The single state in  
$N^{\{1\}} = \{n\}$ will be termed \emph{uncontrollable}.

Denote by $E_i^u[\cdot]$ the expectation under policy $u$ when
starting at $i$. Let
\begin{equation}
\label{eq:acmvu}
v^u_i = E_i^u\left[\int_0^\infty h_{L(t)} \, e^{- \alpha t} \,
  dt\right]
\end{equation}
be the corresponding expected total discounted value of
  holding costs incurred, and let
\begin{equation}
\label{eq:tuiacm}
b^u_i = E_i^u\left[\int_0^\infty \lambda_{L(t)} \, a(t) \, e^{- \alpha t} \,
  dt\right]
\end{equation}
 be the expected total discounted number of
customer rejections.
The \emph{cost objective} is 
\[
 v^u_i(\nu) = v_i^u + \nu \, b_i^u.
\]

The \emph{admission control problem} is to find a
 policy
 minimizing such objective:
\begin{equation}
\label{eq:acp}
v_i(\nu) = \min \, \left\{v_i^u(\nu): u \in \mathcal{U}\right\}.
\end{equation}
We shall refer to (\ref{eq:acp}) as the \emph{$\nu$-charge
  problem}. By standard MDP results, there exists
an optimal policy  that is both
  \emph{deterministic}  and 
independent of the initial state $i$. 

Several variations of problem
(\ref{eq:acp}) have drawn extensive research attention,
aiming to establish
the optimality of \emph{threshold} policies (which shut the entry gate
iff $L(t)$ lies above a critical threshold), and to compute and
optimal threshold. See
\cite{naor,stidh85,chenyao,hordspi,altsti}.

\subsection{Optimal index-based threshold policy}
\label{s:mrcisa}
In contrast with previous analyses, we introduce next 
a novel solution approach grounded on the following observation: 
one would expect that, under ``natural'' regularity conditions, 
as 
rejection charge $\nu$ \emph{grows} from $-\infty$ to
$+\infty$, the
subset $S(\nu)$ of \emph{controllable} states where it is optimal to shut the  gate
in  (\ref{eq:acp})  \emph{decreases monotonically}, from 
$N^{\{0, 1\}}$ to  $\emptyset$, dropping states  in the order
$n-1$,  ..., $0$, consistent with 
threshold policies.

In such case,  say that the $\nu$-charge problem  is
\emph{indexable relative to threshold policies}. Then, 
to each state $j \in N^{\{0, 1\}}$ there corresponds a unique
\emph{critical charge} $\nu_j$ under which it is optimal both to 
admit and to reject
a customer arriving in that state. 
Call $\nu_j$ the \emph{dynamic allocation index} of state $j$. Since 
$$
\nu_0 \leq \nu_1 \leq \cdots \leq \nu_{n-1},
$$
such indices yield an optimal \emph{index policy} for the $\nu$-charge problem:
shut the entry gate in state $j \in N^{\{0, 1\}}$
 iff $\nu \leq \nu_j$. 
The \emph{optimal rejection set} is thus 
\begin{equation}
\label{eq:sgamm}
S(\nu) = \left\{j \in N^{\{0, 1\}}: \nu \leq \nu_j\right\},
\quad
\nu \in \mathbb{R}.
\end{equation}

\subsection{Combinatorial optimization formulation}
\label{s:cof}
The $\nu$-charge problem (\ref{eq:acp}) admits a natural \emph{combinatorial
  optimization}
formulation, which will play a key role in
our solution approach.
Represent each stationary 
deterministic policy  by the  subset $S$ of 
controllable states 
 where it takes the
active action, and call it then
the
\emph{$S$-active policy}, writing $b_i^S$, $v_i^S$, $v_i^S(\nu)$.
This gives a reformulation of
$\nu$-charge problem (\ref{eq:acp}) in terms of 
 finding an optimal \emph{active set}:
\[
v_i(\nu) = \min \, \left\{v_i^S(\nu): S \in 2^{N^{\{0, 1\}}}\right\}.
\]

Represent now
the family of threshold policies
by a \emph{set system} $(N^{\{0, 1\}}, \mathcal{F})$, where
$\mathcal{F} \subseteq 2^{N^{\{0, 1\}}}$ is the  nested family of 
\emph{feasible rejection sets} given by 
\begin{equation}
\label{eq:facp} 
\mathcal{F} = \left\{S_1, \ldots, S_{n+1}\right\},
\end{equation}
with $S_{n+1} = \emptyset$ and 
\begin{equation}
\label{eq:skdefac}
S_{k} = \{k-1, \ldots, n-1\}, \quad 1
\leq k \leq n.
\end{equation}

We shall address and solve  the following problems:
\begin{description}
\item[Problem 1:] Give sufficient conditions on model parameters 
under which $\nu$-charge problem (\ref{eq:acp})
  is indexable relative to threshold policies, so that, in particular, 
$$
v_i(\nu) = 
         \min \, \left\{v_i^S(\nu): S \in \mathcal{F}\right\},
         \quad \nu \in \mathbb{R}.
$$
\item[Problem 2:] Give an efficient algorithm for finding an
  optimal threshold policy/optimal active set; or,  equivalently, for constructing 
  the indices $\nu_j$.
\end{description}

\section{RBs: optimality of index policies with
  a postulated structure}
\label{s:srbosip}
This section extends the approach outlined above
to a general \emph{RB} model. 

\subsection{The $\nu$-charge problem for a single RB}
\label{s:mfsrbp} 
Consider  the problem of 
optimal dynamic effort allocation to a single stochastic
project
modeled as an RB, whose state $X(t)$
evolves over discrete time periods $t = 0, 1, \ldots, $ through the finite 
 state space  $N$.
Its evolution is governed by a policy $u$,
prescribing at each period $t$ which of two actions to take: 
\emph{active} (engage the project; $a(t) = 1$) or 
\emph{passive} (let it rest; $a(t) =0$).
Denote by $\mathcal{U}$ the class of state-dependent, or \emph{stationary} policies.
A policy $u \in \mathcal{U}$ is thus a mapping $u: N \to [0,
1]$, where $u(i)$ (resp. $1 - u(i)$) is the probability that
action $a = 1$ (resp. $a = 0$) is 
taken in state $i$.

Taking action $a$ in
state $i$ has two 
effects:
first, cost $h_i^a$ is incurred in the current period, discounted by
 factor $\beta \in (0, 1)$;
second, the next state changes to $j$ with
probability $p_{ij}^a$.
Write $\mathbf{h}^{a}=(h^a_i)_{i \in N}$ and 
$\mathbf{P}^{a}=(p^a_{ij})_{i, j \in N}$.

We shall partition the states as $N = N^{\{0, 1\}} \cup N^{\{1\}}$. 
Here, $N^{\{0, 1\}}$ is the \emph{controllable state space},
where active and  passive actions differ
in some respect; and
$N^{\{1\}} = N \setminus N^{\{0, 1\}}$ is the \emph{uncontrollable
  state space}, where there is no effective choice.
We shall assume that policies $u \in \mathcal{U}$
take the active action at uncontrollable states, i.e., 
\[
u(i) \equiv 1, \quad i \in N^{\{1\}}.
\] 

Let $v_i^u$ be the expected total
discounted value of costs incurred over an infinite horizon under
policy 
$u$, starting at $i$, i.e., 
\begin{equation*}
v_i^u = E_i^{u}\left[ \sum_{t=0}^{\infty
} h^{a(t)}_{X(t)} \,\beta ^{t}\right]. 
\end{equation*}
Besides such \emph{cost measure}, we shall 
consider the \emph{activity measure} 
\begin{equation}
\label{eq:tup}
b_i^u = 
E_i^u\left[\sum_{t=0}^{\infty}
\theta_{X(t)}^1 \, a(t) \,\beta ^{t}\right], 
\end{equation}
where $\boldsymbol{\theta}^1 = (\theta_j^1)_{j \in N} > \mathbf{0}$ is a given 
\emph{activity weight} vector. 
A convenient interpretation results by
considering that the model is obtained via
\emph{uniformization} (cf. Appendix \ref{s:acpdtr}) from a
continuous-time model, as that in Section
\ref{s:tp}. Suppose in the original model there is a 
\emph{distinguished event} (e.g., rejection of an arriving customer),
 which can only occur under the active
action. 
Let $\theta_j^1$ be the probability of the
event happening during a period in state $j$; then,  $b_i^u$ is
the expected total discounted number of times such event occurs. 

Incorporate further into the model an \emph{activity charge}
$\nu$, incurred each time the active action is taken
\emph{and} the distinguished event occurs. Note that $\nu$
corresponds to the \emph{rejection charge} in the previous section. 
The \emph{total cost objective} is then
\[
v_i^u(\nu) = v_i^u + \nu \,
b_i^u.
\]

The \emph{$\nu$-charge problem} of concern
is to find a policy  minimizing such
objective: 
\begin{equation}
v_i(\nu) = 
\min \, \left\{v_i^u(\nu): u \in \mathcal{U}\right\}.  \label{eq:smdcf}
\end{equation}
Again, there exists  an
optimal deterministic policy 
which is
independent of $i$.

\subsection{DP formulation and polynomial-time solvability}
\label{s:dpes} 
The
 conventional approach to tackle $\nu$-charge problem
 (\ref{eq:smdcf}) is
based on formulating and  solving its DP equations,
 which characterize the \emph{optimal value function} $v_i(\nu)$:
\begin{equation}
\label{eq:rbgdp}
v_i(\nu) = 
\begin{cases}
\displaystyle \min_{a \in \{0, 1\}} \, h_i^a + \nu \, \theta_i^1
                      \, a + \beta \sum_{j \in N} p_{ij}^a \, v_j(\nu) 
  & \text{if } i \in N^{\{0, 1\}} \\
\displaystyle h_i^1 + \nu \, \theta_i^1 + \beta \sum_{j \in N}
                      p_{ij}^1 \, v_j(\nu)
  & \text{if } i \in N^{\{1\}}.
\end{cases}
\end{equation}

In theory,
problem (\ref{eq:smdcf}) can be solved in 
\emph{polynomial time} on the state space's size $|N|$.
This follows from (i) the polynomial size of the standard LP
reformulation of (\ref{eq:rbgdp}); 
and (ii) the polynomial-time solvability of LP by 
the \emph{ellipsoid method}.

In practice, however, solution of (\ref{eq:rbgdp}) through
 general-purpose computational techniques can lead to
prohibitively long running times when $|N|$ is large.
Furthermore, even if such solution is obtained, it is not clear
how it could be used to design heuristics for more complex models, 
where RBs arise as building blocks. 

\subsection{Solution by index policies with a postulated structure }
\label{s:espop} 
We next develop an index solution approach to the $\nu$-charge problem, 
motivated by that outlined in Section \ref{s:mrcisa}, which
extends Whittle's  original approach in \cite{whit}.

As  
in Section \ref{s:cof}, $\nu$-charge problem (\ref{eq:smdcf})
admits a combinatorial optimization formulation. 
Associate to every
$S \subset N^{\{0, 1\}}$ a corresponding  \emph{$S$-active policy}, which
is active over states in $S \cup N^{\{1\}}$  and  passive
over $N^{\{0, 1\}} \setminus S$. 
Write $v_i^S$, $b_i^S$ and $v_i^S(\nu)$.
The $\nu$-charge problem is thus reformulated in terms of finding an
optimal \emph{active set}:
\begin{equation}
v_i(\nu) = 
\min \,\left\{v_i^S(\nu): S \in 2^{N^{\{0, 1\}}}\right\}.  \label{eq:copr}
\end{equation}

As before, we shall be concerned with establishing the existence of
optimal policies
within a \emph{postulated} family, 
given by a 
\emph{set system}
$(N^{\{0, 1\}}, \mathcal{F})$.
Here,
$\mathcal{F}\subseteq 2^{N^{\{0, 1\}}}$
is the corresponding family of \emph{feasible active sets}. 
Let $S \subseteq N^{\{0, 1\}}$.

\begin{definition}[$\mathcal{F}$-policy]
\label{def:fprb}
{\rm We say  that the $S$-active policy is an 
\emph{$\mathcal{F}$-policy} if $S \in \mathcal{F}$.}
\end{definition}

Thus,  in the model of Section \ref{s:tp}, the 
$\mathcal{F}$-policies corresponding to the definition of 
$\mathcal{F}$ in (\ref{eq:facp}) 
are precisely the threshold policies.
We shall require set system $(N^{\{0, 1\}},
 \mathcal{F})$ to be \emph{accessible} and \emph{augmentable}. See
Assumption \ref{ass:condss} in Section \ref{s:egp}.

We next
define a key property of the $\nu$-charge problem.
Let $S(\nu) \subseteq N^{\{0, 1\}}$ be, as before, the corresponding set of
controllable states where 
the active action is optimal.

\begin{definition}[Indexability]
\label{def:ti} 
{\rm We say  that the $\nu$-charge problem is 
\emph{indexable relative to $\mathcal{F}$-policies} if, as 
$\nu$ increases from $-\infty $ to $+\infty
 $,  $S(\nu)$ 
 decreases monotonically from 
$N^{\{0, 1\}}$ to $\emptyset$, with $S(\nu) \in \mathcal{F}$ for 
$\nu \in \mathbb{R}$.}
\end{definition}

Under indexability, to each state $j \in N^{\{0, 1\}}$ is attached
a \emph{critical 
charge} $\nu_{j}$, and
\[
S(\nu) = \left\{j \in N^{\{0, 1\}}: \nu \leq \nu_j\right\} \in
\mathcal{F}, \quad \nu \in \mathbb{R}.
\]

\begin{definition}[Dynamic allocation index]
\label{def:ewi}
{\rm We say  that $\nu_j$ is the \emph{dynamic allocation index} of
  controllable state $j \in N^{\{0, 1\}}$ relative to activity measure
  $b^u$.}
\end{definition}

\begin{remark}
Definitions \ref{def:ti} and \ref{def:ewi} extend Whittle's
\cite{whit} 
notion of \emph{indexability} and his index, 
which are recovered in the case 
 $N^{\{0, 1\}} = N$, $\mathcal{F} = 2^N$, 
$\theta_{j}^1 =1$ for $j \in N$. 
\end{remark}

Regarding problems 1 and 2 in Section \ref{s:md}, 
in light of the above  we
shall address and solve  them as special cases of the  following problems: 

\begin{description}
\item[Problem 1:] Give sufficient conditions on model parameters 
under which the $\nu$-charge problem 
  is indexable relative to $\mathcal{F}$-policies, so that, in particular, 
$$
v_i(\nu) = 
         \min \, \left\{v_i^S(\nu): S \in \mathcal{F}\right\},
         \quad \nu \in \mathbb{R}.
$$ 
\item[Problem 2:] Give an efficient algorithm for finding an
  optimal $\mathcal{F}$-policy; or,  equivalently, for constructing
  the indices $\nu_j$.
\end{description}

We shall solve such problems in Section \ref{s:iprb},
by casting them into the polyhedral
framework developed  in Sections 
\ref{s:egp} and \ref{s:papisp} below.

\section{$\mathcal{F}$-extended polymatroids: properties and optimization}
\label{s:egp} 
This section introduces
a new polytope associated with an accessible set system
$(J, \mathcal{F})$, which
generalizes classic polymatroids  and the
extended polymatroids in \cite{bhagetso,beni}.
As we shall see, the problems
of concern in this paper can be formulated and solved as LPs over such polyhedra.
Most proofs  in this section will remain close to those of analogous
results for extended
polymatroids. The exposition will thus focus on the distinctive features
of the new polyhedra. The reader is referred to 
\cite{beni} to fill the details. 

\subsection{$\mathcal{F}$-extended polymatroids}
\label{s:fep}
Let $J$ be a finite ground set with $|J|=n$ elements, and let $\mathcal{F}
\subseteq 2^{J}$ be a
family of subsets of $J$. 
Given a \emph{feasible set} $S \in \mathcal{F}$, let 
$\partial_{\mathcal{F}}^- S$ and  $\partial_{\mathcal{F}}^+ S$
be the \emph{inner and outer boundaries of 
$S$ relative to $\mathcal{F}$}, defined by
\[
\partial_{\mathcal{F}}^- S  = 
\left\{
j \in S: \, S \setminus \{j\}\in \mathcal{F}\right\} \quad
\text{and} \quad
\partial_{\mathcal{F}}^+ S  = 
\left\{
j \in J \setminus S: \, S \cup \{j\}\in \mathcal{F}\right\},
\]
respectively. We shall
 require \emph{set system} $(J, \mathcal{F})$ to satisfy the
 conditions stated next.

\begin{assumption}
\label{ass:condss}
The following conditions hold: \\
{\rm (i)} $\emptyset \in \mathcal{F}$. \newline
{\rm (ii)} Accessibility: $\emptyset \neq S \in \mathcal{F}
\Longrightarrow
\partial_{\mathcal{F}}^- S \neq \emptyset.$ \\
{\rm (iii)} Augmentability: 
$J \neq S \in \mathcal{F} \Longrightarrow
\partial_{\mathcal{F}}^+ S
\neq \emptyset.$
\end{assumption}

We next introduce the notion of \emph{full $\mathcal{F}$-string}.
Let $\boldsymbol{\pi}=(\pi_{1}, \ldots, \pi_{n})$ be an $n$-vector
spanning $J$, so that
$J = \{\pi_1, \ldots, \pi_n\}$. Let
\begin{equation}
\label{eq:spik}
S_k = \{\pi_k, \ldots, \pi_n\}, \quad 
 1 \leq k \leq n.
\end{equation}

\begin{definition}[Full $\mathcal{F}$-string]
\label{def:os} {\rm We say that
$\boldsymbol{\pi}$ is a 
\emph{full $\mathcal{F}$-string} if 
$$S_k \in \mathcal{F}, \quad 1 \leq k \leq n.$$ }
\end{definition}

We shall denote by $\Pi(\mathcal{F})$ the set of all full
$\mathcal{F}$-strings.
Given coefficients $b^S \geq 0$ and 
$w_j^{S} > 0$ for  $j \in S \in \mathcal{F}$, 
consider the \emph{polytope} $P(\mathcal{F})$ on $\mathbb{R}^J$
defined
by 
\begin{equation}
\label{eq:pab}
\begin{split}
\sum_{j \in S} w_j^{S} \, x_j & \geq b^S, \,  S \in \mathcal{F}
 \setminus \{J\} \\
\sum_{j \in J} w^{J}_j \, x_j & = b^J \\
x_j & \geq 0, \quad j \in J.
\end{split}
\end{equation}

For each $\boldsymbol{\pi} \in \Pi(\mathcal{F})$, 
let $\mathbf{x}^{\boldsymbol{\pi}}= (x^{\boldsymbol{\pi}}_j)_{j\in
J}$ be the unique solution
to
\begin{equation}
\label{eq:tsle}
w^{S_k}_{\pi_{k}} \, x_{\pi_{k}} + \cdots
+ w^{S_k}_{\pi_{n}} \, x_{\pi_{n}} =
b^{S_k}, \quad 1 \leq k\leq n,
\end{equation}

\begin{definition}[$\mathcal{F}$-extended polymatroid]
\label{def:egp} 
{\rm We say that $P(\mathcal{F})$ is
an $\mathcal{F}$-\emph{extended polymatroid} if,
for each 
$\boldsymbol{\pi} \in \Pi(\mathcal{F})$, 
$\mathbf{x}^{\boldsymbol{\pi}}~\in~P(\mathcal{F})$.
}
\end{definition}

\begin{remark} \hspace{1in}
\begin{enumerate}
\item Assumption \ref{ass:condss} ensures the existence of a full
  $\mathcal{F}$-string, hence $P(\mathcal{F}) \neq \emptyset$. 
\item
The extended polymatroids in \cite{bhagetso,beni}
correspond to the  case
$\mathcal{F}=2^{J}$. Classic polymatroids are further recovered
when $w^S_j \equiv 1$ for $j \in S \in 2^J$.
\end{enumerate}
\end{remark}

\subsection{LP over $\mathcal{F}$-extended polymatroids}
\label{s:lpoegp} 
Consider the following LP problem over
$\mathcal{F}$-extended polymatroid $P(\mathcal{F})$: 
\begin{equation}
v^{\textrm{LP}} = \min \, 
\left\{ \sum_{j \in J} c_j \, x_j:
\mathbf{x} \in P(\mathcal{F})\right\}.  \label{eq:lpegp}
\end{equation}

We wish to design an efficient algorithm for solving 
LP (\ref{eq:lpegp}), for which we 
start by investigating the 
vertices of $P(\mathcal{F})$. The next result gives
a \emph{partial} characterization, which is in contrast with the
\emph{complete} one available for  extended
polymatroids.

\begin{lemma}
\label{lma:epegp} For $\boldsymbol{\pi} \in \Pi(\mathcal{F})$, 
$\mathbf{x}^{\boldsymbol{\pi}}$ is a vertex of $P(\mathcal{F})$.
\end{lemma}
\begin{proof} 
The result  follows from Definition \ref{def:egp},
along with the
standard algebraic characterization of a polyhedron's vertices.
\end{proof}
\qed

Lemma \ref{lma:epegp} implies that,
under \emph{some} cost vectors $\mathbf{c} = (c_j)_{j \in J}$, 
 LP (\ref{eq:lpegp}) is solved
by a vertex of the form 
$\mathbf{x}^{\boldsymbol{\pi}}$, so that
\begin{equation}
\label{eq:vpilp}
v^{\textrm{LP}} = \min \, \left\{\sum_{j \in J} c_j \, x^{\boldsymbol{\pi}}_j:  \boldsymbol{\pi} \in
\Pi(\mathcal{F})\right\}.
\end{equation}
We shall thus seek to solve  the
LP for a restricted domain of
\emph{admissible} cost vectors, for which an efficient test for 
property (\ref{eq:vpilp}) is available.

To proceed, consider
the dual LP.
By associating dual variable $y^{S}$ with the primal constraint for feasible set
$S \in \mathcal{F}$, the latter is formulated as
\begin{eqnarray}
v^{\text{LP}} &=&\max \,\sum_{S\in \mathcal{F}} b^S \, y^{S}  \label{eq:obdp} \\
&&\text{subject to}  \notag \\
&&\sum_{S: j \in S \in \mathcal{F}} w^{S}_j \, y^{S}\leq c_j, \quad j \in J
\notag \\
&&y^{S}\geq 0, \quad S \in \mathcal{F} \setminus \{J\}  \notag \\
&&y^{J} \text{ unrestricted}.  \notag
\end{eqnarray}
Note that, 
since $P(\mathcal{F})$ is a nonempty polytope, strong
duality ensures that both the
primal and the dual LP 
have the same finite optimal value $v^{\text{LP}}$.

\subsection{Adaptive-greedy algorithm and allocation indices}
\label{s:aga}
This section discusses the \emph{adaptive-greedy
  algorithm} $\mathrm{AG}_1(\cdot |
\mathcal{F})$,  
described in Figure \ref{fig:ag}, which we introduced in \cite{nmpcl01}.
It defines a tractable domain of
\emph{admissible} cost vectors, under which
 it constructs an optimal \emph{index-based} solution to dual LP
 (\ref{eq:obdp}).

The algorithm is fed with input cost vector $\mathbf{c}$, and produces
as output a 
triplet
$(\mathit{ADMISSIBLE}, \boldsymbol{\pi}, \boldsymbol{\nu})$. 
Here, 
$\mathit{ADMISSIBLE} \in \{\mathit{TRUE}, \mathit{FALSE}\}$ is a Boolean variable; 
$\boldsymbol{\pi} = (\pi_1, \ldots, \pi_n) \in 
\Pi(\mathcal{F})$ is
a full
$\mathcal{F}$-string;  and $\boldsymbol{\nu}~=~(\nu_j)_{j \in J}$
is an \emph{index} vector. 
Since it runs in  $n$ steps,
the algorithm will run in \emph{polynomial time} 
if, for  $S \in \mathcal{F}$, 
 calculation of  $w^S_j$ ($j \in S$)
and  membership test $j \in \partial_{\mathcal{F}}^- S$ are
done in polynomial time.

\begin{figure}[tb]
\centering
\fbox{%
\begin{minipage}{\textwidth}
{\bf ALGORITHM} $\mathrm{AG}_1(\cdot | \mathcal{F})$ \\
{\bf Input:} $\mathbf{c}$ \\
{\bf Output:}
$(\mathit{ADMISSIBLE}, \boldsymbol{\pi}, \boldsymbol{\nu})$

\begin{tabbing}
{\it Initialization:} 
{\bf let } $S_1 := J$ \\
{\bf let} $y^{S_1} := \min
 \, \left\{\frac{\displaystyle c_j}{\displaystyle w^{S_1}_j}:
 j \in \partial_{\mathcal{F}}^- S_1\right\};
 \quad$ \\
{\bf choose}  $\pi_1$ {\bf attaining the minimum above}; $\quad$
{\bf let} $\nu_{\pi_1}  := y^{S_1}$
\end{tabbing}
\begin{tabbing}
{\it Loop:} \\
{\bf for} \= $k := 2$ {\bf to} $n$ {\bf do} \\
 \> {\bf let} $S_{k} := S_{k-1} \setminus \{\pi_{k-1}\}$ \\
 \> {\bf let} 
 $y^{S_{k}} := \min \, 
      \left\{\frac{\displaystyle 1}{\displaystyle w^{S_{k}}_j} \, 
 \left[\displaystyle c_j - \sum_{l=1}^{k-1} y^{S_{l}} \,
   w^{S_{l}}_j\right]:
                 j \in \partial_{\mathcal{F}}^- S_{k}\right\}$ \\ 
 \> {\bf choose} $\pi_{k}$ 
{\bf attaining the minimum above}; $\quad$
 {\bf let} $\nu_{\pi_k} := \nu_{\pi_{k-1}} + y^{S_{k}}$ \\
{\bf end} \{for\} 
\end{tabbing}
{\it Cost admissibility test:} \\
{\bf if} $\nu_{\pi_1} \leq \cdots \leq \nu_{\pi_n}$ {\bf then let}
$\mathit{ADMISSIBLE} := \mathit{TRUE}$
{\bf else let} $\mathit{ADMISSIBLE} := \mathit{FALSE}$
\end{minipage}}
\caption{Adaptive-greedy algorithm $\mathrm{AG}_1(\cdot | \mathcal{F})$
  for LP over $\mathcal{F}$-extended polymatroids.}\label{fig:ag}
\end{figure}

The algorithm 
has two new features
relative to its counterpart for extended polymatroids (cf. 
\cite{bhagetso,beni}),
recovered as 
$\mathrm{AG}_1(\cdot | 2^J)$.
First, the minimization in step $k$ is performed over the set
$\partial_{\mathcal{F}}^- S_{k}$, which is often
much smaller than $\partial_{2^J}^- S_{k} = S_{k}$.
Thus, in the model of Section \ref{s:tp}, $S_{k} = 
\{k-1, \ldots, n-1\}$ and $\partial_{\mathcal{F}}^- S_{k} = 
\{k-1\}$.
Second, the algorithm ends with a \emph{cost admissibility test},
checking whether the generated index sequence is nondecreasing. 
We could instead have implemented such test by checking at each
step $k$ whether
 $\nu_{\pi_k} < \nu_{\pi_{k-1}}$, in which case execution would
be terminated.

\begin{definition}[$\mathcal{F}$-admissible costs]
\label{def:facr}
{\rm 
We say that cost vector $\mathbf{c}$ is
\emph{$\mathcal{F}$-admissible for LP (\ref{eq:lpegp})}
if algorithm $\mathrm{AG}_1(\cdot | \mathcal{F})$, when fed with  
input $\mathbf{c}$, returns an output satisfying
\begin{equation}
\label{eq:gleqpk}
\nu_{\pi_1} \leq \nu_{\pi_2} \leq \cdots \leq \nu_{\pi_n}.
\end{equation}
so that $\mathit{ADMISSIBLE} = \mathit{TRUE}$.
}
\end{definition}
We call the set $\mathcal{C}(\mathcal{F})$ of  $\mathcal{F}$-admissible
$\mathbf{c}$'s  the
\emph{$\mathcal{F}$-admissible cost domain of LP
(\ref{eq:lpegp})}.

\begin{remark} \hspace{1in}
\begin{enumerate}
\item
In the extended
polymatroid case, we have $\mathcal{C}(2^J) =
\mathbb{R}^J$.
\item It is readily verified that Definition \ref{def:facr} is
  consistent, i.e., the values of the outputs $\mathit{ADMISSIBLE}$
and $\boldsymbol{\nu}$ 
 do not
depend on the tie-breaking order in the algorithm.
\end{enumerate}
\end{remark}

\begin{definition}[Allocation index]
\label{def:ailpep}
{\rm We say the $\nu_j$'s are LP (\ref{eq:lpegp})'s \emph{allocation indices}. 
 }
\end{definition}

In the extended polymatroid case, such indices 
give an  optimality criterion.
See \cite{beni}. 
We next extend such result. Let 
$\mathbf{c} \in \mathcal{C}(\mathcal{F})$.
Suppose $\mathrm{AG}_1(\cdot | \mathcal{F})$ is run on 
$\mathbf{c}$, giving output
$(\textit{ADMISSIBLE}, \boldsymbol{\pi}, \boldsymbol{\nu})$. 
Let  
$S_k$ be as in (\ref{eq:spik}), and let  
$\mathbf{y}^{\boldsymbol{\pi}} = \left(y^{\boldsymbol{\pi},
    S}\right)_{S \in \mathcal{F}}$ be given by 
\begin{equation}
\label{eq:ypidef}
y^{\boldsymbol{\pi}, S} = \begin{cases} \nu_{\pi_k} - \nu_{\pi_{k-1}} &
  \text{if $S = S_k$, for some $2 \leq k \leq n$} \\
  \nu_{\pi_1} & \text{if $S = S_1$} \\
  0 & \text{otherwise.}
\end{cases}
\end{equation}
Notice 
$y^{\boldsymbol{\pi}, S_k}$, for $1 \leq k \leq n$, is characterized
as
 the unique solution to
\begin{equation}
\label{eq:ypiskdef}
w_{\pi_{k}}^{S_1} \, y^{S_1} +\cdots
+ w_{\pi_{k}}^{S_k} \, y^{S_k}=
c_{\pi_{k}}, \quad 1 \leq k \leq n.
\end{equation}
The next result is proven as its extended polymatroid counterpart (cf. \cite{beni}).
\begin{theorem}[Index-based objective representation and optimality criterion]
\label{the:oc}
\mbox{}
\begin{itemize}
\item[(a)]
LP (\ref{eq:lpegp})'s objective can be represented as
\begin{equation*}
\sum_{j \in J} c_j \, x_j = \nu_{\pi_1} \sum_{j \in S_1} w^{S_1}_j
\, x_j + \sum_{k=2}^{n} (\nu_{\pi_k} - \nu_{\pi_{k-1}}) \,
\sum_{j \in S_k} w^{S_k}_j \, x_j;
\end{equation*}
furthermore, 
\[
v^{\boldsymbol{\pi}} = \sum_{j \in J} c_j \, x^{\boldsymbol{\pi}}_j
 = \nu_{\pi_1} \, b^{S_1} + \sum_{k=2}^{n} (\nu_{\pi_k} -
\nu_{\pi_{k-1}}) \, b^{S_k}. 
\]
\item[(b)]
If condition (\ref{eq:gleqpk}) holds, so that $\mathbf{c} \in
\mathcal{C}(\mathcal{F})$, 
then
$\mathbf{x}^{\boldsymbol{\pi}}$ and $\mathbf{y}^{\boldsymbol{\pi}}$ is
an optimal 
primal-dual pair for
LPs (\ref{eq:lpegp}) and (\ref{eq:obdp}). 
The  optimal value is then
\begin{equation}
\label{eq:ov}
v^{{\rm LP}} = \nu_{\pi_1} \, b^{S_1} + 
 \sum_{k=2}^{n} (\nu_{\pi_k} - \nu_{\pi_{k-1}}) \,
b^{S_k}.
\end{equation}
\end{itemize}
\end{theorem}

\subsection{Allocation index and admissible cost domain decomposition}
\label{s:id} 
The  allocation indices of extended polymatroids possess a useful
 \emph{decomposition} property (cf. \cite{beni}), which we extend next
to $\mathcal{F}$-extended polymatroids. 

Suppose set system $(J, \mathcal{F})$
is constructed as follows. 
We are given $m$ set systems $(J_k, \mathcal{F}_k)$, for $1 \leq k
\leq m$, satisfying
Assumption \ref{ass:condss}, 
where $J_1$, ..., $J_m$ are disjoint.
Let 
$$ J = \bigcup_{k=1}^m J_k, $$
\begin{equation}
\label{eq:fcupsk}
\mathcal{F} = \left\{S = \bigcup_{k=1}^m S_k: S_k \in \mathcal{F}_k,
  1 \leq k \leq m\right\}.
\end{equation}
It is readily verified that set system $(J, \mathcal{F})$ also  satisfies
Assumption \ref{ass:condss}. 

Suppose we are given $b^S \geq 0$ and $w^S_j > 0$,
for $j \in S \in \mathcal{F}$, such that
$P(\mathcal{F})$ defined by (\ref{eq:pab}) is
an $\mathcal{F}$-extended polymatroid.
Then, Definition \ref{def:egp} implies that each
$P_k(\mathcal{F}_k)$ on $\mathbb{R}^{J_k}$
(with $b^{S_k}$ and $w_{j_k}^{S_k}$, for $j_k \in S_k \in \mathcal{F}_k$)
is an $\mathcal{F}_k$-extended polymatroid.

We shall require 
coefficients $w^S_j$ to satisfy the following requirement.
\begin{assumption}
\label{ass:gd} 
For $1 \leq k \leq m$, 
$$w_{j_k}^{S} = w_{j_k}^{S \cap J_{k}}, \quad S \in \mathcal{F}, j_k \in
S \cap J_{k}.$$
\end{assumption}

Given cost vector $\mathbf{c} = (c_j)_{j \in J}$, let 
 $\mathbf{c}^k = (c_{j_k})_{j_k \in J_k}$ for each $k$. 
Consider the
corresponding LPs given by (\ref{eq:lpegp}) and 
\begin{equation}
\label{eq:sublpp}
v^{k, \textrm{LP}} = \min \, 
\left\{ \sum_{j_k \in J_k} c_{j_k} \, x_{j_k}: \mathbf{x}^k \in P_k(\mathcal{F}_k)\right\},
\end{equation}
having admissible cost domains $\mathcal{C}(\mathcal{F})$ and $\mathcal{C}(\mathcal{F}_k)$, respectively.
Let 
$\boldsymbol{\nu} = (\nu_j)_{j \in J}$ (resp. $\boldsymbol{\nu}^k =
(\nu^k_{j_k})_{j_k \in J_k}$)
be the index vector produced
by the algorithm on input $\mathbf{c}$
(resp. $\mathbf{c}_k$).

We state next the decomposition result without proof, as this follows along the same
lines as Theorem 3's in \cite{beni}.

\begin{theorem}[Admissible cost domain and index decomposition]
\label{the:epid}
Under Assumption \ref{ass:gd}, the following holds: \\
\begin{itemize}
\item[(a)] $\mathbf{c} \in \mathcal{C}(\mathcal{F})$ if and only if 
          $\mathbf{c}^k \in \mathcal{C}(\mathcal{F}_k)$ for $1 \leq k
          \leq m$, i.e., 
$$\mathcal{C}(\mathcal{F}) = \prod_{k=1}^m
\mathcal{C}(\mathcal{F}_k). $$
\item[(b)] For $1 \leq k \leq m$, 
$$\nu_{j_k} = \nu^k_{j_k}, \quad j_k \in J_k.$$
\end{itemize}
\end{theorem}

\begin{remark} \hspace{1in}
\begin{enumerate}
\item Theorem \ref{the:epid}(a) shows that  the admissible cost
  domain of LP (\ref{eq:lpegp}) decomposes as the product of the
  corresponding domains of the $m$ LPs in
  (\ref{eq:sublpp}). The 
  $\mathcal{F}$-admissibility test for $\mathbf{c}$ thus
  decomposes into $m$ simpler tasks, which can be performed \emph{in parallel}.
\item Theorem \ref{the:epid}(b) shows that the calculation of 
 indices $\nu_j$ for LP (\ref{eq:lpegp}) can also be
 decomposed into $m$ simpler parallel tasks, each involving the
 calculation of indices $\nu^k_{j_k}$ for the corresponding LP in (\ref{eq:sublpp}).
\end{enumerate}
\end{remark}

\subsection{A new version of the index algorithm}
\label{s:aaia}
We have found that the algorithm above does not lend
itself well to model \emph{analysis}.
This motivates us to develop the reformulated version
$\mathrm{AG}_2(\cdot | \mathcal{F})$, shown in Figure \ref{fig:ag2}.
This represents an extension of Klimov's \cite{kl} algorithm,
recovered as $\mathrm{AG}_2(\cdot | 2^J)$. We shall later apply 
$\mathrm{AG}_2(\cdot | \mathcal{F})$ to calculate the RB indices
introduced in this paper. We remark that Varaiya et al.\ \cite{vawabu}
first 
applied Klimov's algorithm to calculate the Gittins index for
classic (nonrestless) bandits.

\begin{figure}[tb]
\centering
\fbox{%
\begin{minipage}{\textwidth}
{\bf ALGORITHM} $\mathrm{AG}_2(\cdot | \mathcal{F})$: \\
{\bf Input:} $\mathbf{c}$ \\
{\bf Output:}
$(\mathit{ADMISSIBLE}, \boldsymbol{\pi}, \boldsymbol{\nu})$

\begin{tabbing}
{\it Initialization:} 
{\bf let }
$S_1 = J; \quad$ {\bf let} $\nu^{S_1}_j := c_j/w^{S_1}_j, \quad j \in
J$ \\
{\bf choose} $\pi_1 \in \text{argmin}
 \, \left\{\nu^{S_1}_j:
 j \in \partial_{\mathcal{F}}^- S_1\right\};
 \quad$ 
{\bf let}  $\nu_{\pi_1}  := 
\nu^{S_1}_{\pi_1}$
\end{tabbing}
\begin{tabbing}
{\it Loop:} \\
{\bf for} \= $k := 2$ {\bf to} $n$ {\bf do} \\
 \> {\bf let} $S_{k} := S_{k-1} \setminus \{\pi_{k-1}\}$ \\
 \> {\bf let} 
$\displaystyle \nu^{S_{k}}_j :=  \nu^{S_{k-1}}_j + 
\left(\frac{w^{S_{k-1}}_j}{w^{S_k}_j} - 1\right) \, \left[\nu^{S_{k-1}}_j - 
  \nu^{S_{k-1}}_{\pi_{k-1}}\right], \quad j \in S_{k}$ \\
 \> {\bf choose} 
 $\pi_{k} \in \text{argmin} \, 
      \left\{\nu^{S_{k}}_j:
                j \in \partial_{\mathcal{F}}^- S_{k}\right\}; \quad$
 {\bf let} $\nu_{\pi_k} := 
 \nu^{S_{k}}_{\pi_k}$ \\
{\bf end} \{for\} 
\end{tabbing}
{\it Cost admissibility test:} \\
{\bf if} $\nu_{\pi_1} \leq \cdots \leq \nu_{\pi_n}$ {\bf then let}
$\mathit{ADMISSIBLE} := \mathit{TRUE}$
{\bf else let} $\mathit{ADMISSIBLE} := \mathit{FALSE}$
\end{minipage}}
\caption{Adaptive-greedy algorithm $\mathrm{AG}_2(\cdot |
  \mathcal{F})$ for LP over $\mathcal{F}$-extended polymatroids.}
\label{fig:ag2}
\end{figure}

The latter is based on  the
incorporation of coefficients $c^{S_{k}}_j$, recursively defined
(relative to the full $\mathcal{F}$-string
$\boldsymbol{\pi}$ being generated) by 
\begin{equation}
\label{eq:cjsnk}
\begin{split}
c^{S_1}_j & = c_j, \quad j \in S_1 = J \\
c^{S_{k}}_j & = c^{S_{k-1}}_j - 
   \frac{c^{S_{k-1}}_{\pi_{k-1}}}{w^{S_{k-1}}_{\pi_{k-1}}} 
 \, \left[w^{S_{k-1}}_j -
     w^{S_{k}}_j\right], \quad j \in S_{k}, 2 \leq k \leq n,
\end{split}
\end{equation}
which allows to simplify the expressions in $\mathrm{AG}_1(\cdot |
\mathcal{F})$. We shall further write 
\begin{equation}
\nu^{S_k}_j = \frac{c^{S_k}_j}{w^{S_k}_j}, \quad j \in S_k, 1 \leq
k \leq n.
\end{equation}
From (\ref{eq:cjsnk}), it follows that the ratios $\nu^{S_k}_j$ are
characterized by the recursion
\begin{equation}
\begin{split}
\nu^{S_1}_j & = \frac{c_j}{w^{S_1}_j}, \quad j \in S_1 = J \label{eq:nujsnk} \\
\nu^{S_{k}}_j & = \nu^{S_{k-1}}_j + 
\left(\frac{w^{S_{k-1}}_j}{w^{S_k}_j} - 1\right) \, \left[\nu^{S_{k-1}}_j - 
  \nu^{S_{k-1}}_{\pi_{k-1}}\right], \, j \in S_{k}, 2 \leq k \leq n,
\end{split}
\end{equation}

The next result gives the key relations between both  algorithms.
Let  $\boldsymbol{\pi}$ and $\boldsymbol{\nu}$ be
produced by algorithm $\mathrm{AG}_1(\cdot |
  \mathcal{F})$ on input $\mathbf{c}$, 
and let $S_k$ be given by (\ref{eq:spik}). 

\begin{lemma}
\label{lma:identtwoalg}
For $1 \leq k \leq n$ and $j \in S_{k}$, 
\[
\nu^{S_{k}}_j
 = 
\begin{cases}
\frac{\displaystyle c_j}{\displaystyle w^{S_1}_j}, & \text{if $k = 1$} \\
\nu_{\pi_{k-1}}
+ \frac{\displaystyle c_j - \nu_{\pi_1} \, w^{S_1}_j - 
  \sum_{l=2}^{k-1} 
 \left(\nu_{\pi_{l}} - \nu_{\pi_{l-1}}\right) \, 
   w^{S_{l}}_j}{\displaystyle w^{S_{k}}_j}, & \text{if $k
   \geq 2$}; 
\end{cases}
\]
furthermore, 
\begin{equation}
\label{eq:gampinmkcw}
\nu_{\pi_k} =
\nu^{S_{k}}_{\pi_k}.
\end{equation} 
\end{lemma}
\begin{proof}
We proceed by induction on $k$. 
The case $k = 1$ follows from (\ref{eq:cjsnk}). 

Assume now the result holds for $k-1$, where $k \leq n$, so that
\[
\nu^{S_{k-1}}_j = 
\nu_{\pi_{k-2}} + 
\frac{\displaystyle c_j - \nu_{\pi_1} \, w^{S_1}_j - 
  \sum_{l=2}^{k-2} 
 \left(\nu_{\pi_{l}} - \nu_{\pi_{l-1}}\right) \, 
   w^{S_{l}}_j}{\displaystyle w^{S_{k-1}}_j}, \quad 
j \in S_{k-1},
\]
and $\nu_{\pi_{k-1}} =
\nu^{S_{k-1}}_{\pi_{k-1}}$.
Then, applying the
induction hypothesis and 
(\ref{eq:cjsnk}), yields the following: for $j \in S_{k}$, 
\begin{align*}
\nu^{S_{k}}_j & = 
\frac{c^{S_{k-1}}_j - 
   \nu_{\pi_{k-1}} \, \left[w^{S_{k-1}}_j -
     w^{S_{k}}_j\right]}{w^{S_{k}}_j}
\\
& = \frac{\displaystyle \nu_{\pi_{k-2}} \,
  w^{S_{k-1}}_j + c_j - 
  \nu_{\pi_1} \, 
  w^{S_1}_j
  - \sum_{l=2}^{k-2} 
 \left(\nu_{\pi_{l}} - \nu_{\pi_{l-1}}\right) \,
  w^{S_{l}}_j -\nu_{\pi_{k-1}} \, \left[w^{S_{k-1}}_j -
    w^{S_{k}}_j\right]}{w^{S_{k}}_j} \\
& = 
\nu_{\pi_{k-1}}
 + \frac{\displaystyle c_j - \nu_{\pi_1} \, w^{S_1}_j - 
  \sum_{l=2}^{k-1} (\nu_{\pi_{l}} - \nu_{\pi_{l-1}}) \, 
   w^{S_{l}}_j}{\displaystyle w^{S_{k}}_j}.
\end{align*}
Combining the last identity  with
 (\ref{eq:ypidef})--(\ref{eq:ypiskdef}),
gives $\nu_{\pi_k} = \nu^{S_{k}}_{\pi_k}$,
completing the proof.
\end{proof}
\qed

We are now ready to establish the main result of this section.
\begin{theorem}
\label{the:eialg}
Algorithms $\mathrm{AG}_1(\cdot | \mathcal{F})$ and 
$\mathrm{AG}_2(\cdot | \mathcal{F})$ are equivalent.
\end{theorem}
\begin{proof}
The result follows from Lemma \ref{lma:identtwoalg} and
the description of each algorithm. 
\end{proof}
\qed

\subsection{Properties and interpretation 
of coefficients $c^{S_k}_j$, $\nu^{S_k}_j$, and of
  indices $\nu_j$} 
\label{s:icdjskmc}
Given the central role that coefficients $c^{S_k}_j$, $\nu^{S_k}_j$ and indices
$\nu_j$ play in this paper, it is of interest to  discuss their
properties and interpretation. 
Assume below that 
$\boldsymbol{\pi}$, $\boldsymbol{\nu}$ are
produced by  $\mathrm{AG}_2(\cdot |
  \mathcal{F})$ on input $\mathbf{c} \in \mathcal{C}(\mathcal{F})$.

The 
next result shows that 
the $c^{S_k}_j$'s represent \emph{marginal}, or \emph{reduced costs} of 
LP (\ref{eq:lpegp}). 
The proof follows easily by induction, and is hence omitted.

\begin{proposition}
\label{lma:vpicjs}
For $1 \leq m \leq n-1$, 
\begin{equation}
\label{eq:vpiindex}
v^{\rm LP} = 
\sum_{k=1}^m \nu_{\pi_k} \, 
 \left[b^{S_k} - b^{S_{k+1}}\right] + 
  \sum_{j \in S_{m+1}} c^{S_{m+1}}_j \, x^{\boldsymbol{\pi}}_j.
\end{equation}
\end{proposition}

\begin{remark}
Proposition \ref{lma:vpicjs} sheds further light on 
$\mathrm{AG}_2(\cdot | \mathcal{F})$. Identity
(\ref{eq:vpiindex}) shows that, once the first $m$ elements of
optimal $\mathcal{F}$-string $\boldsymbol{\pi} = (\pi_1, \ldots,
\pi_m, \cdot, \ldots, \cdot)$ have been fixed, its construction proceeds 
by optimizing the \emph{reduced objective} 
$\sum_{j \in S_{m+1}} c^{S_{m+1}}_j \, x_j$.  
\end{remark}

We next address the following issue. 
In step $k$ of algorithm $\mathrm{AG}_2(\cdot |
\mathcal{F})$,
the next
element $\pi_k$ is picked through a minimization over $j \in
\partial_{\mathcal{F}}^- S_k$, so that
\[
\nu_{\pi_k} = \min \, \left\{\nu^{S_k}_j: j \in
\partial_{\mathcal{F}}^- S_k\right\}.
\]
Hence, $\nu_{\pi_k}$ is a \emph{locally
 optimal marginal cost
rate} over $j \in \partial_{\mathcal{F}}^- S_k$. 
We shall next show that $\nu_{\pi_k}$ is an
optimal marginal cost rate over the (typically larger) set
$S_k$, i.e., 
\[
\nu_{\pi_k} = \min \, \left\{\nu^{S_k}_j: j \in
S_k\right\}.
\]
We shall need the following preliminary result, easily proven by
induction on $m$.

\begin{lemma}
\label{lma:gpisnback}
For $1 \leq m \leq n$, 
\[
\nu^{S_m}_{\pi_k} = \nu_{\pi_m} + 
  \frac{1}{w^{S_m}_{\pi_k}} \sum_{l=m+1}^k 
   (\nu_{\pi_{l}} - \nu_{\pi_{l-1}}) \, 
  w^{S_{l}}_{\pi_k},
 \quad m \leq k \leq n.
\]
\end{lemma}

We are now ready to establish the index characterization discussed above.
\begin{proposition}
\label{the:gpimin} 
For $1 \leq m \leq n$, the  index 
$\nu_{\pi_m}$ is characterized as
\begin{equation}
\label{eq:1stic}
\nu_{\pi_m} = \min \, \left\{\nu^{S_m}_j: j \in S_m\right\}.
\end{equation}
\end{proposition}
\begin{proof}
By Lemma \ref{lma:gpisnback}, we have, for $m \leq k \leq n$, 
\begin{align*}
\nu^{S_{m}}_{\pi_k} & = \nu_{\pi_m} + 
  \frac{1}{w^{S_m}_{\pi_k}} \sum_{l=m+1}^k 
   (\nu_{\pi_{l}} - \nu_{\pi_{l-1}}) \, 
  w^{S_{l}}_{\pi_k} \geq  \nu_{\pi_m},
\end{align*}
where the inequality follows from the index ordering
 (\ref{eq:gleqpk}).
\end{proof} 
\qed

\subsection{Index characterization under monotone $w^S_j$'s}
\label{s:icimw}
We
have found that, in applications, coefficients $w^S_j$ are 
often
\emph{nondecreasing on $S$}.

\begin{assumption}
\label{ass:imh}
For $j \in S \subset T$, $S, T \in \mathcal{F}$,
\[
w^S_j \leq w^{T}_j.
\]
\end{assumption}

This section shows that Assumption \ref{ass:imh} implies 
interesting additional properties, including a new index
characterization. 
Let $\boldsymbol{\pi}$, $\boldsymbol{\nu}$, $S_k$ be as in
Section \ref{s:icdjskmc}.

\begin{lemma}
\label{lma:ineqind2}
Under Assumption \ref{ass:imh}, the following holds: \\
\begin{itemize}
\item[(a)] For $1 \leq k \leq n-1$, 
\[
\nu^{S_k}_j \leq
\nu^{S_{k+1}}_j, \quad j \in S_{k+1}.
\]
\item[(b)] For $1 \leq k < l \leq n$, 
\[
\nu_{\pi_k} =  \nu^{S_k}_{\pi_k} 
\leq \nu^{S_k}_{\pi_l} \leq 
  \nu^{S_{k+1}}_{\pi_l} 
  \leq \cdots \leq \nu^{S_{l-1}}_{\pi_{l}}
  \leq \nu^{S_l}_{\pi_l} = \nu_{\pi_l}.
\]
\end{itemize}
\end{lemma}
\begin{proof}
(a)
From (\ref{eq:nujsnk}) and (\ref{eq:gampinmkcw}), together with
$w^{S_k}_j - w^{S_{k+1}}_j \geq 0$, it follows that
\[
\nu_{\pi_k} \leq \nu^{S_k}_j \Longrightarrow
\nu^{S_k}_j \leq \nu^{S_{k+1}}_j.
\]
Since the first inequality holds by Proposition \ref{the:gpimin}, we
 obtain the required result. 

(b) The result follows directly from part (a) and Proposition
\ref{the:gpimin}. 
\end{proof}
\qed

We next give the new index characterization referred to above. 
\begin{theorem}
\label{the:ic}
Under Assumption \ref{ass:imh}, the index $\nu_j$ is
characterized as
\[
\nu_j = \max \, \left\{\nu^{S}_j: j \in S \in \{S_1,
  \ldots, S_n\}\right\}, \quad j \in J.
\]
\end{theorem}
\begin{proof}
The result follows directly from Lemma \ref{lma:ineqind2}(b).
\end{proof}
\qed

\begin{remark}
Theorem \ref{the:ic} characterizes the indices as
  \emph{maximal marginal cost rates} relative to
  feasible sets. This is to be contrasted with the
  result in Proposition \ref{the:gpimin}.
\end{remark}

\subsection{A recursion for the $w_j^S$'s under symmetric marginal
  costs}
\label{s:simw}
Recall that the marginal costs $c^{S_k}_j$
were defined \emph{relative to a given
$\boldsymbol{\pi} \in \Pi(\mathcal{F})$}. 
This prevents us from extending (\ref{eq:cjsnk}) into a definition of 
coefficients $c^S_j$, for $S \in \mathcal{F}$, since the \emph{order}
in which $S$ is constructed might lead to different values.
Yet, in certain applications, including RBs
(cf. Section \ref{s:iprb}), such coefficients 
are \emph{symmetric}.

\begin{definition}[Symmetric marginal costs]
\label{ass:symmc}
{\rm We say that marginal costs are \emph{symmetric} if
the following
recursion gives a consistent definition of $c^S_j$, for 
$j \in S \in \mathcal{F}$:
\begin{equation}
\begin{split}
\label{eq:cjsdef}
c^J_j & = c_j, \quad j \in J \\
c^{S \setminus \{i\}}_j & = c^S_j - \frac{c_i^S}{w_i^S} \,
\left[w^S_j - w^{S \setminus \{i\}}_j\right]
\end{split}
\end{equation}}
\end{definition}

Note that, under marginal cost symmetry, we can further define \emph{marginal cost
rates} $\nu^S_j$, for $j \in S \in \mathcal{F}$, by the natural
extension
of recursion (\ref{eq:nujsnk}).
The next result shows that marginal cost symmetry, under
Assumption \ref{ass:imh}, is equivalent to satisfaction of
 a \emph{second-order
  recursion}
 by the
$w^S_j$'s, useful for their calculation.

\begin{proposition}
\label{pro:wjscalc}  
Under Assumption \ref{ass:imh}, marginal costs are symmetric iff, for 
$S \in \mathcal{F}$, $i_1 \in \partial_{\mathcal{F}}^- S \cap 
 \partial_{\mathcal{F}}^- (S \setminus \{i_2\})$, 
$i_2 \in \partial_{\mathcal{F}}^- S \cap 
 \partial_{\mathcal{F}}^- (S \setminus \{i_1\})$, and 
$j \in S \setminus \{i_1, i_2\}$, it holds that
\begin{equation}
\label{eq:wjscalc}
w^{S \setminus \{i_1, i_2\}}_j = 
\frac{\displaystyle \frac{w^S_{i_1}}{w^{S \setminus
      \{i_2\}}_{i_1}} \, w^{S
    \setminus \{i_2\}}_j + 
 \frac{w^S_{i_2}}{w^{S \setminus \{i_1\}}_{i_2}} \, w^{S
    \setminus \{i_1\}}_j - w^S_j}{\displaystyle 
\frac{w^S_{i_1}}{w^{S \setminus \{i_2\}}_{i_1}} + 
  \frac{w^S_{i_2}}{w^{S \setminus \{i_1\}}_{i_2}} - 1}.
\end{equation}
\end{proposition}
\begin{proof}
The result follows by recursively calculating  
$c^{S \setminus \{i_1, i_2\}}_j$ in
  two different ways, using (\ref{eq:cjsdef}): through the
  sequence $S \to S \setminus \{i_1\} \to S \setminus \{i_1, i_2\}$,
  and through the sequence $S \to S \setminus \{i_2\} \to S \setminus
  \{i_1, i_2\}$. Each gives different expressions for $c^{S
  \setminus \{i_1, i_2\}}_j$. Equating the
  coefficients of corresponding  marginal cost terms
 yields the stated identity.
\end{proof}
\qed

\section{Partial conservation laws}

\label{s:papisp} 
This section  reviews the \emph{partial conservation laws (PCLs)}
framework introduced in \cite{nmpcl01}, emphasizing its grounding 
 on
$\mathcal{F}$-extended
polymatroid theory.

Consider a  scheduling
model involving a finite set $J$  of $n$ \emph{job classes}.
Effort is 
allocated to competing jobs through a 
\emph{scheduling policy} $u$, chosen from
  the space $\mathcal{U}$ of \emph{admissible} policies.
Policy $u$'s performance over
class $j$ is given by
\emph{performance measure }
$x^u_j \geq 0$.
Write $\mathbf{x}^u = (x^u_j)_{j \in J}$.
Associate to every \emph{full string} $\boldsymbol{\pi} = (\pi_1, \ldots, \pi_n)$
spanning the $n$ classes a corresponding \emph{$\boldsymbol{\pi}$-priority
  policy},
assigning higher \emph{priority} to class $\pi_{l}$ over
 $\pi_{k}$ if $l>k$.
Write $x^{\boldsymbol{\pi}}_j$.
Given $S \subseteq J$,
say that a policy gives priority to $S$\emph{-jobs} if it gives priority to
any class $i \in S$ over any class 
$j \in S^{c}=J \setminus S$.

We shall be concerned with solving the \emph{scheduling problem}
\begin{equation}
v = \min \, \left\{\sum\limits_{j \in J} c_j \,
  x^u_j: u \in \mathcal{U}\right\},  \label{eq:sppcl}
\end{equation}
which is to find an admissible policy minimizing
 the stated
 linear cost objective. 
Motivated by applications, we shall 
seek to identify conditions under which an optimal policy exists
within a given family of policies with a postulated structure.  
As in Section \ref{s:espop}, we represent the latter by a
\emph{set system} $(J, \mathcal{F})$ satisfying
Assumption \ref{ass:condss}. 
Let $\boldsymbol{\pi}$ be as above.
Recall the notion of \emph{full $\mathcal{F}$-string} from Definition \ref{def:os}. 

\begin{definition}[$\mathcal{F}$-policy]
\label{def:fp} {\rm We say that the $\boldsymbol{\pi}$-priority policy is an 
\emph{$\mathcal{F}$-policy} if  $\boldsymbol{\pi} \in
\Pi(\mathcal{F})$, i.e., 
$\boldsymbol{\pi}$ is a full $\mathcal{F}$-string of set system $(J,
\mathcal{F})$.}
\end{definition}

\begin{remark}
Sets $S \in \mathcal{F}$  represent
\emph{feasible high-priority class subsets} under
$\mathcal{F}$-policies.
\end{remark}

Consider the following problems:

\begin{enumerate}
\item  Give sufficient conditions under which $\mathcal{F}$-policies are optimal, so that 
\begin{equation*}
v = \min \, \left\{\sum\limits_{j \in J} c_j \,
x^{\boldsymbol{\pi}}_j: \boldsymbol{\pi} \in \Pi(\mathcal{F})\right\}.
\end{equation*}
\item  Give an efficient algorithm for finding
 an optimal $\mathcal{F}$-policy.
\end{enumerate}

To address such problems, consider 
the \emph{achievable performance region} 
\begin{equation*}
\mathcal{X}=\left\{ \mathbf{x}^u: u\in \mathcal{U}\right\},
\end{equation*}
which allows us to 
reformulate (\ref{eq:sppcl}) as
the \emph{mathematical programming problem}
\begin{equation*}
v = \min \,\left\{ \sum\limits_{j \in J} c_j \,x_j:
\mathbf{x} \in \mathcal{X}\right\}.  
\end{equation*}
To proceed, we must assume appropriate
properties on
 $\mathcal{X}$, as discussed next.

\subsection{Partial conservation laws}
\label{s:pcllpr} 
Suppose to each job class and feasible high-priority set 
$j \in S \in \mathcal{F}$ is associated a 
coefficient $w^S_j > 0$, so that $\sum_{j \in S } w^{S}_j
\, x^u_j$ represents
a measure of the system's \emph{workload} corresponding to $S$-jobs,
or \emph{$S$-workload}, under policy $u$.  
We shall refer to $w^S_j$ as the 
\emph{marginal $S$-workload of class $j$}. 
Denote the \emph{minimal $S$-workload}  by 
\begin{equation*}
b^S = \inf \,\left\{\sum_{j \in S } w^{S}_j \, x^u_j: u\in 
\mathcal{U}\right\}, \quad S \in \mathcal{F}.
\end{equation*}

\begin{definition}[Partial conservation laws]
\label{def:pcl}
{\rm We say that performance vector $\mathbf{x}^u$ satisfies 
\emph{partial conservation laws (PCLs)}  relative to
 $\mathcal{F}$-policies if 
the following holds:
 \\
(i)
 for $S \in \mathcal{F} \setminus \{J\}$, 
\begin{equation*}
\sum_{j \in S} w^{S}_j \, x^{\boldsymbol{\pi}}_j = b^S,~\quad
 \text{under any $\boldsymbol{\pi} \in \Pi(\mathcal{F})$ giving
 priority to $S$-jobs.}
\end{equation*}
(ii)
$\displaystyle{
\sum_{j \in J} w^{J}_j \, x^{\boldsymbol{\pi}}_j = b^J,~\quad \text{under any }
\boldsymbol{\pi} \in \Pi(\mathcal{F}).}$}
\end{definition}

\begin{remark} \hspace{1in}
\begin{enumerate}
\item Satisfaction of the above PCLs means that, for each
  $S \in \mathcal{F}$,  the $S$-workload
 is 
minimized by \emph{any} $\mathcal{F}$-policy which gives priority to 
$S$-jobs.
\item The \emph{generalized conservation laws (GCLs)} in \cite{beni}
  are recovered in the
case  $\mathcal{F}=2^{N}$. 
The \emph{strong conservation laws} in \cite{shaya} are further
recovered when $w^S_j \equiv 1$.
\end{enumerate}
\end{remark}

Assume in what follows that
$\mathbf{x}^u$
 satisfies PCLs as above.
This gives  a \emph{partial} characterization of achievable
performance region
$\mathcal{X}$, based on
polytope
$P(\mathcal{F})$ in (\ref{eq:pab}).

\begin{theorem}[Achievable performance]
\label{the:arpc} 
$P(\mathcal{F})$ is an 
$\mathcal{F}$-extended polymatroid, satisfying  $\mathcal{X}
\subseteq P(\mathcal{F})$.
The performance vectors $\mathbf{x}^{\boldsymbol{\pi}}$
of 
$\mathcal{F}$-policies $\boldsymbol{\pi}$ are vertices of 
$P(\mathcal{F})$. 
\end{theorem}
\begin{proof} 
PCLs imply $\mathcal{X} \subseteq
P(\mathcal{F})$. 
Let $\boldsymbol{\pi} \in \Pi(\mathcal{F})$. 
By PCL, performance vector
$\mathbf{x}^{\boldsymbol{\pi}}$ is
the solution of 
(\ref{eq:tsle}). Since
$\mathbf{x}^{\boldsymbol{\pi}} \in \mathcal{X} \subseteq
P(\mathcal{F})$, Definition
\ref{def:egp} implies that $P(\mathcal{F})$ is an 
$\mathcal{F}$-extended polymatroid. By Lemma \ref{lma:epegp}, 
$\mathbf{x}^{\boldsymbol{\pi}}$ is a vertex
of $P(\mathcal{F})$. 
\end{proof}
\qed

\begin{remark}
\hspace{1in}
\begin{enumerate}
\item In the GCL case ($\mathcal{F} = 2^J$), it holds that $\mathcal{X} =
P(2^J)$.
See Theorem 4 in \cite{beni}.
\item
By Theorem \ref{the:arpc},  (\ref{eq:lpegp}) is an 
\emph{LP relaxation} of (\ref{eq:sppcl}), hence
 $v^{\mathrm{LP}} \leq v$.
It further implies optimality of 
$\mathcal{F}$-policies for (\ref{eq:sppcl}) under
\emph{some} cost vectors $\mathbf{c}$, so that
$v^{\mathrm{LP}} = v$; and, in particular, under
 \emph{$\mathcal{F}$-admissible cost vectors}
$\mathbf{c} \in \mathcal{C}(\mathcal{F})$ of
 LP (\ref{eq:lpegp}). 
\end{enumerate}
\end{remark}

We show next that, under PCLs,  the scheduling problem
is solved by an index policy with the postulated structure, under appropriate linear 
objectives.
Let $\mathbf{c} \in \mathcal{C}(\mathcal{F})$, and let 
$\boldsymbol{\pi} \in \Pi(\mathcal{F})$ and $\boldsymbol{\nu} =
(\nu_j)_{j \in J}$ be produced by any index
algorithm in Section \ref{s:egp} on input
$\mathbf{c}$. 
Let
$S_k$ be given by (\ref{eq:spik}), and  
let $v^{\mathrm{LP}}$ be the optimal LP value given by
(\ref{eq:lpegp}). 

\begin{theorem}[Optimality of index $\mathcal{F}$-policies]
\label{the:piupcl} 
The $\boldsymbol{\pi}$-priority policy, giving higher priority to
classes with larger indices $\nu_j$, is optimal.
Its  value is
$v = v^{\mathrm{LP}}$.
\end{theorem}
\begin{proof} 
The result follows directly by combining Theorem \ref{the:oc} and
Theorem \ref{the:arpc}.
\end{proof}
\qed

\begin{remark}
Note that, by Theorem \ref{the:oc}, under any policy $u \in \mathcal{U}$ it holds that
$$
\sum_{j \in J} c_j \, x^u_j = 
\nu_{\pi_1} \sum_{j \in S_1} w^{S_1}_j \, x^u_j + 
\sum_{k=2}^{n} (\nu_{\pi_k} - \nu_{\pi_{k-1}}) \,
\sum_{j \in S_k} w^{S_k}_j \, x^u_j.
$$
\end{remark}

\subsection{Multi-project scheduling and index decomposition}
\label{s:idss} 
This section considers the  case where
problem (\ref{eq:sppcl}) represents a 
\emph{multi-project scheduling} model, which represents the
natural setting for application of the
 decomposition property in Section \ref{s:id}.
We shall apply a special case of the result below in  Section \ref{s:iprb}.
Decomposition results have been previously established in
\cite{beni} (under GCLs), and in \cite{nmpcl01} (under PCLs).
The following is a refined version of the latter.

Consider a finite collection of $m \geq 2$ \emph{projects}, with
project $k \in K = \{1, \ldots, m\}$
evolving through finite \emph{state space} $N_{k}$.
Effort is dynamically allocated to projects through a
\emph{scheduling policy} $u \in \mathcal{U}$, where $\mathcal{U}$ is
the space of \emph{admissible policies}, prescribing which of two actions to take at
each project: 
\emph{engage} it ($a_k = 1$) or \emph{rest} it ($a_k = 0$).
Project $k$'s states are partitioned as $N_k = N_k^{\{0, 1\}} \cup
N_k^{\{1\}}$.
When the
state $i_k$ lies in the \emph{controllable state space} $N_k^{\{0, 1\}}$,
the project can be either engaged or rested, whereas it must be engaged
when it lies in the \emph{uncontrollable state space} $N_k^{\{1\}}$.
We  assume that project state spaces are
disjoint. Write $J_k = N_k^{\{0, 1\}}$.

The performance of policy $u$ over state $j_k \in J_k$ of project $k$ is
given by performance measure $x_{j_k}^{k, u} \geq 0$. 
Write $\mathbf{x}^{k, u} = (x_{j_k}^{k, u})_{j_k \in J_k}$. 

The \emph{multi-project scheduling problem} of concern is
\[
v = \min \, \left\{\sum_{k = 1}^m \sum_{j_k \in J_k} c_{j_k}^k \,
 x_{j_k}^{k, u}: u \in \mathcal{U}\right\},
\]
namely, find a policy that minimizes the stated linear performance objective.
This problem fits 
formulation  (\ref{eq:sppcl}),  by 
letting job classes correspond to project states. 

The PCL framework requires a  notion of priority among 
classes. 
In the current setting, this follows from the natural
notion of priority among projects.
We thus interpret each full string $\boldsymbol{\pi} =
(\pi_1, \ldots, \pi_n)$, where $n = |J|$, as a corresponding
$\boldsymbol{\pi}$-priority policy. 

We must further specify a
set system $(J, \mathcal{F})$, 
defining the family of $\mathcal{F}$-policies. 
Assume we are given a family of policies for operating each project 
$k$ \emph{in isolation}, i.e., prescribing in which controllable
states it should be
engaged, given as an appropriate set system
$(J_k, \mathcal{F}_k)$. 
Project $k$'s \emph{$\mathcal{F}_k$-policies} are obtained by associating
to each set $S_k \in \mathcal{F}_k$  a corresponding 
\emph{$S_k$-active policy}, which engages the project when its state
lies in $S_k \cup N^{\{1\}}$, and rests it otherwise. 
Construct now  $(J, \mathcal{F})$ as in (\ref{eq:fcupsk}). 

Assume further that (i) performance vector
$\mathbf{x}^u = (x^u_j)_{j \in J}$ satisfies PCLs relative to $\mathcal{F}$-policies; and 
that (ii) 
marginal workloads $w^S_j$ satisfy Assumption \ref{ass:gd}.

Suppose every project $k$'s cost vector $\mathbf{c}^k = (c_{j_k}^k)_{j_k \in J_k}$
is $\mathcal{F}_k$-admissible.
Then, Theorem
\ref{the:epid}(a) gives that 
$\mathbf{c} = (c_j)_{j \in J}$, where $c_{j_k} = c_{j_k}^k$,  is 
$\mathcal{F}$-admissible. 
Let $\boldsymbol{\nu}^k = (\nu_{j_k}^k)_{j_k \in J_k}$  be project $k$'s
index vector. 
The following result follows from  Theorem \ref{the:epid}. 

\begin{theorem}[Index decomposition for multi-project scheduling]
\label{the:idps} 
Any $\mathcal{F}$-policy $\boldsymbol{\pi}$ giving higher priority
to projects $k$ whose states $j_k$ have larger indices $\nu_{j_k}^k$
is optimal.
\end{theorem}

\section{PCL-indexable RBs}
\label{s:iprb} 
In this section we return to the RB model discussed in
Section \ref{s:srbosip}. We shall resolve the issues raised in Section 
\ref{s:espop} by deploying the PCL framework. 

\subsection{Standard LP formulation and pure passive-cost
  normalization}
\label{s:slpfppcn}
We review next the standard LP formulation 
of \emph{$\nu$-charge problem} (\ref{eq:smdcf}), 
 arising as
the dual of the LP formulation of its DP equations (\ref{eq:rbgdp}).
We shall use it to reduce the problem
to a
\emph{pure passive-cost} normalized version, on which we
shall focus our analyses.

The standard LP formulation of $\nu$-charge problem (\ref{eq:smdcf})
is 
\begin{eqnarray}
v_i(\nu)
&=& \min \, \mathbf{x}^{0} \, \mathbf{h}^{0} +
\mathbf{x}^{1}(\mathbf{h}^{1} + \nu \,
\boldsymbol{\theta}^1)  \label{eq:dclpf} \\
&&\text{subject to}  \notag \\
&&\mathbf{x}^{0}\,(\mathbf{I}-\beta \,\mathbf{P}^{0})+\mathbf{x}^{1}\,
(\mathbf{I}-\beta \,\mathbf{P}^{1}) = \mathbf{e}_i  \label{eq:slpc} \\
&& x^0_j = 0, \quad j \in N^{\{1\}}, \label{eq:xj00} \\
&&\mathbf{x}^{0}, \mathbf{x}^{1}\geq \mathbf{0} \notag
\end{eqnarray}
where $\mathbf{x}^{a}=(x^{a}_j)_{j \in N}$, 
$\boldsymbol{\theta}^1 = (\theta_j^1)_{j \in N}$, and $\mathbf{e}_i$ is the $i$th unit coordinate vector in $\mathbb{R}^N$.
 Vectors are in row or
column form as required. In such LP,
variable $x^a_j$ corresponds
 to the standard \emph{state-action occupation measure}
\begin{equation*}
x_{ij}^{a, u} = E_i^u\left[\sum_{t = 0}^\infty 1\{X(t) = j, a(t) = a\}
  \, \beta^t\right],
\end{equation*}
giving the expected total discounted number of times
action $a$ is taken in state $j$ under policy $u$, starting at $i$.
Thus, (\ref{eq:xj00}) says  the
project must be active at uncontrollable states.

The LP constraints (\ref{eq:slpc}) imply that
\begin{equation*}
\mathbf{x}^{1}= \mathbf{e}_i \,
(\mathbf{I}-\beta \,\mathbf{P}^{1})^{-1} -
 \mathbf{x}^{0} \, (\mathbf{I}- \beta \,\mathbf{P}^{0}) \, 
(\mathbf{I} - \beta \,\mathbf{P}^{1})^{-1},
\end{equation*}
and hence its objective (\ref{eq:dclpf}) can be
reformulated  as 
\begin{equation}
\label{eq:objrefptc}
\mathbf{e}_i \,
(\mathbf{I}-\beta \,\mathbf{P}^{1})^{-1} \, \mathbf{h}^{1} + 
\mathbf{x}^0 \, \widehat{\mathbf{h}}^0 + \nu \, \mathbf{x}^1 \,
\boldsymbol{\theta}^1 = 
v_i^{N^{\{0, 1\}}} + 
\mathbf{x}^0 \, \widehat{\mathbf{h}}^0 + \nu \, \mathbf{x}^1 \,
\boldsymbol{\theta}^1, 
\end{equation}
where $\widehat{\mathbf{h}}^0 = \left(\widehat{h}^0_j\right)_{j \in N}$
 is the \emph{normalized passive-cost vector} given
by 
\begin{equation}
\label{eq:hatc0}
\widehat{\mathbf{h}}^0 = 
\mathbf{h}^{0} - (\mathbf{I} - \beta \, \mathbf{P}^{0}) \,
(\mathbf{I} - \beta \,\mathbf{P}^{1})^{-1} \, \mathbf{h}^{1}.
\end{equation}
Note further that identity (\ref{eq:hatc0}) and the definition of
uncontrollable states gives
\begin{equation*}
\widehat{h}^0_j = 0, \quad j \in N^{\{1\}}.
\end{equation*}

We shall focus henceforth 
on the following \emph{normalized
$\nu$-charge problem}
\begin{equation}
\label{eq:normgcp}
\widehat{v}_i(\nu)  = \min \, \left\{\sum_{j \in N^{\{0, 1\}}} \widehat{h}_j^0
\, x_{ij}^{0, u} + \nu \, b_i^u: u \in \mathcal{U}\right\},
\end{equation}
whose optimal value is related to $v_i(\nu)$ by 
\begin{equation*}  
v_i(\nu)  = 
v_i^{N^{\{0, 1\}}} + 
\widehat{v}_i(\nu).
\end{equation*}

\subsection{PCLs for normalized $\nu$-charge problem}
\label{s:rtpsp}
We shall next cast 
  problem (\ref{eq:normgcp}) into the
 multi-project scheduling case of the PCLs in Section
 \ref{s:idss}.
Reinterpret (\ref{eq:normgcp}) as a 
\emph{two-project} scheduling model, by
adding to the original project  a \emph{calibrating
project} with a \emph{single state} $\ast $.
One project must be 
engaged at each time, where the calibrating project is
 engaged when the original project is
 rested. 

As in Section \ref{s:idss}, let
the controllable state space of the two-project
model be
$$J^* = N^{\{0, 1\}} \cup \{\ast \}.$$
We shall seek to establish PCLs for performance vector
 $\mathbf{x}_i^u = (x_{ij}^u)_{j \in J^*}$, where 
\begin{equation}
\label{eq:xuijdef}
x_{ij}^u = 
\begin{cases}
x_{ij}^{0, u} & \text{if $j \in N^{\{0, 1\}}$} \\ 
  b_i^u & 
\text{if $j =\ast$.}
\end{cases}
\end{equation}
Note that the normalized $\nu$-charge problem
can then be formulated as
\begin{equation*}
\widehat{v}_i(\nu) =  \min \, 
\left\{\sum_{j \in N^{\{0, 1\}}} \widehat{h}_j^0 \, x_{ij}^{u} 
+ \nu \, x_{i*}^u: u \in \mathcal{U}\right\}.
\end{equation*}
Regarding \emph{priorities}'s interpretation, note  that, e.g., 
 giving higher priority to calibrating project's state $*$ over
 original project's state
$j$ means that the latter is
rested in state $j$.

Recall from Section \ref{s:espop} that we are given 
an appropriate set system $(N^{\{0, 1\}}, \mathcal{F})$ 
defining the family
of \emph{$\mathcal{F}$-policies} (cf. Definition \ref{def:fprb}).
Proceeding as in Section \ref{s:idss}, construct a 
set system $(J^*, \mathcal{F}^*)$ for the two-project
model by letting
$$
\mathcal{F}^* = \left\{S^* = S_1 \cup S_2: S_1 \in \mathcal{F}, S_2 \in
\{\emptyset, \{*\}\}\right\}.
$$

We shall seek to establish that performance vector
$\mathbf{x}_i^u$ satisfies PCLs relative to $\mathcal{F}^*$,
 for which suitable coefficients
$w^{S^*}_j$ and $b^{S^*}$ must be defined. 

We start by defining \emph{marginal workloads}
$w^S_j$, for $j \in N$, $S  \subseteq N^{\{0, 1\}}$, in terms of
\emph{activity measures} $b_i^S$ (cf. 
Section \ref{s:mfsrbp}). The latter are characterized by
\begin{equation}
\label{eq:tiscalc}
b_i^{S}  =  
\begin{cases}
\displaystyle
\theta_i^1 + \beta \, \sum_{j \in N} p_{ij}^1 \, b_j^S & \text{if } i
\in S \cup N^{\{1\}}  \\
\displaystyle
\beta \, \sum_{j \in N} p_{ij}^0 \, b^S_j, & \text{if } i
\in N^{\{0, 1\}} \setminus S;
\end{cases}
\end{equation}

We shall use below the following
notation: 
given $\mathbf{d} = (d_j)_{j \in N}$, 
$\mathbf{A} = (a_{i, j})_{i, j \in N}$, and 
$S, T \subseteq N$, we shall write $\mathbf{d}_S = (d_j)_{j \in S}$ and 
$\mathbf{A}_{S T} = (a_{i j})_{i \in S, j \in T}$.
We can thus reformulate the above equations as
\begin{equation*}
\begin{split}
\mathbf{b}_{S \cup N^{\{0, 1\}}}^S & =  \boldsymbol{\theta}_{S \cup N^{\{0, 1\}}}^1 + \beta \, 
  \mathbf{P}_{S \cup N^{\{0, 1\}}, N}^1 \, \mathbf{b}^S \\
\mathbf{b}_{N^{\{0, 1\}} \setminus S}^S & =  
\beta \, \mathbf{P}_{N^{\{0, 1\}} \setminus S, N}^0 \,  \mathbf{b}^S.
\end{split}
\end{equation*}
Let now 
\begin{equation}
\label{eq:wisdef}
w_i^S = \theta_i^1 \, 1\{i \in N^{\{0, 1\}}\}
  + \beta \, \sum_{j \in N} (p_{ij}^1 - p_{ij}^0) \, 
  t^S_j, \quad i \in N; 
\end{equation}
i.e.,
\begin{equation}
\begin{split}
\label{eq:vecwisdef}
\mathbf{w}_{N^{\{0, 1\}}}^S & = \boldsymbol{\theta}_{N^{\{0, 1\}}}^1 + \beta \, 
   \left(\mathbf{P}_{N^{\{0, 1\}}, N}^1 - \mathbf{P}_{N^{\{0, 1\}},
   N}^0\right)
 \, \mathbf{b}^S \\
\mathbf{w}_{N^{\{1\}}}^S & = \beta \, 
   \left(\mathbf{P}_{N^{\{1\}}, N}^1 - \mathbf{P}_{N^{\{1\}},
   N}^0\right) \, \mathbf{b}^S = \mathbf{0},
\end{split}
\end{equation}
where the last identity follows from the assumption
 $p_{ij}^1 = p_{ij}^0$ for $i \in N^{\{1\}}$.
Coefficient $w_i^S$ thus represents the
\emph{marginal increment in activity measure $b^S$ resulting from
a passive-to-active action interchange
 in initial state
$i$.} 

We proceed with a preliminary result, 
giving further relations between 
$\mathbf{b}^S$ and $\mathbf{w}^S$. The proof is omitted, as it follows by straightforward
algebra from the above.
\begin{lemma}
\label{lma:reltws}
The following identities hold: 
\begin{equation}
\label{eq:imbp0bp1}
\begin{split}
(\mathbf{I} - \beta \, \mathbf{P}^0) \, \mathbf{b}^S & = 
\left[\begin{array}{c} \mathbf{w}_S^S \\ \mathbf{0}_{N^{\{0, 1\}} \setminus S}
       \\
 \boldsymbol{\theta}_{N^{\{1\}}}^1 \end{array}\right]
   \\
\boldsymbol{\theta}^1 - 
(\mathbf{I} - \beta \, \mathbf{P}^1) \, \mathbf{b}^S & = 
\left[\begin{array}{c} \mathbf{0}_{S \cup N^{\{1\}}} \\
    \mathbf{w}_{N^{\{0, 1\}} \setminus
  S}^S\end{array}\right].
\end{split}
\end{equation}
\end{lemma}

Motivated by Assumption \ref{ass:gd}, we 
complete the marginal workload definitions by letting,  for 
$j \in J^*$ and $S^* = S \cup \{\ast \}$,
 with $S \subseteq N^{\{0, 1\}}$, 
\begin{equation*}
w^{S^*}_j = 
\begin{cases}
w^{S}_j & \text{if $j \in N^{\{0, 1\}}$} \\ 
1 & \text{if $j = *$.}
\end{cases}
\end{equation*}

It remains to define the 
function $b_i^{S^*}$  arising in the right-hand side of the
PCLs (which now depends on initial state 
$i$). 
Let, for $S^* \subseteq J^*$, 
\begin{equation*}
b_i^{S^{\ast}}= 
\begin{cases}
b_i^S & \text{if $S^{\ast }=S \cup \{\ast \}, 
\emptyset \neq S \subseteq N^{\{0, 1\}}$} \\ 
0 & \text{otherwise}.
\end{cases}
\end{equation*}

The next result gives a set
of \emph{workload decomposition laws}, i.e., linear
equations 
 relating 
workload terms corresponding to the active and the passive action.

\begin{proposition}[Workload decomposition laws]
\label{lem:keys} 
For $u \in \mathcal{U}$ and $S \subseteq N^{\{0, 1\}}$, 
\begin{equation*} 
b_i^u + 
\sum_{j \in S} w^{S}_j \, x_{ij}^{0, u} 
=
b_i^S + 
\sum_{j \in N^{\{0, 1\}} \setminus S} w_j^{S} \, x_{ij}^{1, u}.
\end{equation*}
\end{proposition}
\begin{proof}
Using in turn 
equations (\ref{eq:slpc}) and (\ref{eq:imbp0bp1}), we have
\begin{equation*}
\begin{split}
0 & = \left[
\mathbf{x}_i^{0, u} \,  (\mathbf{I} - \beta \, \mathbf{P}^0) + 
\mathbf{x}_i^{1, u} \,  (\mathbf{I} - \beta \,
\mathbf{P}^1) - \mathbf{e}_i\right] \, 
\mathbf{b}^S \\
& =  
\mathbf{x}_i^{0, u} \, (\mathbf{I} - \beta \, \mathbf{P}^0)
\, \mathbf{b}^S + 
\mathbf{x}_i^{1, u} \,
\left[(\mathbf{I} - \beta \, \mathbf{P}^1) \, \mathbf{b}^S
 - \boldsymbol{\theta}^1\right] - 
\mathbf{e}_i \, \mathbf{b}^S + 
\mathbf{x}_i^{1, u} \, \boldsymbol{\theta}^1 \\
& = 
\mathbf{x}_{i, S}^{0, u} \, \mathbf{w}_S^S -
\mathbf{x}_{i, N \setminus S}^{1, u} \, \mathbf{w}_{N
  \setminus S}^S - 
b_i^S + 
b_i^u,
\end{split}
\end{equation*}
which gives the required result, after simplification
using  Lemma \ref{lma:reltws}.
\end{proof}
\qed

The relation 
between coefficients $b^S_j$'s and $w^S_j$'s is further clarified
next.

\begin{corollary}
\label{cor:tsjwsj}
For $i \in N$ and $S \subseteq N^{\{0,
  1\}}$, 
\begin{align*}
b_i^{S \cup \{j\}} & = b_i^S + w^S_j \, x_{ij}^{1, S \cup \{j\}},
\quad j \in N^{\{0, 1\}} \setminus S \\
b_i^S & = b_i^{S \setminus \{j\}} + w^S_j \, x_{ij}^{0, S \setminus \{j\}},
\quad j \in S.
\end{align*}
\end{corollary}
\begin{proof}
It follows by letting $u = S \cup \{j\}$ and 
$u = S \setminus \{j\}$ in Proposition \ref{lem:keys}, respectively.
\end{proof}
\qed

Proposition \ref{lem:keys} suggests the following conditions
for satisfaction of PCLs. 

\begin{assumption}
\label{ass:pcl1}
Marginal workloads $w^S_j$ satisfy the following: for $S \in \mathcal{F}$, 
$$
w^S_j > 0, \quad j \in
N^{\{0, 1\}}.
$$
\end{assumption}

Assumption
\ref{ass:pcl1} represents a monotonicity property of 
$b^u$, as shown next.

\begin{proposition}
\label{pro:iuwjsg0} 
Assumption \ref{ass:pcl1} is equivalent to the following: 
for $S \in \mathcal{F}$,
\begin{equation}
\begin{split}
\label{eq:tjstjscpj}
b^S_j & < b^{S \cup \{j\}}_j, \quad j \in N^{\{0, 1\}} \setminus S \\
b^S_j & > b^{S \setminus \{j\}}_j, \quad j \in S.
\end{split}
\end{equation}
\end{proposition}
\begin{proof}
The result follows from Corollary \ref{cor:tsjwsj}, by noting 
that 
$x_{jj}^{1, S \cup \{j\}} > 0$, for $j \in N^{\{0, 1\}} \setminus S$,
and 
$x_{jj}^{0, S \setminus \{j\}} > 0$, for $j \in S$.
\end{proof}
\qed

We are now ready to established the required PCLs.
\begin{theorem}[PCLs]
\label{the:pit} Under Assumption \ref{ass:pcl1}, 
performance vector 
$\mathbf{x}_i^u$ satisfies PCLs relative to $\mathcal{F}^*$-policies. 
\end{theorem}
\begin{proof} 
The result follows by combining Proposition 
\ref{lem:keys} with  Assumption \ref{ass:pcl1}.
Consider, e.g., the case $S^* = S \cup \{*\}$, where $\emptyset \neq S \in \mathcal{F}$. 
Under any policy $u \in \mathcal{U}$, 
\begin{align*}
\sum_{j \in \mathcal{F}^*} w^{S^*}_j \, x_{ij}^u & = 
b_i^u +   \sum_{j \in S} w^S_j \, x_{ij}^{0, u} \\
& = b_i^S + 
\sum_{j \in N^{\{0, 1\}} \setminus S} w^{S}_j \, x_{ij}^{1, u} \\
& \geq b_i^S = b_i^{S^*},
\end{align*}
with equality attained in the last inequality if 
\emph{priority} is given to \emph{$S^*$-jobs}, i.e., if the passive action is
taken at states $j \in N^{\{0, 1\}} \setminus S$. 
Other cases follow similarly.
\end{proof}
\qed

We next define a class of RBs that will be
shown to be indexable. Let $n = |N^{\{0, 1\}}|$.

\begin{definition}[PCL-indexable RBs]
\label{def:pclirb}
{\rm We say the RB is
  \emph{PCL-indexable} relative to activity measure
$b^u$ and $\mathcal{F}$-policies if the following conditions holds: \\
(i) \emph{Positive marginal workloads:} Assumption \ref{ass:pcl1} holds. \\
(ii) \emph{Index monotonicity:} Let 
$(\mathit{ADMISSIBLE}, \boldsymbol{\pi}, \boldsymbol{\nu})$ be the
  output of any index algorithm in Section \ref{s:egp} on input
  $\widehat{\mathbf{h}}_{N^{\{0, 1\}}}^0$. Then, the indices satisfy
\begin{equation}
\label{eq:indxmon}
\nu_{\pi_1} \leq \cdots \leq \nu_{\pi_n},
\end{equation}
i.e., $\mathit{ADMISSIBLE} = \mathit{TRUE}$, or 
$\widehat{\mathbf{h}}_{N^{\{0, 1\}}}^0 \in \mathcal{C}(\mathcal{F})$}.
\end{definition}

\begin{remark}
The
definition of PCL-indexability  in
\cite{nmpcl01} is recovered
in the case
$\theta_j^1 \equiv 1$.
\end{remark}

Assume below that the RB
is PCL-indexable. 
Feed any index algorithm
 with input
$\widehat{\mathbf{h}}_{N^{\{0, 1\}}}^0 $
to get $\mathcal{F}$-string 
$\boldsymbol{\pi}$ 
and index vector $\boldsymbol{\nu}$. 
Let $S_k = \{\pi_k, \ldots,
\pi_n\}$, for $1 \leq k \leq n$. 
The next result shows that PCL-indexability implies
indexability
(cf. Definition \ref{def:ti}).

\begin{theorem}[PCL-indexability $\Longrightarrow$ indexability]
\label{the:indexrb}
The RB is
  indexable, and the dynamic allocation index of state 
$j$ is $\nu_j$, for $j \in N^{\{0, 1\}}$.
\end{theorem}
\begin{proof}
Theorem \ref{the:idps} applies to the two-project formulation
of the normalized
$\nu$-charge problem.
It follows that (i) the priority index of the calibrating project's 
state is $\nu_* = \nu$; and (ii) 
the dynamic allocation index for the original project's controllable state $j$
is $\nu_j$. 
The result now follows by interpreting Theorem \ref{the:idps} in terms of
Definition \ref{def:ti}.
\end{proof}
\qed

Several consequences follow  from the above, 
starting with a reformulation of the $\nu$-charge problem
as an LP 
over an $\mathcal{F}^*$-extended polymatroid. 
Consider the polyhedron $P_i(\mathcal{F}^*) \subset
\mathbb{R}^{J^*}$ 
defined as in
(\ref{eq:pab}), relative to parameters $w^{S^*}_j$ and $b_i^{S^*}$
as above.

\begin{corollary}
\label{cor:lprgcp} 
$P_i(\mathcal{F}^*)$ is an $\mathcal{F}^*$-extended
polymatroid.
The 
$\nu$-charge problem can be reformulated as the LP
$$
v_i(\nu) = v_i^{N^{\{0, 1\}}} + 
  \min \, \left\{\sum_{j \in N^{\{0, 1\}}} \widehat{h}^0_j \, x_j + \nu
    \, x_*: \mathbf{x}  \in
    P_i(\mathcal{F}^*)\right\}.
$$
\end{corollary}

\begin{figure}[ht]
\centering
\psfragscanon
\psfrag{K}[u]{$ $}
\psfrag{I}[r]{$v_i(\nu)$}
\psfrag{B}[ul]{$\nu_{\pi_n}$}
\psfrag{H}[c]{$\nu_{\pi_1}$}
\psfrag{D}[u]{$\nu_{\pi_2}$}
\psfrag{F}[u]{$\nu$}
\includegraphics[height=3in,width=4.8in,keepaspectratio]{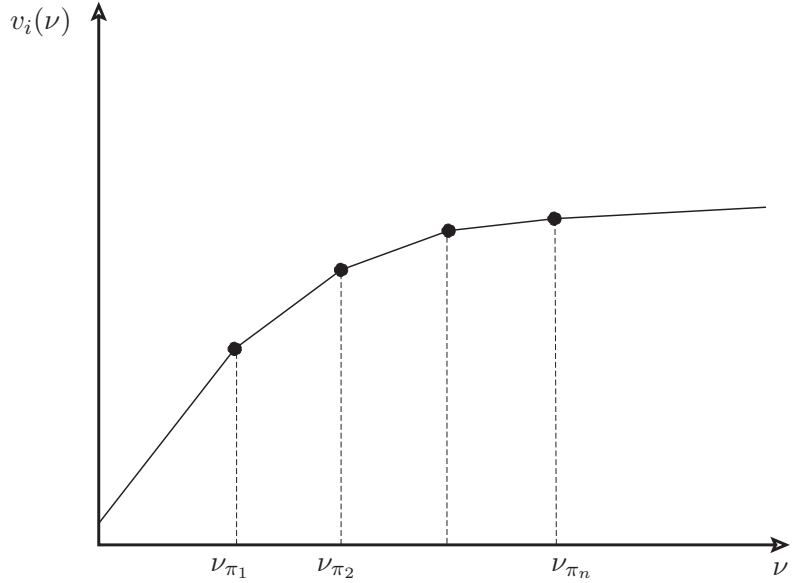}
\caption{Dependence on activity charge $\nu$ of optimal value function $v_i(\nu)$.}
\label{fig:ovgcp}
\end{figure}

The next result, illustrated in Figure \ref{fig:ovgcp}, 
characterizes  
$v_i(\nu)$.
Let $S_{n+1}~=~\emptyset$.

\begin{corollary}
\label{cor:rbvr}
Function $v_i(\nu)$ 
 is continuous, concave and piecewise linear on
$\nu$, and
\begin{align*}
v_i(\nu) & = \min \, \left\{v_i^{S_k}(\nu): 0 \leq k \leq n\right\}  \\
& = \begin{cases} v_i^{S_1}(\nu) = v_i^{S_1} + \nu \, b_i^{S_1}, &
  \text{if } \nu \in
  (-\infty,    \nu_{\pi_1}] \\
    v_i^{S_{k}}(\nu) = v_i^{S_k} + \nu \, b_i^{S_k}, & \text{if }
    \nu \in [\nu_{\pi_{k-1}}, \nu_{\pi_k}], 
 2 \leq k \leq n  \\
    v_i^{S_{n+1}}(\nu) = v_i^{S_{n+1}} + \nu \, b_i^{S_{n+1}}, &
  \text{if } \nu \in 
[\nu_{\pi_n}, +\infty).
    \end{cases}
\end{align*}
\end{corollary}
\begin{proof} 
The identities follow from Theorem
\ref{the:indexrb} and Definition \ref{def:ti}. They imply
$v_i(\nu)$ is continuous concave piecewise linear on
$\nu$,  being the minimum of linear functions of $\nu$. 
\end{proof}
\qed

\subsection{Marginal costs}
\label{s:mcia}
Recall that in Section \ref{s:aaia} we introduced  
\emph{marginal costs} $c^S_j$'s to simplify
index calculations.
This section discusses further properties of
such coefficients in the RB setting. 

We start by defining coefficients $c^S_j$, for 
$j \in N$, $S \subseteq N^{\{0, 1\}}$, in terms of  
value measure $v_i^S$.
For every $S$,
 the $v_i^{S}$'s are characterized by the linear equations
\begin{equation*}
v_i^S =
\begin{cases} 
\displaystyle h_i^1 +
 \beta \,\sum_{j \in N} p_{ij}^1 \, v^{S}_j & \text{if } i\in S \\
\displaystyle h_i^0 + 
 \beta \, \sum_{j \in N} p_{ij}^0 \, v^{S}_j & \text{if } i \in N
 \setminus S;
\end{cases}
\end{equation*}
or, in vector notation, 
\begin{equation}
\label{eq:jmn1}
\begin{split}
\mathbf{v}_{S}^{S} &= \mathbf{h}_{S}^1
 +
 \beta \, \mathbf{P}_{S N}^{1} \, 
\mathbf{v}^{S}
 \\
\mathbf{v}_{N \setminus S}^{S} & =  \mathbf{h}_{N \setminus S}^0 + 
 \beta \, \mathbf{P}_{N \setminus S, N}^{0} \, 
\mathbf{v}^{S}.
\end{split}
\end{equation} 
Define now
\begin{equation}
c_i^S = h_i^0 - h_i^1 
+
\beta \, \sum_{j \in N} (p_{ij}^0 - p_{ij}^1) \, v^{S}_j, 
\quad i \in N,
\label{eq:jaiscdef}
\end{equation}
i.e.,
\begin{equation}
\label{eq:jasvec}
\mathbf{c}^S = \mathbf{h}^0 - \mathbf{h}^1
 + \beta \, \left(\mathbf{P}^0 - 
  \mathbf{P}^1\right) \, \mathbf{v}^S,
\end{equation}
Coefficient $c_i^S$ thus represents the \emph{marginal
increment in cost measure $v^S$
resulting from a passive-to-active action interchange in initial state
$i$.}

It immediately follows that 
\begin{equation}
\label{eq:ciseq0}
c^S_j = 0, \quad j \in N^{\{1\}}. 
\end{equation}
Furthermore, 
(\ref{eq:jmn1})--(\ref{eq:jasvec}) readily yields the following 
counterpart of Lemma \ref{lma:reltws}.
\begin{lemma}
\label{lma:civisrel}
The following identities hold:
\begin{equation}
\begin{split}
\mathbf{h}^0 - (\mathbf{I} - \beta \, \mathbf{P}^0) \, 
 \mathbf{v}^S  & =  
\left[\begin{array}{c} \mathbf{c}_S^S \\ \mathbf{0}_{N \setminus S}
\end{array}\right]
 \label{eq:faisrb} \\
(\mathbf{I} - \beta \, \mathbf{P}^1) \, 
 \mathbf{v}^S - 
\mathbf{h}^1
 & = 
\left[\begin{array}{c}
 \mathbf{0}_{S} \\
 \mathbf{c}_{N \setminus S}^S\end{array}\right]. 
\end{split}
\end{equation}
\end{lemma}

The next result is a cost analog of Proposition \ref{lem:keys}.

\begin{proposition}[Cost decomposition laws]
\label{lem:ckeys} 
For $u \in \mathcal{U}$ and $S \subseteq N^{\{0, 1\}}$, 
\begin{equation}  \label{eq:clcssc}
v_i^S + 
\sum_{j \in S} c^{S}_j \, x_{ij}^{0, u} 
=
v_i^u + 
\sum_{j \in N^{\{0, 1\}} \setminus S} c^{S}_j \, x_{ij}^{1, u}.
\end{equation}
\end{proposition}
\begin{proof}
Using in turn 
equations (\ref{eq:slpc}) and (\ref{eq:faisrb}), we have
\begin{equation*}
\begin{split}
0 & = \left[
\mathbf{x}_i^{0, u} \,  (\mathbf{I} - \beta \, \mathbf{P}^0) + 
\mathbf{x}_i^{1, u} \,  (\mathbf{I} - \beta \,
\mathbf{P}^1) - \mathbf{e}_i\right] \, 
\mathbf{v}^S \\
& =  
\mathbf{x}_i^{0, u} \, \left[(\mathbf{I} - \beta \, \mathbf{P}^0)
\, \mathbf{v}^S - \mathbf{h}^0\right] + 
\mathbf{x}_i^{1, u} \,
\left[(\mathbf{I} - \beta \, \mathbf{P}^1) \, \mathbf{v}^S
 - \mathbf{h}^1\right] - \\
& \qquad
\mathbf{e}_i \, \mathbf{v}^S + 
\mathbf{x}_i^{1, u} \, \mathbf{h}^1 + 
\mathbf{x}^{0, u}(i) \, \mathbf{h}^0 \\
& = 
-\mathbf{x}_{i, S}^{0, u} \, \mathbf{c}_S^S +
\mathbf{x}_{i, N \setminus S}^{1, u} \, \mathbf{c}_{N
  \setminus S}^S - 
v_i^S + 
v_i^u,
\end{split}
\end{equation*}
which yields the result, using (\ref{eq:ciseq0}).
\end{proof}
\qed

The relation 
between coefficients $v^S_j$'s and $c^S_j$'s is  clarified next
(cf. Corollary \ref{cor:tsjwsj}).

\begin{corollary}
\label{cor:cjsvjs}
The following identities hold: for $i \in N$ and $S \subseteq N^{\{0,
  1\}}$, 
\begin{align*}
v_i^S & = v_i^{S \cup \{j\}} + c^{S}_j \, x_{ij}^{1, S
  \cup \{j\}}, \quad j \in N^{\{0, 1\}} \setminus S \\
v_i^{S \setminus \{j\}} & = v_i^S + 
  c^{S}_j \, x_{ij}^{0, S \setminus \{j\}}, \quad j
  \in S.
\end{align*}
\end{corollary}
\begin{proof}
It follows by letting $u = S \cup \{j\}$ and 
$u = S \setminus \{j\}$ in Proposition \ref{lem:ckeys}, respectively.
\end{proof}
\qed

The next result sheds further light on the relation between time and value measures,
and between marginal workloads and marginal costs. 

\begin{proposition}
\label{lma:keyidvec}
Under Assumption \ref{ass:pcl1}, the following holds: for $j \in S \in \mathcal{F}$,
 \\
\begin{itemize}
\item[(a)] 
$ \mathbf{v}^{S \setminus \{j\}} - \mathbf{v}^{S}  = 
\frac{\displaystyle{c^{S}_j}}{\displaystyle{w^{S}_j}} \, 
  \left(\mathbf{b}^{S} - \mathbf{b}^{S \setminus \{j\}}\right) = 
\frac{\displaystyle{c^{S \setminus \{j\}}_j}}{\displaystyle{w^{S
      \setminus \{j\}}_j}} \, 
  \left(\mathbf{b}^{S} - \mathbf{b}^{S \setminus \{j\}}\right).
$ 
\\
\item[(b)] 
$\frac{\displaystyle{c^{S}_j}}{\displaystyle{w^{S}_j}}
=
\frac{\displaystyle{c^{S \setminus \{j\}}_j}}{\displaystyle{w^{S \setminus \{j\}}_j}}.$ \\
\item[(c)]
$
\mathbf{c}^{S} - \mathbf{c}^{S \setminus \{j\}} = 
 \frac{\displaystyle{c^{S}_j}}{\displaystyle{w^{S}_j}} \,
\left(\mathbf{w}^{S} - \mathbf{w}^{S \setminus \{j\}}\right).
$
\end{itemize}
\end{proposition}
\begin{proof} 
(a) This part follows from Proposition \ref{lem:keys} and 
Proposition \ref{lem:ckeys}.

(b) The result follows from (a) and Proposition \ref{pro:iuwjsg0}. 

(c) From (\ref{eq:vecwisdef}), we readily obtain
\begin{equation}
\label{eq:wsktsk}
\mathbf{w}^{S} - \mathbf{w}^{S \setminus \{j\}} = 
\beta \, 
   \left(\mathbf{P}^1 - \mathbf{P}^0\right) \,
 \left(\mathbf{b}^{S} - \mathbf{b}^{S \setminus \{j\}}\right).
\end{equation}
Similarly, by (\ref{eq:jasvec}), we have 
\begin{equation}
\label{eq:cskvsk}
\mathbf{c}^{S} - \mathbf{c}^{S \setminus \{j\}} = 
 \beta \, \left(\mathbf{P}^0 - 
  \mathbf{P}^1\right) \,
\left(\mathbf{v}^{S} - \mathbf{v}^{S \setminus \{j\}}\right).
\end{equation}
The result now follows by combining part (a) with
(\ref{eq:wsktsk})--(\ref{eq:cskvsk}). 
\end{proof}
\qed

\begin{remark} \hspace{1in}
\begin{enumerate}
\item
Proposition \ref{lma:keyidvec} shows that the $c^S_j$'s
defined by (\ref{eq:jaiscdef}) extend those defined by (\ref{eq:cjsdef}).
\item It follows by construction that the $c^S_j$'s 
are \emph{symmetric} (cf. Definition \ref{ass:symmc}). Therefore, 
marginal workloads satisfy the recursion in Proposition
\ref{pro:wjscalc}.
\item
Note that, by combining identities
(\ref{eq:hatc0}), (\ref{eq:jmn1}) and (\ref{eq:faisrb}), it
 follows
that
\begin{equation}
\label{eq:h0cj}
\widehat{\mathbf{h}}^0 = \mathbf{c}^{J}.
\end{equation}
\end{enumerate}
\end{remark}

\subsection{PCL-indexability as a law of diminishing marginal returns}
\label{s:pclidmr}
This section discusses the intuitive interpretation of 
PCL-indexability (cf. Definition \ref{def:pclirb}) as a form of the 
classic economic law of \emph{diminishing
marginal
returns}.
Suppose  the project is PCL-indexable as above, 
 and let
$\boldsymbol{\pi}$, $\boldsymbol{\nu}$ and $S_k$ be as in Section
\ref{s:rtpsp}.
Assume the initial state 
is drawn from a probability distribution assigning a
\emph{positive mass $p_i > 0$ to each state $i \in N$}. 
Write $\mathbf{p} = (p_i)_{i \in N}$, 
$b^S = \sum_{i \in N} p_i \, b_i^S$, and $v^S = \sum_{i \in N} p_i \, v_i^S$.

\begin{theorem}[Index characterization and diminishing marginal returns]
\label{the:dmr}
\begin{itemize}
\item[(a)] 
\[
b^{S_{n+1}}  < b^{S_n} < \cdots <
b^{S_1}.
\]
\item[(b)]
For $1 \leq k \leq n$, dynamic allocation index $\nu_{\pi_k}$ is given by
\begin{align*}
\nu_{\pi_k} & = \frac{v^{S_{k+1}} -
  v^{S_k}}{b^{S_{k}} -
b^{S_{k+1}}} \\
& = 
\min \, \left\{\frac{v^{S_{k} \setminus \{j\}} -
  v^{S_k}}{b^{S_{k}} - b^{S_{k} \setminus \{j\}}}:  j \in S_k\right\} \\
& = \max \, 
 \left\{\frac{v^{S_k} - v^{S_k \cup \{j\}}}{b^{S_k \cup \{j\}} -
  b^{S_k}}: j \in N^{\{0, 1\}} \setminus S_k\right\}.
\end{align*}
\item[(c)] Diminishing marginal returns: 
\[
\frac{v^{S_{2}} -
  v^{S_1}}{b^{S_{1}} - b^{S_{2}}} \leq 
\frac{v^{S_{3}} -
  v^{S_2}}{b^{S_{2}} - b^{S_{3}}} \leq \cdots \leq 
\frac{v^{S_{n+1}} -
  v^{S_n}}{b^{S_{n}} - b^{S_{n+1}}}.
\]
\end{itemize}
\end{theorem}
\begin{proof}
(a) This part follows from Proposition \ref{pro:iuwjsg0} and $\mathbf{p}
> \mathbf{0}$. 

(b) The first identity follows from Proposition \ref{lma:keyidvec}, identity
(\ref{eq:gampinmkcw}), and part (a). The second identity then follows from
(\ref{eq:1stic}) in Proposition \ref{the:gpimin}. The third identity
further follows from Corollary \ref{cor:rbvr}.

(c) The result follows from parts (a), (b) and the inequalities in (\ref{eq:gleqpk}).
\end{proof}
\qed

\begin{remark} \hspace{1in}
\begin{enumerate}
\item Part (a) shows that the \emph{busy} or \emph{active} 
time (as
  measured by $b^u$) is strictly increasing along the set/policy sequence 
$\emptyset = S_{n+1} \subset S_n \subset \cdots \subset S_1 = N^{\{0, 1\}}$.
\item Part (b) characterizes index $\nu_{\pi_k}$
  as a \emph{locally optimal marginal cost rate}:  it
  is the \emph{minimal rate of marginal cost increase from $v^{S_k}$ 
 per unit marginal activity decrease from $b^{S_k}$ resulting from  
 an active-to-passive action interchange on some 
  state $j \in S_k$}. Furthermore, $\nu_{\pi_k}$ is the 
 \emph{maximal rate of marginal cost decrease from $v^{S_k}$ per unit
 marginal activity increase from $b^{S_k}$ resulting from a
 passive-to-active action interchange on some state $j \in N^{\{0,
   1\}} \setminus S_k$}.
\item Part (c) shows that the optimal rate of marginal cost decrease per unit
  marginal active time increase  diminishes on the base active time. It
thus represents a form of the law of diminishing
  marginal returns. 
\end{enumerate}
\end{remark}

\begin{figure}[ht]
\centering
\psfragscanon
\psfrag{F}[u]{$b^u$}
\psfrag{H}[u]{$b^{S_1}$}
\psfrag{E}[ul]{$b^{S_{k-1}}$}
\psfrag{D}[u]{$b^{S_4}$}
\psfrag{C}[u]{$b^{S_k}$}
\psfrag{D}[u]{$b^{S_k \cup \{j\}}$}
\psfrag{B}[llu]{$b^{S_{k+1}}$}
\psfrag{A}[u]{$b^{S_k \setminus \{i\}}$}
\psfrag{K}[u]{$b^{S_{n+1}}$}
\psfrag{I}[c]{$v^u$}
\includegraphics[height=3in,width=4.8in,keepaspectratio]{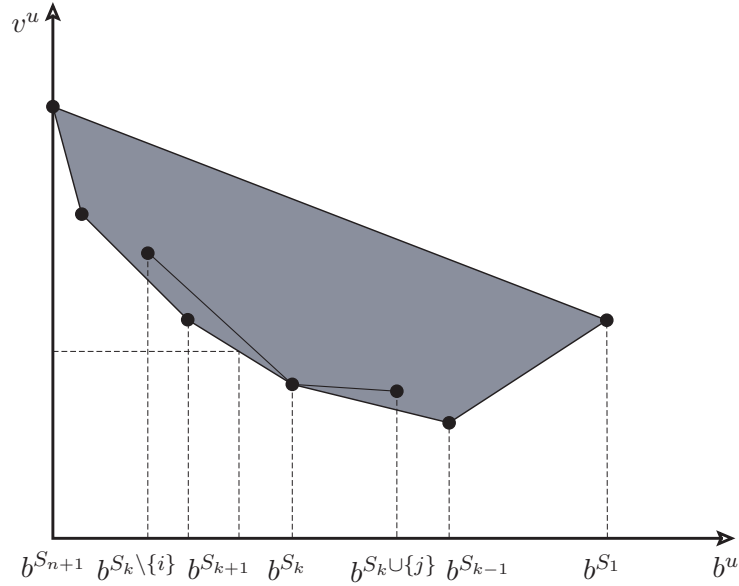}
\caption{Activity-cost plot: PCL-indexability and diminishing marginal returns.}
\label{fig:dmr}
\end{figure}

Figure \ref{fig:dmr} illustrates the result by an 
\emph{activity-cost plot}, where
 the shaded area represents the \emph{region of achievable
activity-cost pairs $(b^u, v^u)$}. 
We further have the following index characterization under
Assumption \ref{ass:imh} on \emph{nondecreasing
marginal workloads}.
\begin{theorem}
\label{the:icvt}
Under Assumption \ref{ass:imh}, 
\begin{equation}
\label{eq:nuj}
\nu_j = \max \, \left\{\frac{v^{S \setminus \{j\}}_j - v^{S}_j}{b^S_j
    - b^{S \setminus \{j\}}_j}: j \in S \in \{S_1, \ldots,
    S_n\}\right\}, \quad j \in N^{\{0, 1\}}.
\end{equation}
\end{theorem}
\begin{proof}
The result follows directly from Theorem \ref{the:ic} and Proposition
\ref{lma:keyidvec}.
\end{proof}
\qed

\begin{remark} \hspace{1in}
\begin{enumerate}
\item Theorem \ref{the:icvt} represents an RB counterpart
  of the \emph{Gittins index} characterization for classic bandits (
$\mathbf{P}^0 = \mathbf{I}$ and $\mathbf{h}^0 = \mathbf{0}$)
 as an \emph{optimal average cost (reward) rate per unit
  time}, first given in \cite{gi79}. 
In the classic case Theorem \ref{the:icvt} gives
\[
\nu_j = \max \, \left\{\frac{- v^{S}_j}{b^S_j}: j \in S \in \{S_1, \ldots,
    S_n\}\right\}, \quad j \in N^{\{0, 1\}},
\]
since $b^{S \setminus \{j\}}_j = v^{S \setminus \{j\}}_j = 0$.
Actually, the Gittins index characterization in \cite{gi79} is
\[
\nu_j = \max \, \left\{\frac{- v^{S}_j}{b^S_j}: j \in S \in 2^N\right\}, \quad j \in N,
\]
\item We pose the \emph{open problem}: Find conditions under
  which (\ref{eq:nuj}) extends to
\[
\nu_j = \max \, \left\{\frac{v^{S \setminus \{j\}}_j - v^{S}_j}{b^S_j
    - b^{S \setminus \{j\}}_j}: j \in S \in \mathcal{F}\right\}, \quad j \in N^{\{0, 1\}}.
\]
\end{enumerate}
\end{remark}

\subsection{Extension to the long-run average criterion}
\label{s:ocsac}
The results above for the time-discounted criterion readily 
extend to the long-run average criterion under suitable \emph{ergodicity}
conditions, by standard limiting (Tauberian) arguments. 
Assume that the model is \emph{communicating}, i.e., every state can
be reached from every other state under some stationary policy. 
Assume further that, for every $S \in \mathcal{F}$, the $S$-active
policy is \emph{unichain}, i.e., it induces a single recurrent class
plus a (possibly empty) set of transient states.
 Then, it is well known that
measures $b_{i}^u(\beta)$, $v_{i}^u(\beta)$, $x_{ij}^{a, u}(\beta)$ 
(where we have  made explicit the dependence on $\beta$), when scaled by
factor $1 - \beta$,  converge to limiting
values independent of the initial state
$i$, given by 
\[
\bar{b}^u = 
\lim_{T \to \infty} \, 
\frac{1}{T} \, E_i^u\left[\sum_{t=0}^T \theta_{X(t)}^1 \, a(t)\right]
= 
\lim_{\beta \nearrow 1} \, 
(1-\beta) \, b_i^u(\beta),
\]
\[
\bar{v}^u = 
\lim_{T \to \infty} \, 
\frac{1}{T} \, E_i^u\left[\sum_{t=0}^T h_{X(t)}^{a(t)}\right]
= 
\lim_{\beta \nearrow 1} \, 
(1-\beta) \, v_i^u(\beta),
\]
\[
\bar{x}^{a, u}_j = \lim_{T \to \infty} \, 
\frac{1}{T} \, E_i^u\left[\sum_{t=0}^T 1\{X(t) = j, a(t) = a\}\right]
= 
\lim_{\beta \nearrow 1} \, 
(1-\beta) \, x_{ij}^{a, u}(\beta).
\]
Hence, $\bar{b}^u$, $\bar{v}^u$ and $\bar{x}^{a, u}_j$ are the corresponding
\emph{long-run average}, or \emph{steady-state}, measures. 

We next argue that 
 the \emph{unscaled} quantities $w_i^S(\beta)$ and 
$c_i^S(\beta)$ converge to finite limits
$\bar{w}_i^S$ and $\bar{c}_i^S$ as $\beta \nearrow 1$.
Start with marginal workload $w_i^S(\beta)$. 
We can write, for 
$i \in N$ and $S \in \mathcal{F}$, 
\begin{equation}
\label{eq:tsund}
b_i^S(\beta) = \frac{\bar{b}^S}{1 - \beta} + a_i^S + O(1-\beta), \quad \text{as }
\beta \nearrow 1,
\end{equation}
where the values $a_i^S$ are determined, up to an additive constant, by
the equations
\[
\bar{b}^S + a_i^S = 
\begin{cases}
\displaystyle \theta_i^1 + \sum_{j \in N} p_{ij}^1 \, a^S_j & 
\text{if } i \in S \cup N^{\{1\}} \\ \\
\displaystyle \sum_{j \in N} p_{ij}^0 \, a^S_j & \text{if } i \in N^{\{0, 1\}}
  \setminus S.
\end{cases}
\]
Now, substituting for $b_i^S(\beta)$ as given by (\ref{eq:tsund}) 
in (\ref{eq:wisdef}), and letting 
$\beta \nearrow 1$, gives
\[
\bar{w}_i^S = \lim_{\beta \nearrow 1} w_i^S(\beta) = 
\theta_i^1 \, 1\{i \in N^{\{0, 1\}}\}
  + \sum_{j \in N} (p_{ij}^1 - p_{ij}^0) \, 
  a^S_j, \quad i \in N.
\]

We proceed analogously with  marginal costs $c_i^S(\beta)$. 
Write
\begin{equation}
\label{eq:vsund}
v_i^S(\beta) = \frac{\bar{v}^S}{1 - \beta} + f_i^S + O(1-\beta), \quad \text{as }
\beta \nearrow 1,
\end{equation}
where the values $f_i^S$ are determined, up to an additive  constant, by
the equations
\[
\bar{v}^S + f_i^S = 
\begin{cases}
\displaystyle h_i^1 + \sum_{j \in N} p_{ij}^1 \, f^S_j & 
\text{if } i \in S \\ \\
\displaystyle h_i^0 + 
\sum_{j \in N} p_{ij}^0 \, f^S_j & \text{if } i \in N
  \setminus S.
\end{cases}
\]
Now, substituting for $v_i^S(\beta)$ as given by (\ref{eq:vsund}) 
in (\ref{eq:jaiscdef}), and letting 
$\beta \nearrow 1$, gives
\[
\bar{c}_i^S = \lim_{\beta \nearrow 1} c_i^S(\beta) = 
h_i^0 - h_i^1 
  + \sum_{j \in N} (p_{ij}^0 - p_{ij}^1) \, 
  f^S_j, \quad i \in N.
\]

Thus, previous results carry over to the
long-run average case. 

\subsection{Optimal control subject to an activity constraint}
\label{s:ocsac2}
In applications, it is often of interest to impose a constraint 
on the mean rate of activity. See, e.g., 
\cite{hordspi} and the references therein. 
This is particularly relevant under the 
\emph{long-run average} criterion discussed above, on
which
we focus next.

The 
\emph{constrained control problem} of concern is 
to find a stationary policy 
 minimizing cost measure $\bar{v}^u$, among those whose
long-run average activity rate is $\bar{b}^u = t$:
\begin{equation}
\label{eq:ccp}
\bar{v}_t = \min \, \left\{\bar{v}^u: \bar{b}^u = t, u \in \mathcal{U}\right\}.
\end{equation}

Assume the project is PCL-indexable as in Section \ref{s:pclidmr}, and let 
$\boldsymbol{\pi}$, $\boldsymbol{\nu}$ and $S_k$ be its
optimal $\mathcal{F}$-string, index vector and 
active sets. 
Suppose that, for some $1 \leq k \leq n$,  
\[
\bar{b}^{S_{k+1}} < t < \bar{b}^{S_k},
\]
and let 
\[
p = \frac{t - \bar{b}^{S_{k+1}}}{\bar{b}^{S_k} - \bar{b}^{S_{k+1}}}, \qquad
q = 1 - p.
\]
Denote by $(S_{k+1}, \pi_k, p)$ the stationary policy that
is: active on states $j \in S_{k+1} \cup N^{\{1\}}$;
active on state $\pi_k$ with probability $p$; and 
passive otherwise. 
The next result follows immediately from  Section
\ref{s:pclidmr}, and hence its proof is omitted. See also Figure 
\ref{fig:dmr}. 

\begin{proposition}
\label{pro:occ}
The following holds: \\
\begin{itemize}
\item[(a)] Policy $(S_{k+1}, \pi_k, p)$ is optimal for problem
(\ref{eq:ccp}); its  optimal value is
\[
\bar{v}_t = (1-p) \, \bar{v}^{S_{k+1}} + p \, \bar{v}^{S_k}.
\]
\item[(b)] Function $\bar{v}_t$ is piecewise linear concave on $t$, with 
\[
\frac{d}{dt} \bar{v}_t = \nu_{\pi_k}, \quad \bar{b}^{S_{k+1}} < t < \bar{b}^{S_k}.
\]
\end{itemize}
\end{proposition}

\begin{remark}
Proposition \ref{pro:occ}(b) characterizes
the index $\nu_{\pi_k}$ as a derivative of the optimal 
constrained value function $v_t$ with respect to the required 
activity level $t$. 
\end{remark}

\section{Admission control problem: PCL-indexability analysis}
\label{s:qaca}
This section returns to the 
admission control model introduced in
Section \ref{s:tp}.
We shall resolve the issues raised in Section \ref{s:cof}
by deploying a PCL-indexability analysis.
See the Appendix  for
important yet ancillary material relevant to this section.
 
In what follows, we shall write $\Delta x_i =
x_i - x_{i-1}$, 
$d_i = \mu_i - \lambda_i$,
and 
\[
\rho_i = \frac{\lambda_i}{\mu_{i+1}}.
\]
We next state the regularity conditions we shall require of model
parameters.

\begin{assumption}
\label{ass:costconv}
The following conditions hold: \\
{\rm (i)} Concave nondecreasing $d_i$:
$ 0 \leq \Delta d_{i+1} 
   \leq \Delta d_i, \quad 1 \leq i \leq
   n-1$, and $ \Delta d_1 > 0$. \\
{\rm (ii)} Convex nondecreasing $h_i$:
$ \Delta h_{i+1} \geq \Delta h_i \geq 0, \quad 1 \leq i \leq n-1.$
\end{assumption}

\begin{remark}
Assumption \ref{ass:costconv} is significantly
less restrictive than Chen and Yao's conditions in \cite{chenyao}.
Besides requiring $d_i$ to be nondecreasing in their condition
 (5.5a), they require $\mu_i$ to be concave nondecreasing in their
 condition (5.5b). They  further impose additional conditions, including
 \emph{linearity} of holding costs.
\end{remark}

\subsection{PCL-indexability analysis under the discounted criterion}
\label{s:adc}
We shall establish PCL-indexability of the model  relative to the family of
\emph{threshold policies},  given by set system $(N^{\{0, 1\}}, \mathcal{F})$, where 
$\mathcal{F} = \left\{S_1, \ldots, S_{n+1}\right\}$
 is given by (\ref{eq:facp})--(\ref{eq:skdefac}).
\emph{Activity
measure} $b_i^u$ is  given by (\ref{eq:tuiacm}), which corresponds to
letting
$\theta_j^1 = \lambda_j / (\alpha + \Lambda)$ in (\ref{eq:tup}), 
 for $j \in N$, where $\Lambda$ is the \emph{uniformization rate} (cf. Appendix \ref{s:acpdtr}). 

\begin{figure}[tbp]
\fbox{%
\begin{minipage}{\textwidth}
{\bf Calculation of $w_i^{S_1}$:}
\begin{equation*}
\frac{w_0^{S_1}}{\lambda_0}  = \frac{\alpha + \Delta d_1 }{\alpha + \mu_1}; \quad
\frac{w_i^{S_{1}}}{\lambda_i}  = 
\frac{\displaystyle \alpha + \Delta d_{i+1}  +  
  \frac{w^{S_{1}}_{i-1}}{\rho_{i-1}}}{\alpha + \mu_{i+1}}, 
\quad 1 \leq i \leq n-1
\end{equation*}
{\bf Calculation of $w^{S_{2}}_i$:}
\begin{equation*}
\frac{w^{S_{2}}_0}{\lambda_0}  = 
\frac{\alpha + \Delta d_1}{\alpha + \lambda_0 + \mu_1}; \quad
\frac{w^{S_{2}}_i}{\lambda_i}  = 
\frac{\displaystyle \alpha + \Delta d_{i+1} +  
  \frac{w^{S_{2}}_{i-1}}{\rho_{i-1}}}{\alpha + \mu_{i+1}},
\quad 1 \leq i \leq n-1
\end{equation*}
{\bf Calculation of $w^{S_{k+1}}_i$'s, for $2 \leq k \leq n$:}
\begin{equation*}
\begin{split}
\frac{w^{S_{k+1}}_{k-1}}{\lambda_{k-1}} & = \frac{1}{a_k} \,
\frac{\displaystyle \alpha + \Delta d_k + 
  \frac{w^{S_k}_{k-2}}{\rho_{k-2}}}{\alpha + \lambda_{k-1} + \mu_k}; 
\quad \frac{w^{S_{k+1}}_{k-2}}{\rho_{k-2}}  =
-(\alpha + \Delta d_k) + 
  \frac{\alpha + \lambda_{k-1} + \mu_k}{\lambda_{k-1}} \, 
  w^{S_{k+1}}_{k-1}  \\
\frac{w^{S_{k+1}}_i}{\lambda_i} & =
\frac{\displaystyle \alpha + \Delta d_{i+1}  + 
\frac{w^{S_{k+1}}_{i-1}}{\rho_{i-1}}}{\alpha + \mu_{i+1}}, 
\quad k \leq i \leq n-1  \\
\frac{w^{S_{k+1}}_i}{\rho_i} & =
- (\alpha + \Delta d_{i+2})  + 
\frac{\alpha + \lambda_{i+1} + \mu_{i+2}}{\lambda_{i+1}} \, 
w^{S_{k+1}}_{i+1}  - 
w^{S_{k+1}}_{i+2}, 
\quad 0 \leq i \leq k-3
\end{split}
\end{equation*}
\end{minipage}}
\caption{Recursive calculation of marginal workloads 
$w^{S_k}_i$.}
\label{fig:aibcomp2}
\end{figure}

We must first calculate marginal workload coefficients $w^{S_k}_i$,
for which a complete recursion is given in Figure
\ref{fig:aibcomp2}.
It involves coefficients $a_i$, given by 
(\ref{eq:akdef}).

\begin{proposition}
\label{the:aisrp} 
Marginal workloads $w_i^S$, for $i \in N^{\{0, 1\}}$ and 
$S \in \mathcal{F}$, are calculated by
the recursion shown in Figure \ref{fig:aibcomp2}.
\end{proposition}
\begin{proof}
The result follows by reformulating 
in terms of the $w_i^S$'s the equations on
terms $\Delta b^S(i)$'s given in Lemma \ref{lma:1actdeltap} and
Lemma \ref{lma:keyreldeltat} in Appendix \ref{s:cmw}, using identity (\ref{eq:wistisrel}).
\end{proof}
\qed

\begin{remark}
The recursion in Figure \ref{fig:aibcomp2} further
yields coefficients $w_i^S$ when 
$\alpha = 0$. These are the 
\emph{long-run average marginal workloads} discussed in Section \ref{s:pclitac}.
\end{remark}

The next result establishes the required properties of marginal workloads.

\begin{proposition}[Positive nondecreasing $w_i^S$'s]
\label{pro:pcl1mod}
Under Assumption \ref{ass:costconv}(i): \\
(a)
$
w_i^S > 0$, for  $i \in N^{\{0, 1\}}, S \in \mathcal{F},
$
and hence Assumption
\ref{ass:pcl1} holds.  \\
(b) $w_i^S$ is nondecreasing on $S \in
\mathcal{F}$, for $i \in S$ fixed, 
and hence
Assumption \ref{ass:imh} holds.
\end{proposition}
\begin{proof}
Both parts follow directly from Lemma \ref{lma:aiigt0} in Appendix \ref{s:cmw}.
\end{proof}
\qed

Figure \ref{fig:caibc} illustrates the recursions
and  inequalities established in Appendix \ref{s:cmw} 
on marginal workloads (arrows indicate the direction of
calculations). 
\emph{Pivot} terms, forming the backbone of the recursion,  are enclosed in boxes. 

\begin{figure}[tbp]
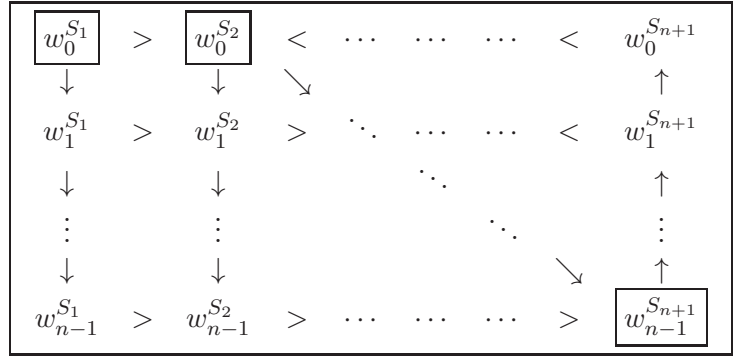

\begin{center}
\boxed{ \begin{tabular}{ccccccccc}  
$\boxed{w^{S_1}_{0}}$ & $>$ 
& $\boxed{w^{S_{2}}_0}
$ & $<$ & $\cdots$ & $\cdots $ & $\cdots$
& $< $   & $w^{S_{n+1}}_0$
 \\
$\downarrow$ & & $\downarrow$ 
 & $\searrow$ & & & &    & $\uparrow$ \\
$w^{S_1}_1$ & $>$ & $w^{S_{2}}_1$ & $>$
 & $\ddots $ & $\cdots$ & $\cdots$ & $< $  & 
  $w^{S_{n+1}}_1$   \\ 
   $\downarrow $ & & $\downarrow$ & &  & $\ddots$ & & 
 & $\uparrow$ \\
   $\vdots $ & & $\vdots$ & &  & & $\ddots$ & &   $\vdots$
 \\
$\downarrow $ & & $\downarrow$ & & &  & & $\searrow $ 
& $\uparrow$  \\
$w^{S_1}_{n-1}$  & $>$ & $w^{S_{2}}_{n-1}$ & $>$ & $\cdots$  &
$\cdots$ & $\cdots$ & $>$   & $\boxed{w^{S_{n+1}}_{n-1}}$ 
\end{tabular}}
\end{center}
\caption{Relations between marginal workloads $w^{S_k}_i$.}
\label{fig:caibc}
\end{figure}

Marginal cost analyses are given in Appendix 
\ref{s:mccc}, yielding the following recursion.

\begin{proposition}
\label{pro:mcc}
Marginal costs $c^{S_{k+2}}_k$, for $0 \leq k \leq n-1$,
are calculated by 
\begin{align*}
c^{S_2}_0  & = 
  \frac{\lambda_0}{\alpha + \lambda_0 + \mu_1} \, \Delta h_1 \\
c^{S_{k+2}}_k
 & = \frac{\lambda_k}{a_{k+1}} \, \frac{\displaystyle \Delta h_{k+1} + 
  \frac{c^{S_{k+1}}_{k-1}}{\rho_{k-1}}}{\alpha + \lambda_k +
  \mu_{k+1}}, \quad 1 \leq k \leq n-1.
\end{align*}
\end{proposition}

We are now ready to establish the model's
PCL-indexability, and to 
calculate its indices.
Construct  $\nu_0, \ldots, \nu_{n-1}$ recursively by
\begin{equation}
\label{eq:windexqs}
\begin{split} 
\nu_0 & = \frac{\Delta h_1}{\alpha + \Delta d_1}  \\
\nu_j & = 
\nu_{j-1} + \frac{\displaystyle{\Delta h_{j+1} - \nu_{j-1} \, 
  (\alpha + \Delta d_{j+1}})}{\displaystyle{\alpha + \Delta d_{j+1} + 
  \frac{w^{S_{j+1}}_{j-1}}{\rho_{j-1}}}}, \quad 1 \leq j \leq n-1.
\end{split}
\end{equation}
We shall need the following preliminary result.

\begin{lemma}
\label{lma:cwineq}
Under Assumption \ref{ass:costconv}, the following holds:
\\
\begin{itemize}
\item[(a)]
$
\frac{\displaystyle \Delta h_j}{\displaystyle \alpha + \Delta d_j} \leq 
\frac{\displaystyle \Delta h_{j+1}}{\displaystyle \alpha + \Delta d_{j+1}}, \quad 1 \leq j \leq n-1.
$ \\
\item[(b)]
$
\nu_j \leq \frac{\displaystyle 
\Delta h_{j+1}}{\displaystyle \alpha + \Delta d_{j+1}}, \quad 0 \leq j \leq n-1.
$ \\
\item[(c)]
$ \nu_0 \leq \nu_1 \leq \cdots \leq \nu_{n-1}.$
\end{itemize}
\end{lemma}
\begin{proof} 
(a)
The result follows directly from 
Assumption \ref{ass:costconv}. 

(b)
Proceed by induction on $j$.
The case $j = 0$ holds by (\ref{eq:windexqs}).
Suppose now
$$
\nu_{j-1} 
 \leq \frac{\Delta h_j}{\alpha + \Delta d_j}.
$$
It then follows, by part (a), that
$$
\nu_{j-1} \leq \frac{\Delta h_{j+1}}{\alpha + \Delta d_{j+1}}.
$$
Notice now that  the last identity in
 (\ref{eq:windexqs}) can be reformulated as
\begin{equation*}
\nu_j  =
\frac{\Delta h_{j+1}}{\alpha + \Delta d_{j+1}} + 
\frac{\displaystyle{\frac{w^{S_{j+1}}_{j-1}}{\rho_{j-1}}}}{
  \displaystyle{\alpha + \Delta d_{j+1} + \frac{w^{S_{j+1}}_{j-1}}{\rho_{j-1}}}} \, 
\left[\nu_{j-1} - \frac{\Delta h_{j+1}}{\alpha + \Delta d_{j+1}}\right].
\end{equation*}
Since $\alpha + \Delta d_{j+1} > 0$ and $w^{S_{j+1}}_{j-1} > 0$, 
it follows from the last identity 
that
\begin{equation*}
  \nu_{j-1} \leq \frac{\Delta h_{j+1}}{\alpha + \Delta d_{j+1}}
 \Longleftrightarrow
\nu_j \leq \frac{\Delta h_{j+1}}{\alpha + \Delta d_{j+1}},
\end{equation*}
which completes the induction.

(c) This follows from parts (a) and (b), together with  (\ref{eq:windexqs}).
\end{proof}
\qed

We are now ready to establish the main result of this section.
\begin{theorem}[PCL-indexability: discounted criterion]
\label{the:wiservcon}
Under Assumption \ref{ass:costconv}, 
the admission control model is PCL-indexable relative to
threshold policies and
rejection measure $b^u$.
Its dynamic allocation indices are the $\nu_j$'s given by 
(\ref{eq:windexqs}), and satisfy (\ref{eq:nuj}). \\
\end{theorem}
\begin{proof} 
Using (\ref{eq:gampinmkcw}) and 
Proposition \ref{lma:keyidvec}(b), 
we must show that 
\[
\nu_j = \frac{c^{S_{j+2}}_j}{w^{S_{j+2}}_j}, \quad 
0 \leq j \leq n-1.
\]
This readily follows  by induction on $j$, drawing on
Proposition \ref{fig:aibcomp2} and Proposition \ref{pro:mcc}. 
Furthermore, Proposition \ref{pro:pcl1mod} and Theorem \ref{eq:nuj}
imply that the index satisfies (\ref{eq:nuj}).
\end{proof}
\qed

\subsection{PCL-indexability under the long-run average criterion}
\label{s:pclitac}
As in  Section \ref{s:ocsac}, the PCL-indexability analysis above 
extends to the long-run average version of the admission control model. 
The relevant rejection and cost measures are
\[
\bar{b}^u = \lim_{T \to \infty} \, \frac{1}{T} \, E^u\left[\int_0^T \lambda_{L(t)}
  \, a(t) \, dt\right],
\]
\[
\bar{v}^u = \lim_{T \to \infty} \, \frac{1}{T} \, E^u\left[\int_0^T 
  h_{L(t)} \, dt\right].
\]
The limiting values of $w^{S_k}_j$, $c^{S_k}_j$ and $\nu_j$ as $\alpha
\searrow 0$ (equivalent to letting $\beta \nearrow 1$ in
Section \ref{s:ocsac})  are
obtained by setting $\alpha = 0$ in the given recursions. The next
result follows.

\begin{corollary}[PCL-indexability: long-run average criterion]
\label{cor:wiservcon}
Under Assumption \ref{ass:costconv}, 
the admission control model, under the long-run average criterion, is PCL-indexable relative to
threshold policies and
rejection measure $\bar{b}^u$.
Its indices satisfy 
\begin{equation*}
\bar{\nu}_j = \max \, \left\{\frac{\bar{v}^{S \setminus \{j\}} - \bar{v}^{S}}{\bar{b}^S
    - \bar{b}^{S \setminus \{j\}}}: j \in S \in \{S_1, \ldots,
    S_n\}\right\}, \quad j \in N^{\{0, 1\}}.
\end{equation*}
\end{corollary}

\subsection{The case $\lambda_j = \lambda$, $\mu_j = \mu$, $\alpha =
  0$}
\label{s:ccasr}
This section derives the long-run average 
indices
when 
$\lambda_j = \lambda$, $\mu_j = \mu$. 
Note that $\rho_j = \rho = \lambda/\mu$. 
As we shall see in Section \ref{s:appl}, the case $\rho > 1$ is 
often of interest in applications. 

The following results follow easily by induction, and hence we omit their proof. 
Note first that the coefficients $a_j$, defined by
(\ref{eq:akdef}), are  given, for $1 \leq j \leq n-1$,   by 
 
\[
a_j = \frac{1}{1 + \rho} \, \frac{1 + \cdots + \rho^j}{1 + \cdots +
  \rho^{j-1}} = \begin{cases}
\displaystyle \frac{1}{1 + \rho} \, \frac{\rho^{j+1} -1}{\rho^j-1} &
\text{if } \rho \neq 1 \\ \\
\displaystyle \frac{1}{2} \, \frac{j+1}{j} & \text{if } \rho = 1.
\end{cases}
\]

Regarding marginal workloads, we have 
\begin{align*}
w^{S_1}_j & = \lambda, \quad 1 \leq j \leq n-1 \\
w^{S_2}_0 & = 
    \frac{1}{1 + \rho} \, \lambda \\
w^{S_{j+1}}_{j-1} & = \frac{1}{1 + \rho} \, \frac{w^{S_j}_{j-2}}{a_j},
\quad
  2 \leq j \leq n.
\end{align*}
Such recursion gives
\[
w^{S_{j+1}}_{j-1} = \frac{\lambda}{\displaystyle (1 + \rho)^j \,
  \prod_{i=1}^j a_i}
= \frac{\lambda}{1 + \cdots + \rho^j},
\quad
1 \leq j \leq n.
\]

Hence, index recursion (\ref{eq:windexqs}) 
reduces to
\begin{align*}
\nu_0 & = \frac{\Delta h_1}{\mu} \\
\nu_j & = \nu_{j-1} + \Delta h_{j+1} \, 
  \frac{1 + \cdots + \rho^j}{\mu}, \quad 1 \leq j \leq n-1,
\end{align*}
which yields
\begin{equation}
\label{eq:nukta}
\nu_j = \frac{1}{\mu} \sum_{i=1}^{j+1} \Delta h_i \, \left(1 + \cdots +
  \rho^{i-1}\right) = 
\begin{cases}
\displaystyle \frac{1}{\mu} \sum_{i=1}^{j+1} \Delta h_i \,
 \frac{\rho^i -1}{\rho-1}
 & \text{if } \rho \neq 1 \\ \\
\displaystyle \frac{1}{\mu} \sum_{i=1}^{j+1} i \, \Delta h_i & \text{if } \rho = 1.
\end{cases}
\end{equation}

\begin{remark}
\label{re:ndindx}
A consequence of (\ref{eq:nukta}) is that, in this setting,
\emph{the index is monotonic, and hence the model is PCL-indexable,
 under the relaxed assumption that
cost rates $h_i$ be only nondecreasing}:  they need \emph{not}
be convex as in Assumption \ref{ass:costconv}(ii).
\end{remark}

In the \emph{linear cost}
case $h_j = h \, j$, we  obtain
\begin{equation}
\label{eq:nukplhc}
\nu_j = 
\frac{h}{\mu} \sum_{i=1}^{j+1} \left(1 + \cdots +
  \rho^{i-1}\right) =
\begin{cases}
\displaystyle \frac{h}{\mu} \, \left[\frac{\rho^{j+2}-1}{(\rho-1)^2} -
 \frac{j+2}{\rho-1} 
 \right]
 & \text{if } \rho \neq 1 \\ \\
\displaystyle \frac{h}{\mu} \, \frac{(j+1) \, (j+2)}{2} & \text{if } \rho = 1.
\end{cases}
\end{equation}

In the \emph{quadratic cost} case $h_j = h \, j^2$, we obtain, when $\rho \neq 1$,

\begin{align}
\label{eq:nujqq}
\nu_j & = \frac{h}{\mu} \, \left[ 
\left( \frac{2j+1}{(\rho -1)^2}-\frac{2}{\left( \rho -1\right)
    ^{3}}\right) \rho ^{j+2} - \frac{j\left(j+2\right)}{\rho-1}
 +\frac{3}{(\rho -1)^2}+\frac{2}{
\left( \rho -1\right) ^{3}}\right]
\end{align}
and, when $\rho = 1$,
\begin{align*}
\nu_j &  =  \frac{h}{\mu} \, \frac{\left( j+1\right) \, \left( j+2\right) \,
  \left( 4j+3\right)}{6}.
\end{align*}

\section{Applications to routing and make-to-stock scheduling in queueing systems}
\label{s:appl}
In this section we apply the admission control  index  obtained in
Section \ref{s:qaca} to develop new \emph{heuristic} index policies for
two hard queueing control problems.

\subsection{An index policy for admission control and routing
  to parallel queues}
\label{s:irpmq}
Consider a system at which customers arrive as 
a Poisson stream with rate $\lambda$.
Upon arrival, a customer may be either rejected, or routed to one of
$m$ queues for service. 
Queue $k$ has a finite buffer holding at most $n_k$ customers.
 Its service times are exponential, with rate
$\mu_k(j_k)$ when it holds
$L_k(t) = j_k$ customers at time $t \geq 0$, for $j_k \in
N_k = \{0, \ldots, n_k\}$. When all buffers are full, an arriving
customer is lost. 

Customers in queue $k$ incur holding costs at rate $h_{k}(j_k)$ while $L_k(t) =
j_k$, discounted in time at rate $\alpha > 0$.
Furthermore, a \emph{rejection charge} $\nu$ is incurred per lost
customer.
The problem of concern is to find a stationary 
\emph{admission control and routing policy} prescribing
whether to admit each arriving customer and, if so, to which nonfull
queue to route it, in order to minimize the 
expected total discounted  sum of holding costs and rejection charges incurred
over an infinite horizon.

We assume  model parameters satisfy the following conditions.

\begin{assumption}
\label{ass:2costconv}
For $1 \leq k \leq m$, the following holds: \\
{\rm (i)} Concave nondecreasing $\mu_{k}(j_k)$:
$0 \leq \Delta \mu_{k}(j_k+1) \leq \Delta \mu_{k}(j_k), \quad 1 \leq j_k < n_k$.\\
{\rm (ii)} Convex nondecreasing $h_{k}(j_k)$:
$0 \leq \Delta h_{k}(j_k) \leq \Delta h_{k}(j_k+1), \quad 1 \leq j_k < n_k$.
\end{assumption}

We aim 
to design a well-grounded and tractable \emph{heuristic} policy, 
for which we shall use
the \emph{admission control  index} developed
in Section \ref{s:qaca}.
The idea is to note that this model is an RBP made up of $m$ single-queue
\emph{admission control RBs} as studied before where,
\emph{at each time, at most one of the $m$ entry gates must be open.}

Let $\nu_{k}(j_k)$ be queue $k$'s admission control index,
representing the 
\emph{fair rejection charge} for a customer finding queue $k$ in state $j_k < n_k$.
Such interpretation leads to the following \emph{admission control and
  routing index policy}: 
\begin{enumerate}
\item Route an arriving customer to a nonfull queue $k$ whose current state
 $j_k < n_k$ has the \emph{smallest
  index} $\nu_{k}(j_k)$ satisfying $\nu_k(j_k) < \nu$, if any is available.
\item Otherwise, reject the customer.
\end{enumerate}

In the case where queues are symmetric (ignoring possibly
  different buffer lengths), and the admission control capability is
  removed (by letting $\nu = \infty$),
 such policy reduces to the celebrated \emph{shortest queue
  routing} policy. The latter is known to be optimal 
under appropriate assumptions. See \cite{winston,horko,johri}.

In the case of constant service rates
$\mu_{k}(j_k) =
\mu_k$ and linear holding costs $h_{k}(j_k) = h_k \, j_k$,
under the long-run average criterion ($\alpha = 0$), identity
(\ref{eq:nukplhc})
yields the \emph{routing index}

\begin{equation}
\label{eq:ri}
\nu_k(j_k) = \frac{h_k}{\mu_k} \, \left[ 
\frac{\rho_k^{j_k+2}-1}{(\rho_k -1)^2}  - \frac{j_k+2}{\rho_k - 1}\right],
\end{equation}
where $\rho_{k} = \lambda/\mu_{k}$. 
The \emph{heavy
  traffic} case $\rho_k > 1$, where each queue  lacks
the capacity to process all the traffic, is of considerable interest
in applications; in such case, when there are 2 queues, the 
\emph{switching curve} in state space $(j_1, j_2)$ determined by such 
policy is \emph{asymptotically linear with limiting slope} $\ln \rho_1 / \ln \rho_2$ as
$j_1, j_2 \to \infty$.
The index policy above readily extends to models
with \emph{infinite buffers}.  

Note that the
standard heuristic in the linear cost case routes customers to the 
queue with smallest index $\widehat{\nu}_{k}(j_k) = h_k \, (j_k +1)/ \mu_{k}$.

\subsection{An index policy for scheduling a
multiclass make-to-stock queue with lost sales}
\label{s:smmtsq}
We next consider a model for scheduling a multiclass
make-to-stock queue (cf. \cite[Ch. 4]{buzshan}) in the \emph{lost sales} case, which
extends a simpler model studied by
Veatch and Wein in \cite{vewe} (having constant production and
demand rates, and linear holding costs).

A flexible production facility makes $m$
 products, labeled by $k = 1, \ldots, m$,  in a make-to-stock mode. 
The facility can work on at most an item at a time.
Finished product $k$ items are stored in
a dedicated stock, holding up to $n_{k}$ items. 
When this contains
$L_k(t) = j_k$ units, the facility can work at rate
$\mu_k(j_k)$ on such products, and corresponding \emph{customer
 orders} arrive at rate $\lambda_k(j_k)$. We assume mutually 
independent, exponential production and interarrival times. 
A product $k$'s order is immediately filled from stock 
if $j_k \geq 1$, and is otherwise
 lost. 
At each time, the facility can either stay idle, or 
engage in production of an item, by following a stationary policy. 

Product $k$ incurs state-dependent \emph{stock holding costs}, at rate
$c_k(j_k)$ per unit time;
\emph{stockout costs}, at rate $s_k$ per lost order; and is sold for
a state-dependent \emph{price} $r_k(j_k)$. 
The resulting product $k$'s \emph{net cost rate} per unit time in state $j_k$ is thus 
\[
h_k(j_k) = c_k(j_k) + s_k \, \lambda_k(0) \, 1\{j_k = 0\} -
r_k(j_k) \, \lambda_k(j_k) \, 1\{j_k > 0\}.
\]
We further assume that \emph{production is subsidized at rate 
$\nu$ per completed item}. 
Costs and rewards are discounted in time at rate
$\alpha > 0$.

We shall assume that model parameters satisfy the following
conditions (cf. Assumption \ref{ass:costconv}).
Let $d_k(j_k) = \lambda_k(j_k) - \mu_k(j_k)$ for $j_k \geq 1$. 
For consistency with previous analyses,
write $\Delta d_k(1) = \lambda_k(1) - \Delta \mu_k(1)$.

\begin{assumption}
\label{ass:mpmtsq}
For $1 \leq k \leq m$, the following holds: \\
{\rm (i)} Concave nondecreasing $d_k(j_k)$: 
$ 0 \leq \Delta d_k(j_k+1) \leq \Delta d_k(j_k), \quad 1 \leq j_k <
n_k$, 
 and   $\Delta d_k(1) > 0.$ \\
{\rm (iii)} Convex nondecreasing $h_k(j_k)$:
$0 \leq \Delta h_{k}(j_k) \leq \Delta h_{k}(j_k+1), \quad 0 \leq j_k < n_k.$
\end{assumption}

The goal is to design a state-dependent \emph{production scheduling policy}, 
which dynamically prescribes whether to engage in production and, if so,
of which product, 
so as to minimize
the expected total discounted value of costs accrued over an
infinite horizon.

The admission control index derived before readily
yields a heuristic index policy for such problem. 
The idea is to note that the present model is an RBP made up of $m$
single-queue admission control projects as studied
before, \emph{with the roles of parameters 
$\lambda$'s and $\mu$'s interchanged}.
Thus, opening queue $k$'s
\emph{entry gate} corresponds to making product $k$.
One must then, \emph{at each time, open at most one entry gate}. 
  
Let $\nu_{k}(j_k)$ be queue $k$'s admission control index, representing
the \emph{critical production subsidy} under which one should be indifferent between
idling and making product $k$ in state $j_k$.
Such interpretation leads to the following \emph{production control index policy}: 
\begin{enumerate}
\item Make a product $k$ with a nonfull stock level $j_k < n_k$ having the \emph{smallest
  index} $\nu_{k}(j_k)$ satisfying $\nu_k(j_k) < \nu$, if any is available.
\item Otherwise, idle the facility.
\end{enumerate}

Note that one may equivalently regard $-\nu$ as a \emph{production cost
  rate per completed item}. Hence, the indices $-\nu_k(j_k)$ represent 
\emph{critical production costs} for product $k$. 
Note further that, in the case of identical products, such policy
prescribes to make the product $k$ having the 
\emph{least stock} $j_k$ available, as long as $\nu_k(j_k) <
\nu$.

We next draw on the results in Section \ref{s:ccasr} to give 
explicit formulae for the index in some special cases, corresponding
to constant arrival and service rates
$\lambda_k(j_k) = \lambda_k$, $\mu_k(j_k) = \mu_k$, under the
\emph{long-run average} criterion $\alpha = 0$.
Let $\rho_k = \lambda_k/\mu_k \neq 1$.

Consider first the case of \emph{linear stock holding costs}  and
\emph{constant selling prices}, 
\[
h_k(j_k) = c_k \, j_k + s_k \, \lambda_k \, 1\{j_k = 0\} -
r_k \, \lambda_k \, 1\{j_k > 0\}.
\]
The results in Section \ref{s:ccasr}
then
yield the \emph{production index}

\begin{equation}
\label{eq:ri2}
\nu_k(j_k) = 
\frac{c_{k}}{\mu _{k}} \, \left[ \frac{\rho
    _{k}^{-j_k-1}-1}{(1-\rho _{k})^2 }- \frac{j_{k}+1}{1 - \rho_k} \right]
- r_k - s_k.
\end{equation}

\begin{remark} 
\label{rem:mtsind} 
\hspace{1in}
\begin{enumerate}
\item
The index (\ref{eq:ri2}) equals Whittle's  in 
\cite{vewe} \emph{scaled by  factor $1/\mu_k$}. Yet, although both indices
give the same (optimal) policy for a single-product problem, such
factor causes them to give
\emph{distinct} policies for the multi-product problem if the 
$\mu_k$'s differ.
\item The
index policy idles the facility when the number of units in stock for each product
lies at or above a corresponding critical \emph{base-stock} level. The
\emph{idling policy} is thus characterized by the
\emph{hedging-point} (cf. \cite{vewe}) consisting of such
base-stocks.
\item The index in (\ref{eq:ri2}) also gives a policy
  for a model with
\emph{unlimited storage capacity} ($n_k = \infty$). In
  such setting, if $\rho_k > 1$ for \emph{some} product $k$, then
  $\nu_k(j_k) < 0$. Hence, in the case $\nu = 0$, the facility will never
idle.
\end{enumerate}
\end{remark}

Consider next the case where stock holding costs are 
\emph{quadratic}, so that
\[
h_k(j_k) = c_k \, j_k^2 + s_k \, \lambda_k \, 1\{j_k = 0\} -
r_k \, \lambda_k \, 1\{j_k > 0\}.
\]
One then obtains, via (\ref{eq:nujqq}),  the production index

\begin{align}
\label{eq:nujq}
\nu_k(j_k) & = \frac{c_{k}}{\mu _{k}} \, 
\left[\left(\frac{2j_{k}+3}{(1-\rho_k)^2} -\frac{2}{(1-\rho _{k})^3}\right) \rho
  _{k}^{-j_k-1} \right.\\
& \quad \left. - \frac{(j_{k}+1)^{2}}{1-\rho_k} 
 -\frac{1}{(1-\rho_k)^2} +\frac{2}{(1-\rho _{k})^3}\right] - r_k - s_k. \nonumber
\end{align}

\section{Concluding remarks}
\label{s:c}
We have developed a polyhedral approach to the development of
dynamic allocation indices in a variety of stochastic scheduling
problems. 
In our view, such results offer a glimpse of the 
untapped potential 
which polyhedral methods have to offer in the field of stochastic 
optimization. 
We highlight two avenues for further research, which are the
subject of ongoing work: 
test empirically the proposed heuristic index policies, as in
\cite{benirb}; 
and provide approximate and asymptotic analyses of their performance, 
as in \cite{glanm01}. 

\appendix

\section{Discrete-time reformulation}
\label{s:acpdtr}
We reformulate the model of concern into
discrete time
by deploying the standard
\emph{uniformization}
technique (cf. \cite{lippman}), which
proceeds in two steps: (i) 
the original process $L(t)$ is reformulated into an
equivalent \emph{uniformized process} $\tilde{L}(t)$,
 having \emph{uniform
transition rate} $\Lambda$; process $\tilde{L}(t)$ is obtained by
sampling $L(t)$
 at time epochs corresponding to a Poisson process
with rate $\Lambda$; these includes \emph{real} as well as \emph{virtual} transitions,
in which no state change occurs; 
and (ii) process  $\tilde{L}(t)$ is  reformulated into
a \emph{discrete-time process} $X(t)$,
 by viewing inter-transition intervals as
discrete time periods.

Note that $\Lambda > 0$ is a valid \emph{uniform transition rate} iff it satisfies 
$$
\lambda_i + \mu_i \leq \Lambda, \quad i \in N.
$$
The resulting 
discrete-time process $X(t)$, for $t = 0, 1, \ldots$,
is  an RB (cf. Section \ref{s:srbosip})
 characterized by the following
elements: 
\begin{description}
\item[-State space: ] $N = \{0, 1, \ldots, n\}$; 
$N^{\{0, 1\}} = \{0, \ldots, n-1\}$; $N^{\{1\}} = \{n\}$.
\item[-Actions:] 
$a = 0$ (passive; open entry gate) and $a = 1$ (active; shut entry gate). 
\item[-Transition probability matrices:] Under action $a =
  1$, 
\begin{equation*}
\mathbf{P}^{1}  = \frac{1}{\Lambda} \, \left[
\begin{array}{ccccc}
\Lambda &  &  &  &    \\ 
\mu_1 & \Lambda - \mu_1 &  &  &    \\ 
& \ddots & \ddots &  &    \\ 
&  & \ddots & \ddots &    \\ 
  &  &  & \mu_{n} & \Lambda - \mu_{n}
\end{array}
\right];
\end{equation*}
and, under action $a = 0$, 

{\small
\begin{equation*}
\mathbf{P}^{0}  = \frac{1}{\Lambda} \, \left[
\begin{array}{ccccc}
\Lambda-\lambda_0 & \lambda_0 &  &   & \\ 
\mu_1 & \Lambda - \lambda_1 - \mu_1 & \lambda_1 &    &\\ 
& \ddots & \ddots & \ddots &  \\
 &  & \ddots & \ddots & \ddots   \\ 
&   &  & \mu_{n} & \Lambda - \mu_{n}
\end{array}
\right].
\end{equation*}}
\item[-One-period holding costs:]
$\mathbf{c}^0 = \mathbf{c}^1 = \displaystyle{
\frac{1}{\alpha + \Lambda}} \, \mathbf{h}$.
\item[-Discount factor:]
$\displaystyle{
\beta = \frac{\Lambda}{\alpha + \Lambda}.}
$
\end{description}

\section{Marginal workload and cost analysis}
\label{a:2}
\subsection{Marginal workloads: calculation and properties}
\label{s:cmw}
We next address the tasks of calculating 
marginal
workloads  $w^{S_k}_i$ for the admission control model,
and of establishing
their required properties.

\subsection*{Calculation of scaled $w^{S_k}_i$'s}
To avoid dependence on uniformization rate
$\Lambda$, 
the coefficients $w_i^S$ we shall calculate correspond to those
defined by (\ref{eq:wisdef}) after \emph{scaling} 
by factor $\alpha +
\Lambda$. 
Since 
\begin{equation}
\label{eq:p1mp0}
\mathbf{P}^1 - \mathbf{P}^0 = \frac{1}{\Lambda} \, \left[ 
\begin{array}{cccccc}
\lambda_0 & -\lambda_0 &  &  &  & \\ 
 & \lambda_1 & -\lambda_1 &  &  & \\ 
&  & \ddots & \ddots &  & \\ 
&  &  & \lambda_{n-1} & -\lambda_{n-1} & \\
&  &  &               & 0 & 0
\end{array}
\right],
\end{equation}
we have
\begin{equation}
\label{eq:wistisrel} 
w_i^S  = \begin{cases}
\displaystyle{\lambda_i \, \left[1
 - \Delta  b^S_{i+1}\right]} &  \text{if } 0 \leq i \leq n-1 \\
0 & \text{if } i = n.
\end{cases}
\end{equation}

Calculation of the 
$w_i^S$'s thus reduces to that of the $\Delta b^S_i$'s. 
To study the latter, we  start by characterizing the
coefficients $b^{S_k}_i$, through their defining equations in
(\ref{eq:tiscalc}).
We shall denote by $\lambda_i^S$ the  \emph{birth rate}
 in state $i$ under the
$S$-active policy, i.e., 
$$
\lambda_i^S = \lambda_i \, 1\{i \in N^{\{0, 1\}} \setminus S\}, \quad i \in N.
$$
Note that $\lambda^{S_{k}}_i = \lambda_i \, 1\{0 \leq i < k-1\}$, for 
$1 \leq k \leq n+1$.

\begin{lemma}
\label{lma:1actp} For  $1 \leq k \leq n+1$, coefficients
$b^{S_{k}}_i$ are characterized
by the  equations
\begin{equation*}
\begin{split}
(\alpha + \Lambda) \, b^{S_{k}}_0 & = 
 \lambda_0 - \lambda^{S_{k}}_0
 + 
 (\Lambda - \lambda^{S_{k}}_0) \, b^{S_{k}}_0 + 
 \lambda^{S_{k}}_0 \, b^{S_{k}}_1 \\
(\alpha + \Lambda) \, b^{S_{k}}_i & = 
 \lambda_i - \lambda^{S_{k}}_i
 + 
 \mu_i \,  b^{S_{k}}_{i-1} 
 + 
 (\Lambda - \lambda^{S_{k}}_i - \mu_i) \, b^{S_{k}}_i + 
 \lambda^{S_{k}}_i \, b^{S_{k}}_{i+1}, \,  
 1 \leq i \leq n-1 \\
(\alpha + \Lambda) \, b^{S_{k}}_n & = \lambda_n + \mu_n \,
  b^{S_{k}}_{n-1} + 
 (\Lambda  - \mu_n) \, b^{S_{k}}_n.
\end{split}
\end{equation*}
\end{lemma}

The next result, characterizing coefficients
$\Delta b^{S_k}_i$, follows immediately.

\begin{lemma}
\label{lma:1actdeltap}
For $1 \leq k \leq n+1$,
coefficients $\Delta b^{S_{k}}_i$ are characterized by the 
equations
\begin{equation*}
\begin{split}
(\alpha + \lambda_0^{S_k} + \mu_1) \, 
  \Delta b^{S_{k}}_1 & = 
\Delta \lambda_1 - \Delta \lambda^{S_{k}}_1  + 
  \lambda^{S_{k}}_1 \, \Delta b^{S_{k}}_2 \\
(\alpha + \lambda^{S_{k}}_{i-1} + \mu_i) \, 
  \Delta b^{S_{k}}_i & = \Delta \lambda_i -
  \Delta \lambda^{S_{k}}_i
 + 
  \mu_{i-1} \, \Delta b^{S_{k}}_{i-1}  + 
  \lambda^{S_{k}}_i \, \Delta b^{S_{k}}_{i+1}, \, 
 2 \leq i \leq n-1 \\
(\alpha  + \lambda^{S_{k}}_{n-1} + \mu_n) \, 
  \Delta b^{S_{k}}_n & = \Delta \lambda_{n} + \lambda^{S_{k}}_{n-1} +
  \mu_{n-1} \, \Delta b^{S_{k}}_{n-1}.
\end{split}
\end{equation*}
\end{lemma}

We next develop a recursive procedure to solve
the equations in Lemma \ref{lma:1actdeltap},
based on the following observations: (i) 
the equations give
$$ \Delta b^{S_1}_1 = \frac{\Delta \lambda_1}{\alpha + \mu_1},$$ 
from which remaining $\Delta b^{S_1}_i$'s are 
calculated;
(ii) for $1 \leq k \leq n$, 
once \emph{pivot coefficient} $\Delta b^{S_{k+1}}_k$ is
available, 
they give the remaining $\Delta b^{S_{k+1}}_i$'s; and (iii)
the first  pivot is
$$\Delta b^{S_{2}}_1 = \frac{\lambda_1}{\alpha + \lambda_0 +
  \mu_1}.$$ 
Hence, if we can express
pivot 
 $\Delta b^{S_{k+2}}_{k+1}$ in terms
 of  $\Delta b^{S_{k+1}}_k$, for $1 \leq k \leq n-1$, this would complete
 a recursion to calculate all coefficients 
$\Delta b^{S_k}_i$.

We next seek to relate successive pivots,  
drawing on \cite{chenyao}.
Consider, for  $1 \leq k \leq n-1$, the vectors
(where $\mathbf{x}^{T}$ denotes the transpose of vector 
$\mathbf{x}$) 
\begin{equation*}
\begin{split}
\Delta \mathbf{b}^k & = \left(  \Delta b^{S_{k+1}}_1, \ldots,
\Delta b^{S_{k+1}}_k\right)^{T} \\
\Delta \widehat{\mathbf{b}}^k & = \left(\Delta b^{S_{k+2}}_1,
\ldots,  \Delta b^{S_{k+2}}_k\right)^{T} \\
\mathbf{b}^k  & =  
\frac{\lambda_k}{\alpha + \lambda_{k-1} +
\mu_k} \, \mathbf{e}_k \\
\widehat{\mathbf{b}}^k & = 
\frac{\lambda_k \, \Delta b^{S_{k+2}}_{k+1}}{\alpha + \lambda_{k-1}
  + \mu_k} \, \mathbf{e}_k, 
\end{split}
\end{equation*}
where $\mathbf{e}_k$ is the $k$th unit coordinate vector in
 $\mathbb{R}^k$.
Let further $\mathbf{B}^k$ be  the $k \times k$ matrix 
\begin{equation*}
\mathbf{B}^k = \left[ 
\begin{array}{ccccc}
0 & \frac{\lambda_1}{\alpha + \lambda_0 + \mu_1} &  &  &    \\ 
\frac{\mu_1}{\alpha + \lambda_1 + \mu_2} & 0 & \frac{\lambda_2}{\alpha +
\lambda_1 + \mu_2} &  &    \\ 
& \ddots & \ddots & \ddots &    \\ 
&  & \ddots & \ddots & \ddots   \\
  &  &  & \frac{\mu_{k-1}}{\alpha + \lambda_{k-1} + \mu_k} & 0
\end{array}
\right],
\end{equation*}
with
$\mathbf{B}^1 = 0$. 
The next result reformulates  some  equations
in Lemma \ref{lma:1actdeltap}.

\begin{lemma}
\label{lma:deltatvecrec} For 
$1 \leq k \leq n-1$: 
\begin{itemize}
\item[(a)] $\Delta \mathbf{b}^k = \mathbf{b}^k + \mathbf{B}^k \, \Delta 
\mathbf{b}^k $.
\item[(b)] $\Delta \widehat{\mathbf{b}}^k = \widehat{\mathbf{b}}^k + 
\mathbf{B}^k \, \Delta \widehat{\mathbf{b}}^k $.
\end{itemize}
\end{lemma}

Recall that
$d_i = \mu_i - \lambda_i$, and that 
we require that the $d_i$ satisfy Assumption \ref{ass:costconv}.
To proceed, introduce
coefficients 
\begin{equation}
\label{eq:akdef}
a_k = 
\begin{cases}
1 & \text{if $k = 1$} \\ \\
\displaystyle{\frac{\det \left(\mathbf{I} - \mathbf{B}^{k}\right)}{\det \left(%
\mathbf{I} - \mathbf{B}^{k-1}\right)}} &  \text{if $2 \leq k \leq n$.}
\end{cases}
\end{equation}

\begin{lemma}
\label{lma:gicalc} 
Under
Assumption \ref{ass:costconv}(i), the following holds: \\
\begin{itemize}
\item[(a)]  $a_k > 0$, for $1 \leq k \leq n$. \\
\item[(b)] The $a_k$'s can be computed recursively by letting $a_1 = 1$ and 
\begin{equation*}
a_k =  1 - \frac{\lambda_{k-1} \, \mu_{k-1}}{(\alpha + \lambda_{k-2} +
\mu_{k-1}) \, (\alpha + \lambda_{k-1} + \mu_k)} \, \frac{1}{a_{k-1}}, \quad
2 \leq k \leq n.
\end{equation*}
\item[(c)]
$
\frac{\displaystyle \alpha + \mu_k}{\displaystyle \alpha + \lambda_{k-1} + \mu_k} < a_k
< 1, \quad 2 \leq k \leq n.
$
\end{itemize}
\end{lemma}
\begin{proof}
(a)
Under Assumption \ref{ass:costconv}(i) 
the row sums of $\mathbf{B}^{k}$ are less than unity,
and hence so is its spectral radius.
It follows that 
$\det \left(\mathbf{I} - \mathbf{B}^{k}\right) > 0$, 
which proves the result.

(b)
The recursion follows from the definition of $a_k$ and the
identity 
\begin{align*}
\det (\mathbf{I} - \mathbf{B}^{k}) & = \det (\mathbf{I} - 
\mathbf{B}^{k-1})
- \frac{\lambda_{k-1} \, \mu_{k-1}}{ (\alpha +
\lambda_{k-2} + \mu_{k-1}) \, (\alpha + \lambda_{k-1} + \mu_k)} \, \det
(\mathbf{I} - \mathbf{B}^{k-2}) 
\end{align*}

(c)
Let $2 \leq k \leq n$. It follows from (a) and (b) that 
$a_k < 1$. 
We  next show that
$$ a_k > \frac{\alpha + \mu_k}{\alpha + \lambda_{k-1} +
  \mu_k}, \quad 1 \leq k \leq n,$$
by induction on $k$. The case $k = 1$ is trivial.
Assume the result holds for $k - 1$, i.e., 
$$a_{k-1} > \frac{\alpha + \mu_{k-1}}{\alpha + \lambda_{k-2} + \mu_{k-1}}.$$
Then, part (b) and the induction hypothesis yield 
\begin{equation*}
\begin{split}
a_k &= 1 - \frac{\lambda_{k-1}}{\alpha + \lambda_{k-1} + \mu_k} \, 
     \frac{\displaystyle{\frac{\mu_{k-1}}{\alpha + \lambda_{k-2} +
     \mu_{k-1}}}}{a_{k-1}} \\
    &> 1 - \frac{\lambda_{k-1}}{\alpha + \lambda_{k-1} + \mu_k} 
     = \frac{\alpha + \mu_k}{\alpha + \lambda_{k-1} + \mu_k},
\end{split}
\end{equation*}
which completes the  proof.
\end{proof}
\qed

We are now ready to relate
successive pivots.

\begin{lemma}
\label{lma:keyreldeltat} For $1 \leq k \leq n-1$, 
\begin{equation*}
a_{k+1} \, \left[1 - 
\Delta b^{S_{k+2}}_{k+1}\right] = 
\frac{\alpha + \Delta d_{k+1} +
\mu_k \, \left[1 - \Delta b_k^{S_{k+1}}\right]}{\alpha + \lambda_k + \mu_{k+1}}; 
\end{equation*}
or, equivalently, 
\begin{equation*}
  w^{S_{k+2}}_k = \frac{\lambda_k}{a_{k+1}}
\frac{\displaystyle \alpha + \Delta d_{k+1} +  
  \frac{w^{S_{k+1}}_{k-1}}{\rho_{k-1}}}{\alpha + \lambda_k + \mu_{k+1}}.
\end{equation*}
\end{lemma}
\begin{proof}
Fix $1 \leq k \leq n-1$.
By Lemma \ref{lma:deltatvecrec} and the definitions of 
$\mathbf{h}^k$, $\widehat{\mathbf{h}}^k$, we have 
\begin{equation} 
\begin{split}
\Delta \mathbf{b}^k - \Delta \widehat{\mathbf{b}}^k & = (\mathbf{I} - 
\mathbf{B}^k)^{-1} \, (\mathbf{b}^k - \widehat{\mathbf{b}}%
^k)  \label{eq:diftihatti} \\
& = \frac{ 
\lambda_k \, \left[1 -   \Delta b^{S_{k+2}}_{k+1}\right]}{
\alpha + \lambda_{k-1} + \mu_k} \, 
(\mathbf{I} - \mathbf{B}^k)^{-1} \, \mathbf{e}^k. 
\end{split}
\end{equation}
Now, noting that the element in position $(k, k)$ of matrix 
$\left(\mathbf{I}- \mathbf{B}^k\right)^{-1}$ is 
$$\frac{\det\left(\mathbf{I} - \mathbf{B}^{k-1}\right)}{ 
\det \left(\mathbf{I} - \mathbf{B}^k\right)}, $$
which by
definition equals $1/a_k$, it follows from 
the last identity above that

\begin{equation}
\label{eq:dtisitisip1}
\Delta b^{S_{k+1}}_k - \Delta b^{S_{k+2}}_k = 
\frac{1}{a_k} \, \frac{ \lambda_k \, (1 -  
\Delta b^{S_{k+2}}_{k+1})}{
\alpha + \lambda_{k-1} + \mu_k}. 
\end{equation}
Combining the previous identity with 
\begin{equation*}
\Delta b^{S_{k+2}}_{k+1} = \frac{\lambda_{k+1}}{\alpha +
\lambda_k + \mu_{k+1}} + \frac{\mu_k}{\alpha + \lambda_k +
\mu_{k+1}} \,  \Delta b^{S_{k+2}}_k, 
\end{equation*}
(cf. Lemma \ref{lma:1actdeltap}), and substituting  for $a_k$ in terms of $a_{k+1}$ 
(cf. Lemma \ref{lma:gicalc}), yields
the required identities (after straightforward algebra).
\end{proof}
\qed

\subsection*{Properties of marginal workloads}
We next set out to establish 
properties of
marginal workloads which are invoked in Section \ref{s:qaca}.

\begin{lemma}
\label{lma:aiigt0} 
Under Assumption \ref{ass:costconv}(i), 
the following holds, for $\alpha \geq 0$: \\
\begin{itemize}
\item[(a)] $w^{S_{k+1}}_{k-1} > 0, \quad 1 \leq k \leq n$. \newline
\item[(b)] $w^{S_{k+2}}_{i-1} > w^{S_{k+1}}_{i-1}, \quad 1
\leq i \leq k \leq n-1$. \newline
\item[(c)] $w^{S_{k+1}}_{i-1} > 0 \Longrightarrow  w^{S_{k+1}}_i
 > 0, \quad  1 \leq k \leq i \leq n-1$.  \newline
\item[(d)]
$w^{S_{k+1}}_i > w^{S_{k+2}}_i, \quad 1 \leq k \leq i \leq n-1.$ 
\end{itemize}
\end{lemma}
\begin{proof}
(a) Proceed by  induction on $k$.
 The case $k =1$ holds by the expression for 
$w^{S_2}_0$ in Figure \ref{fig:aibcomp2} and Assumption \ref{ass:costconv}(i).
Suppose now $w^{S_k}_{k-2} > 0$. We  have 
\begin{equation*}
\begin{split}
\frac{w^{S_{k+1}}_{k-1}}{\lambda_{k-1}} & = \frac{1}{a_k} \,
\frac{\displaystyle \alpha + \Delta d_k + 
  \frac{w^{S_k}_{k-2}}{\rho_{k-2}}}{\alpha + \lambda_{k-1} + \mu_k}
> 0,
\end{split}
\end{equation*}
where the identity is taken from Figure \ref{fig:aibcomp2}, and the
inequality follows from the induction hypothesis, along with
Assumption \ref{ass:costconv}(i) and $a_k > 0$ (Lemma \ref{lma:gicalc}).

(b) Using (\ref{eq:wistisrel}),  we
 can rewrite identity (\ref{eq:diftihatti}) as 
\begin{equation*}
\Delta \mathbf{b}^k - \Delta \widehat{\mathbf{b}}^k = 
\frac{w^{S_{k+2}}_k}{\alpha + \lambda_{k-1} + \mu_k} \,
 \, 
(\mathbf{I} - \mathbf{B}^k)^{-1} \, \mathbf{e}_k. 
\end{equation*}
Now,  since the spectral radius of $\mathbf{B}^k$ is
less than unity (cf. Lemma \ref{lma:gicalc}'s proof),
matrix $\left(\mathbf{I} - \mathbf{B}^k\right)^{-1}$ is positive componentwise, and 
hence $\left(\mathbf{I} - \mathbf{B}^k\right)^{-1} \, \mathbf{e}_k > \mathbf{0}$.
Combining this with part (b) and the last identity above yields 
$\Delta \mathbf{b}^k - \Delta \widehat{\mathbf{b}}^k > \mathbf{0}$, i.e., 
\begin{equation*}
\Delta b^{S_{k+1}}_i > \Delta b^{S_{k+2}}_i,  \quad
1 \leq i \leq k.
\end{equation*}
By (\ref{eq:wistisrel}), these inequalities give
the required result.

(c) The result follows from
\begin{equation*}
\frac{w^{S_{k+1}}_i}{\lambda_i} =
\frac{\displaystyle \alpha + \Delta d_{i+1} + 
\frac{w^{S_{k+1}}_{i-1}}{\rho_{i-1}}}{\alpha + \mu_{i+1}}, 
\quad k \leq i \leq n-1
\end{equation*}
(cf. Figure \ref{fig:aibcomp2}), and  Assumption \ref{ass:costconv}(i).

(d)
By (\ref{eq:wistisrel}), the result is equivalent to
\begin{equation}
\label{eq:equivreltw}
\Delta b^{S_{k+1}}_{i+1} < \Delta b^{S_{k+2}}_{i+1}, \quad 1 \leq k \leq
i \leq n-1.
\end{equation}

Now,  it follows from Lemma \ref{lma:1actdeltap}
that, for $1 \leq k \leq n-1$, 
\begin{equation*}
\begin{split}
(\alpha + \mu_i) \, \Delta b^{S_{k+1}}_i  & =
   \mu_{i-1} \, \Delta b^{S_{k+1}}_{i-1}, 
\quad k+1 \leq i \leq n-1 \\
(\alpha + \mu_i) \, \Delta b^{S_{k+2}}_i  & =
   \mu_{i-1} \, \Delta b^{S_{k+2}}_{i-1}, 
\quad k+2 \leq i \leq n-1, 
\end{split}
\end{equation*}
 hence
$$
\Delta b^{S_{k+2}}_i - \Delta b^{S_{k+1}}_i =
   \frac{\mu_{i-1}}{\alpha + \mu_i} \, 
 (\Delta b^{S_{k+2}}_{i-1} - \Delta b^{S_{k+}}_{i-1}) , 
\quad k+2 \leq i < n.
$$
In light of the last identity, to prove (\ref{eq:equivreltw})
it is enough to show that
$$
\Delta b^{S_{k+2}}_{k+1} - \Delta b^{S_{k+1}}_{k+1} > 0,
$$
which we establish next. 
Consider the case $k = 0$. By Lemma \ref{lma:1actdeltap}, we have
\begin{equation*}
\begin{split}
\Delta b^{S_{2}}_1 - \Delta b^{S_1}_1 & = 
\frac{\lambda_1}{\alpha + \lambda_0 + \mu_1} - 
\frac{\Delta \lambda_1}{\alpha + \mu_1} \\
& = \frac{\lambda_0}{\alpha + \mu_1} \, \frac{\alpha + \Delta d_1}{
\alpha + \lambda_0 + \mu_1}  >
0,
\end{split}
\end{equation*}
where the  inequality follows by Assumption \ref{ass:costconv}(i).
Consider now the case $k \geq 1$. 
Drawing again on Lemma \ref{lma:1actdeltap}, we have
\begin{equation*} 
\begin{split}
(\alpha + \lambda_k + \mu_{k+1}) \, \Delta b^{S_{k+2}}_{k+1} & = 
\lambda_{k+1} + \mu_k \, \Delta b^{S_{k+2}}_k \\
(\alpha + \mu_{k+1}) \, \Delta b^{S_{k+1}}_{k+1} & = 
\Delta \lambda_{k+1}
  + 
  \mu_k \, \Delta b^{S_{k+1}}_k.
\end{split}
\end{equation*}
Using in turn the last two identities, 
(\ref{eq:dtisitisip1}) and (\ref{eq:wistisrel}), 
 part (a) 
and Lemma \ref{lma:gicalc}(c), yields
\begin{equation*}
\begin{split}
(\alpha + \mu_{k+1}) \, 
(\Delta b^{S_{k+2}}_{k+1} - \Delta b^{S_{k+1}}_{k+1}) 
& = 
\lambda_k \, \left[1 -  \Delta b^{S_{k+2}}_{k+1}\right] + 
   \mu_k \, (\Delta b^{S_{k+2}}_k - \Delta b^{S_{k+2}}_k) \\
& = 
\left[1 - \frac{\mu_k / a_k}{\alpha +\lambda_{k-1} + \mu_k}\right]
 \,  w^{S_{k+2}}_k > 0,
\end{split}
\end{equation*}
as required. This completes the proof.
\end{proof}
\qed

\subsection{Marginal cost calculation}
\label{s:mccc}
We set out in this section to calculate 
marginal
costs  $c_i^{S_k}$, proceeding similarly as before for
marginal workloads. Again, to eliminate
the dependence on uniformization rate $\Lambda$, the terms 
$c_i^S$ below correspond to those defined by
(\ref{eq:jaiscdef}) after \emph{scaling} by factor $\alpha + \Lambda$.

We start by relating coefficients $v_i^S$'s and $c_i^S$'s. From
 (\ref{eq:jaiscdef}) and (\ref{eq:p1mp0}),
we obtain
\begin{equation*}
%\label{eq:cistisrel} 
c_i^S = 
\lambda_i \, \Delta
  v^S_{i+1}, \quad
0 \leq i \leq n-1.
\end{equation*}
We must thus calculate the
$\Delta v^{S_{k}}_i$'s.
Start by calculating the $v_i^{S_k}$'s through 
(\ref{eq:jmn1}).

\begin{lemma}
\label{lma:2actp} 
For $1 \leq k \leq n+1$, coefficients $v^{S_{k}}_i$ are characterized
by the equations
\begin{equation*}
\begin{split}
(\alpha + \Lambda) \, v^{S_{k}}_0 & = h_0  + 
 (\Lambda - \lambda^{S_{k}}_0) \, v^{S_{k}}_0 + 
 \lambda_0^{S_{k}} \, v^{S_{k}}_1 \\
(\alpha + \Lambda) \, v^{S_{k}}_i & = h_i + 
 \mu_i \,  v^{S_{k}}_{i-1} + 
 (\Lambda - \lambda^{S_{k}}_i - \mu_i) \, v^{S_{k}}_i  + 
 \lambda^{S_{k}}_i \, v_{i+1}^{S_k}, \quad 
 1 \leq i \leq n-1 \\
(\alpha + \Lambda) \, v^{S_{k}}_n & = h_n + \mu_{n} \,  v^{S_{k}}_{n-1} + 
 (\Lambda  - \mu_n) \, v^{S_{k}}_n.
\end{split}
\end{equation*}
\end{lemma}

It follows that coefficients 
$\Delta v_i^{S_k}$ are characterized as shown next.
\begin{lemma}
\label{lma:2actdeltap}
For $1 \leq k \leq n+1$,
coefficients $\Delta v^{S_{k}}_i$ are characterized by the 
equations
\begin{equation*}
\begin{split}
(\alpha + \lambda_0^{S_k} + \mu_1) \, 
  \Delta v_1^{S_k} & = \Delta h_1   + 
  \lambda^{S_{k}}_1 \, \Delta v^{S_{k}}_2 \\
(\alpha + \lambda^{S_{k}}_{i-1} + \mu_i) \, 
  \Delta v^{S_{k}}_i & = \Delta h_i  + 
  \mu_{i-1} \, \Delta v^{S_{k}}_{i-1} + 
  \lambda^{S_{k}}_i \, \Delta v_{i+1}^{S_k}, \quad 
 2 \leq i \leq n-1 \\
(\alpha  + \lambda^{S_{k}}_{n-1} + \mu_n) \, 
  \Delta v^{S_{k}}_n & = \Delta h_n + 
  \mu_{n-1} \, \Delta v^{S_{k}}_{n-1}.
\end{split}
\end{equation*}
\end{lemma}

We next develop  a recursion to calculate  
\emph{pivot} terms $\Delta v^{S_{k+1}}_k$, 
along the lines followed in Appendix \ref{s:cmw} to calculate the
$\Delta b^{S_{k+1}}_k$'s.
Note that Lemma \ref{lma:2actdeltap} yields 
\[
\Delta v^{S_1}_1 = \frac{\Delta h_1}{\alpha + \mu_1},
\]
and hence
\[
c^{S_1}_0 = \frac{\lambda_0}{\alpha + \mu_1} \, \Delta h_1.
\]
It further yields 
the first such pivot as
$$
\Delta v^{S_2}_1 = \frac{\Delta h_1}{\alpha + \lambda_0 + \mu_1},
$$
so that
\begin{equation*}
c^{S_2}_0  = 
  \frac{\lambda_0}{\alpha + \lambda_0 + \mu_1} \, \Delta h_1.
\end{equation*}

We next set out to relate successive pivots.
Associate with  $1 \leq k \leq n-1$ the vectors
\begin{equation*}
\begin{split}
\Delta \mathbf{v}^k & = \left(  \Delta v^{S_{k+1}}_1, \ldots,
\Delta v^{S_{k+1}}_k\right)^{T} \\
\Delta \widehat{\mathbf{v}}^k & = \left(\Delta v^{S_{k+2}}_1,
\ldots,  \Delta v^{S_{k+2}}_k\right)^{T} \\
\mathbf{h}^k  & =  
\left(\frac{\Delta h_1}{\alpha + \lambda_0 + \mu_1}, 
   \ldots, \frac{\Delta h_k}{\alpha +\lambda_{k-1} + \mu_k}\right) \\
\widehat{\mathbf{h}}^k & = \mathbf{h}^k + 
\frac{\lambda_k \, \Delta v^{S_{k+2}}_{k+1}}{\alpha +\lambda_{k-1} + \mu_k} \, \mathbf{e}_k. 
\end{split}
\end{equation*}

The next result is a counterpart to Lemma \ref{lma:deltatvecrec}.

\begin{lemma}
\label{lma:2deltatvecrec} For 
$1 \leq k \leq n-1$, \newline
\begin{itemize}
\item[(a)] $\Delta \mathbf{v}^k = \mathbf{h}^k + \mathbf{B}^k \, \Delta 
\mathbf{v}^k $; \newline
\item[(b)] $\Delta \widehat{\mathbf{v}}^k = \widehat{\mathbf{h}}^k + 
\mathbf{B}^k \, \Delta \widehat{\mathbf{v}}^k$.
\end{itemize}
\end{lemma}

The  relation between pivots
 $\Delta v^{S_{k+1}}_k$ and $\Delta v^{S_{k+2}}_{k+1}$, and its 
marginal cost reformulation is given next. 
The proof is similar to that of Lemma
 \ref{lma:keyreldeltat}, and is hence omitted.

\begin{lemma}
\label{lma:2keyreldeltat} For $1 \leq k \leq n-1$, 
\begin{equation*}
a_{k+1} \, \Delta v^{S_{k+2}}_{k+1} 
 = 
\frac{
\Delta h_{k+1}}{\alpha + \lambda_k + \mu_{k+1}} + 
\frac{\mu_k}{\alpha + \lambda_k + \mu_{k+1}} \,
 \Delta v^{S_{k+1}}_k;
\end{equation*}
or, equivalently, 
\begin{equation*} 
c^{S_{k+2}}_k
 = \frac{\lambda_k}{a_{k+1}} \, \frac{\displaystyle \Delta h_{k+1} + 
  \frac{c^{S_{k+1}}_{k-1}}{\rho_{k-1}}}{\alpha + \lambda_k + \mu_{k+1}}.
\end{equation*}
\end{lemma}

\section{Possible inconsistency of the Whittle index relative to 
threshold policies}
\label{a:1}
The reader may wonder whether the extra flexibility provided by
parameters $\theta_j^1$ in the new index introduced in
Definition \ref{def:ewi}
significantly expands the scope of the original Whittle index. 
We argue next that such is the case by showing, in
the
setting of the admission control model,  that
the Whittle index does \emph{not} rank the states in a manner
consistent with threshold policies, under the
parameter range 
given by Assumption \ref{ass:costconv}.

Recall that the Whittle index arises from the appropriate
$\nu$-charge problem  obtained by charging costs at 
rate $\nu$ \emph{while the entry gate is shut}. Namely, the
corresponding activity
measure $b^u$ obtains by letting $\theta_j^1 = 1$, for 
$j \in N = \{0, \ldots, n\}$

We shall consider that an index policy for the
admission control model
is \emph{consistent with threshold policies} if index $\nu_j$ 
is \emph{nondecreasing}  on
$j \in \{0, \ldots, n-1\}$.

Consider  the case where the buffer size is $n = 2$,   
service rates are
$\mu_j = \mu$, 
and cost rates are 
$h_j =  h \, j$. 
Suppose
 arrival rates $\lambda_j$ 
are strictly decreasing on $j$, namely
\begin{equation}
\label{eq:lamb2nicv}
\Delta \lambda_2 < 0, \Delta \lambda_1 < 0.
\end{equation}
It then follows that Assumption \ref{ass:costconv}  holds.

Take, in particular,  
$\lambda_0=1$,
$\lambda_1=\frac{1}{2}$, $\lambda_2 = \frac{1}{4}$, 
$\mu =\frac{3}{2}$, $\alpha = \frac{1}{33}$, $h = 1$. 
Pick the uniformization rate
$\Lambda =3$, so that
$\beta = \Lambda/(\alpha + \Lambda) = \frac{99}{100}$.

The corresponding RB is
indexable, in Whittle's sense, and has Whittle indices
$$
\nu_2 = 0 < 
\nu_1 = \frac{3300}{6767}    < 
\nu_0 = \frac{11\,022}{19\,111}.
$$
They thus give a state ranking which is \emph{inconsistent} with 
threshold policies.

Such inconsistency only arises, however,
under state-dependent arrival rates; 
under a constant
arrival rate $\lambda$, the extended index equals Whittle's scaled
by factor $1/\lambda$.

\subsubsection*{Acknowledgement} 
The author has presented preliminary versions of this paper
at the 
Conference on Stochastic Networks (Madison, WI, June 2000), the 27th
Conference on Stochastic Processes and Their Applications (Cambridge,
UK, July 2001), the 11th INFORMS Applied Probability Conference
(New York, July 2001), the Dagstuhl 
Seminar on Scheduling in
Computer and Manufacturing Systems (Wadern, Germany, June 2002), and
a Eurandom seminar (Eindhoven, Netherlands, February 2001).
The author wishes to thank Professors K. Glazebrook, G. Weiss, and
P. Whittle for 
interesting discussions, and the
two referees for their suggestions, which helped improve the 
paper's exposition.

\end{document}